\documentclass[a4paper]{amsart}
\usepackage{amsmath,amssymb,amsthm}
\usepackage{mathrsfs}
\usepackage{aliascnt}
\usepackage{cleveref}
\usepackage{setspace}
\usepackage{needspace}
\usepackage{url}
\allowdisplaybreaks
\usepackage[textwidth=156mm,textheight=240mm,centering,
  headheight=12pt,headsep=18pt,footskip=24pt]{geometry}
\newcommand{\newaliasthm}[3]{%
  \newaliascnt{#1}{theorem}%
  \newtheorem{#1}[#1]{#2}%
  \aliascntresetthe{#1}%
  \crefname{#1}{#2}{#3}%
  \Crefname{#1}{#2}{#3}%
}
\theoremstyle{plain}
\newtheorem{theorem}{Theorem}[section]
\crefname{theorem}{Theorem}{Theorems}
\Crefname{theorem}{Theorem}{Theorems}
\newaliasthm{proposition}{Proposition}{Propositions}
\newaliasthm{lemma}{Lemma}{Lemmas}
\newtheorem*{lemma*}{Lemma}
\newaliasthm{corollary}{Corollary}{Corollaries}
\theoremstyle{definition}
\newaliasthm{definition}{Definition}{Definitions}
\newaliasthm{problem}{Problem}{Problems}
\newtheorem*{bannaiproblem}{Problem (Bannai)}
\newtheorem*{cxxproblem}{Problem (Cui--Xia--Xiang)}
\newtheorem*{xiangconjecture}{Conjecture (Xiang)}
\newaliasthm{example}{Example}{Examples}
\newtheorem*{example*}{Example}
\theoremstyle{remark}
\newaliasthm{remark}{Remark}{Remarks}
\crefname{section}{Section}{Sections}
\Crefname{section}{Section}{Sections}
\crefname{subsection}{Section}{Sections}
\Crefname{subsection}{Section}{Sections}
\crefname{subsubsection}{Section}{Sections}
\Crefname{subsubsection}{Section}{Sections}
\crefname{appendix}{Section}{Sections}
\Crefname{appendix}{Section}{Sections}
\crefname{equation}{equation}{equations}
\Crefname{equation}{Equation}{Equations}
\numberwithin{equation}{section}
\newcommand{\R}{\mathbb R}
\newcommand{\Z}{\mathbb Z}
\newcommand{\N}{\mathbb N}
\newcommand{\C}{\mathbb C}
\newcommand{\Q}{\mathbb Q}
\newcommand{\TQ}{\mathbb T_{\Q}}
\ExplSyntaxOn
\NewDocumentCommand{\compactindexlist}{mm}
  {
    \par\smallskip
    \begingroup
    \fontsize{6}{7}\selectfont
    \setlength{\parindent}{0pt}
    \raggedright
    \noindent
    \(#1=\{\clist_use:nn {#2} {,\allowbreak\mkern2mu}\}\)\par
    \endgroup
    \smallskip
  }
\ExplSyntaxOff
\newcommand{\Hom}{\operatorname{Hom}}
\renewcommand{\d}{\mathrm d}
\renewcommand{\Re}{\operatorname{Re}}
\renewcommand{\Im}{\operatorname{Im}}

\title{Rational Designs and Prouhet--Tarry--Escott Problems}

\author[R.~Misawa]{Ryutaro Misawa}
\address[R.~Misawa]{Division of Mathematics\\
Graduate School of Information Sciences\\
Tohoku University\\
6-3-09 Aramaki-Aza-Aoba, Aoba-ku, Sendai 980-8579\\
Japan}
\email{misawa.ryutaro.q2@dc.tohoku.ac.jp}

\author[M.~Sawa]{Masanori Sawa}
\address[M.~Sawa]{Department of System Informatics\\
Graduate School of System Informatics\\
Kobe University\\
1-1 Rokkodai-cho, Nada-ku, Kobe, Hyogo 657-8501\\
Japan}
\email{sawa@people.kobe-u.ac.jp}

\author[F.~\v{S}piranec]{Filip \v{S}piranec}
\address[F.~\v{S}piranec]{Independent Researcher}
\email{spiranec.filip@gmail.com}
\date{}

\begin{document}

\begin{abstract}
For every strength, we prove the existence of rational designs on spheres, orthogonal and unitary groups, real and complex Grassmannians, and ellipsoids that are defined by positive-definite quadratic forms with rational coefficients and have a rational point.
This answers Bannai's Problem for these spaces, solves the problem of Cui--Xia--Xiang on rational spherical $4$-designs on the unit circle, and proves Xiang's conjecture in the stronger form of existence in every dimension.
To this end, we introduce the $(G,\mathcal S)$-Prouhet--Tarry--Escott problem for an abelian group $G$ and a finite set of characters $\mathcal S$ and develop its existence and minimality theory.
We also prove a general theorem transforming integer-weighted rational designs into
equi-weighted designs
on distinct rational points.
On the unit circle $\mathbb{S}^1$, we determine the minimum sizes of rational spherical designs of strengths $4,5,6,7$ and show that they are respectively $2304,4608,9216,9216$.
We also show that the minimum size of a rational spherical $5$-design on $\mathbb S^2$ is $192$.
We further apply our framework to prove a conjecture of Thuswaldner--Ziegler on sums of distinct units and to answer Mahabaduge's problem on convex equilateral lattice polygons affirmatively.
\end{abstract}

\maketitle

\section{Introduction}\label{sec:introduction}

Let $Z\subseteq\mathbb{R}^d$ be a nonempty set equipped with a probability measure $\mu$, and let $\mathcal{V}$ be a finite-dimensional real vector space of $\mu$-integrable real-valued functions on $Z$ containing the constant function $\mathbf{1}$.
A nonempty finite subset $\mathfrak{X} \subseteq Z$ is called a \emph{$\mathcal V$-design} if
\begin{equation}\label{eq:equal-weight-design}
\frac{1}{|\mathfrak X|}\sum_{\mathbf x\in\mathfrak X}\varphi(\mathbf x)
=\int_Z\varphi\,\mathrm d\mu
\qquad(\varphi\in\mathcal V).
\end{equation}
Thus, a $\mathcal V$-design provides a finite-point analogue of the mean value theorem for integrals: instead of representing the integral average of a single function by its value at one point, it reproduces the integrals of all functions in $\mathcal{V}$ by equal-weight averaging over the same finite set of points.
More generally, a nonempty finite subset $\mathfrak{X} \subseteq Z$ equipped with positive weights $m'_{\mathbf x}$, $\mathbf{x} \in \mathfrak{X}$, is called a \emph{weighted $\mathcal{V}$-design} if
\begin{equation}\label{eq:weighted-design}
\sum_{\mathbf x\in\mathfrak X}m'_{\mathbf x}\varphi(\mathbf x)
=\int_Z\varphi\,\mathrm d\mu
\qquad(\varphi\in\mathcal V).
\end{equation}
Since $\mathbf{1} \in \mathcal{V}$, the weights of every weighted $\mathcal{V}$-design sum to $1$, and every $\mathcal{V}$-design is a weighted $\mathcal{V}$-design with equal weights $m'_{\mathbf x}=1/|\mathfrak X|$.

A fundamental existence theorem of Seymour and Zaslavsky~\cite{SeymourZaslavsky1984} states that a \(\mathcal V\)-design exists whenever \(Z\) is path-connected, \(\mu\) assigns positive measure to every nonempty relatively open subset of \(Z\), and every function in \(\mathcal V\) is continuous.
This result provides a general existence principle for designs under natural topological and measure-theoretic hypotheses.

For a nonnegative integer \(t\), put
\[
\mathcal P_t(Z;\mathbb R)
:=\{f|_Z:f\in\mathbb R[X_1,\ldots,X_d],\ \deg f\leq t\}.
\]
When \(Z=\mathbb S^{d-1}\) and \(\mu\) is the normalized spherical measure, a \(\mathcal P_t(\mathbb S^{d-1};\mathbb R)\)-design is called a \emph{spherical \(t\)-design}~\cite{DelsarteGoethalsSeidel1977}.
Spherical designs are related to cubature formulas allowing general weights, frame theory~\cite[Section~6.4]{Waldron2018}, and isometric embeddings of finite-dimensional Banach spaces~\cite{LyubichVaserstein1993}.
Designs have also been studied on real Grassmannians~\cite{BachocCoulangeonNebe2002,BachocBannaiCoulangeon2004}, unitary groups~\cite{GrossAudenaertEisert2007,RoyScott2009}, complex Grassmannians~\cite{Roy2008}, and orthogonal and special orthogonal groups~\cite{GraefPotts2009,Okuda2014,HashagenFlammiaGrossWallman2018}.
For a given $t$, basic problems of designs are to determine whether $\mathcal P_t(Z;\mathbb R)$-designs exist, to construct them, to find their minimum size, and to describe the designs of that size.
For a given number of points $N$, the problem is to find the largest $t$ for which a $\mathcal P_t(Z;\mathbb R)$-design with $N$ points exists~\cite[Sections~1.5 and~1.6]{HarpePache2005}.
We study the existence and explicit construction of $\mathcal V$-designs whose points all have rational coordinates.

\begin{definition}\label{def:rational-v-design}
Put $Z(\mathbb{Q}) := Z\cap\mathbb{Q}^d$.
A $\mathcal{V}$-design $\mathfrak{X}$ is called \emph{rational} if $\mathfrak{X} \subseteq Z(\mathbb{Q})$.
\end{definition}

In complex matrix spaces, we use the real and imaginary parts of the entries as coordinates, so rationality means that every matrix entry belongs to $\mathbb{Q}(\mathrm i)$.
Although rational points are dense in many spaces of interest, density guarantees approximation rather than exact integration.
Since a design must reproduce prescribed integrals exactly, the existence of a design over \(\mathbb R\) does not, in general, imply the existence of a rational one.
A rational $\mathcal{P}_t(Z;\mathbb{R})$-design requires the integrals of all polynomials with rational coefficients and degree at most $t$ to be rational.

Arithmetic design theory explores the connections between geometric design problems and arithmetic, including the existence and construction of rational designs and cubature formulas.
Around the same time, rational designs began to be studied from different perspectives, as illustrated by the work of Sawa--Uchida and Cui--Xia--Xiang.
Sawa--Uchida~\cite{SawaUchida2019,SawaUchida2020} approached the subject from the perspectives of algebraic geometry and the theory of orthogonal polynomials, focusing primarily on one-dimensional problems, whereas Cui--Xia--Xiang~\cite{CuiXiaXiang2019} emphasized geometric and number-theoretic aspects, with a focus on higher-dimensional spaces.
A fundamental question in this area, posed by Bannai~\cite[Section~1.2]{CuiXiaXiang2019}, is the following.

\begin{bannaiproblem}
For which spaces $Z$ and measures $\mu$ does there exist a rational $\mathcal{P}_t(Z;\mathbb{R})$-design for every integer $t\geq1$?
\end{bannaiproblem}

\theoremstyle{definition}
\newtheorem*{tzconjecture}{Conjecture (Thuswaldner--Ziegler)}
\newtheorem*{mahabadugeproblem}{Problem (Mahabaduge)}
\theoremstyle{remark}

Cui--Xia--Xiang~\cite[Lemmas~3.4 and~3.7]{CuiXiaXiang2019} gave sufficient conditions for the existence, at every strength, of designs supported on rational points with positive rational weights, called \emph{rationally weighted designs}.
Under additional assumptions, including openness of the space, they deduced the existence of
equi-weighted rational designs
from that of rationally weighted designs~\cite[Theorem~1.4]{CuiXiaXiang2019}.

On spheres, their results imply the existence of rationally weighted designs of every strength. Such designs also have applications beyond design theory: for example, Gromov~\cite[Section~13]{Gromov2025} used rationally weighted spherical $4$-designs in a curvature optimization problem in differential geometry.

By contrast, the equal-weight existence theorem of Cui--Xia--Xiang does not apply directly to spheres because of its openness assumption.
The existence of rational spherical designs poses a nontrivial arithmetic problem.
Cui et al. investigated the existence of rational spherical $4$-designs through polynomial equations.
For the relation between rational designs and Hilbert--Kamke-type equations, see \cite[Sections~5 and~7.1]{CuiXiaXiang2019} and Mishima et al.~\cite{MishimaLuSawaUchida2026}.
Cui et al. noted that parametrizing rational points on the circle leads to a Diophantine equation whose degree exceeds the number of variables, and reported that they could not apply the classical circle method of analytic number theory to this equation.
This led them to pose the following open question~\cite[Section~1.3]{CuiXiaXiang2019}:

\begin{cxxproblem}
Does a rational spherical $4$-design exist on the unit circle $\mathbb{S}^1$?
\end{cxxproblem}

Xiang~\cite[Theorem~1.1]{Xiang2022} subsequently established the existence of rational spherical designs of strengths $5$, $7$, and $11$ in all sufficiently high dimensions, using an explicit construction that builds on known designs in lower dimensions.
He also proposed the following conjecture for arbitrary strength.

\begin{xiangconjecture}[{\cite[Conjecture~1.4]{Xiang2022}}]
For every positive integer $t$, there exists an integer $d\geq2$ such that $\mathbb{S}^{d-1}$ admits a rational spherical $t$-design.
\end{xiangconjecture}

We prove the existence of rational designs at every strength on spheres, real and complex Grassmannians, and orthogonal, special orthogonal, and unitary groups.
Our proof uses a general theorem that converts rationally weighted designs into
equi-weighted rational designs
under suitable hypotheses on a group action (\cref{thm:integer-to-equal-weight}).
For ellipsoids defined over $\mathbb Q$ in ambient dimension at least $2$, we also prove that the existence of a rational point implies the existence of rational designs of every strength (\cref{thm:rational-ellipsoidal-existence}).
On $\mathbb{S}^1$ and $\mathbb{S}^2$, we also give explicit constructions attaining the minimum sizes at several low strengths.
In particular, for every pair of integers $d\geq2$ and $t\geq1$, there exists a rational spherical $t$-design on $\mathbb{S}^{d-1}$ (\cref{thm:rational-spherical-existence}).
This solves the existence problem for rational spherical $4$-designs on the circle and strengthens Xiang's conjecture to existence in every dimension.

We start with the rationally weighted designs of Cui--Xia--Xiang.
To obtain designs with equal weights, we replace repeated points with distinct rational points while preserving the design conditions.
For this step, we generalize the classical Prouhet--Tarry--Escott problem to find finite sets with equal character sums.

For given integers $n>t\geq1$, the classical Prouhet--Tarry--Escott (PTE) problem is to find two disjoint multisets of $n$ integers with equal sums of $r$th powers for $r=1,\ldots,t$.
Prouhet~\cite{Prouhet1851} showed how to make these power sums agree for any number of sets of equal size.
Jacroux~\cite[Lemma~1]{Jacroux1995} developed this idea further by showing how to combine such constructions to obtain solutions of higher degree.

Alpers and Tijdeman~\cite{AlpersTijdeman2007} introduced a multidimensional generalization in which the elements of the two multisets are integer vectors and the sums of all monomials up to a prescribed total degree agree.
We formulate a common framework encompassing these variants by replacing power functions with characters on abelian groups and allowing an arbitrary number of blocks.

Let $\mathbb{K}$ be a field of characteristic $0$.
A character is a group homomorphism $\chi : G \longrightarrow \mathbb{K}^\times$, where $G$ is an abelian group and $\mathbb{K}^\times$ is the multiplicative group of nonzero elements of $\mathbb{K}$.
Let $\mathcal{S}$ be a finite set of characters containing the trivial character $g\mapsto 1$ and satisfying $\bigcap_{\chi\in\mathcal{S}}\ker\chi=\{1_G\}$, where $1_G$ is the identity of $G$.
For an integer $q\ge 2$, let $A_1,\ldots,A_q$ be nonempty finite multisets in $G$ with pairwise disjoint supports.
Suppose that, counting multiplicities,
\[
\begin{gathered}
|A_1| = \cdots = |A_q|, \\
\sum_{g\in A_1} \chi(g) = \cdots = \sum_{g\in A_q} \chi(g) \qquad (\chi \in \mathcal{S}).
\end{gathered}
\]
We call $(A_1,\ldots,A_q)$ a \emph{$(G,\mathcal{S})$-PTE solution}.
Each $A_i$ is called a \emph{block}, and a $(G,\mathcal{S})$-PTE solution in which every block is a set is called a \emph{$(G,\mathcal{S})$-PTE set solution}.

We prove a general existence theorem for $(G,\mathcal S)$-PTE set solutions on infinite abelian groups $G$.
Under the hypotheses of \cref{thm:infinite-existence}, with the values of the characters in $\mathcal S$ lying in a common number field, such a solution exists for every prescribed number $q\ge 2$ of blocks; moreover, the blocks can be chosen to avoid any prescribed finite subset of $G$.
The proof uses the pigeonhole principle, comparing the number of candidate subsets with the number of possible combinations of their sizes and character sums.
For free abelian groups $G$ of finite rank at least $2$, with the values of the characters in $\mathcal S$ lying in a common number field, the same counting argument gives $(G,\mathcal{S})$-PTE set solutions with any prescribed number of blocks at least $2$ (\cref{prop:free-abelian-rank-two-existence}).
One conceptual innovation of this paper is the use of the sunflower lemma of Erd\H{o}s--Rado~\cite{ErdosRado1960}, a result in extremal set theory, to obtain pairwise disjoint sets of rational points with equal moment sums.
This gives an alternative proof of the existence of $(G,\mathcal S)$-PTE set solutions on the rational unit circle (\cref{prop:circle-moment-pte,prop:circle-pte-lifting}).
We also use this method to construct a rational spherical $11$-design on the unit circle, invariant under rotation by $\pi/2$, and obtain an explicit upper bound on its size (\cref{prop:circle-strength-eleven-quantitative}).

We apply this PTE framework to the construction of rational designs and to problems in algebraic number theory and discrete geometry.

For rational designs, the PTE framework allows repeated points to be replaced by distinct rational points while preserving the design conditions.
In addition to the existence results above, we study explicit constructions and arithmetic restrictions on size.
On the circle, the vertices of a regular $(t+1)$-gon form a spherical $t$-design~\cite{DelsarteGoethalsSeidel1977}, so $t+1$ points suffice when real coordinates are allowed.
Let $\mathcal N_{1,\mathbb Q}(t)$ be the set of sizes of rational spherical $t$-designs on the circle. For each fixed positive squarefree integer $\ell$, \cref{thm:conic-exponential-lower-bound} gives an exponential lower bound for the size of rational designs on the ellipse $x_1^2+\ell x_2^2=1$. For $\ell=1$, it gives
\[
\log\min\mathcal N_{1,\mathbb Q}(t)\geq t+o(t)
\qquad(t\to\infty).
\]

At strengths $4$ through $7$, we determine all possible sizes by combining local divisibility conditions with constructions attaining the corresponding lower bounds. Specifically,
\[
\mathcal N_{1,\mathbb Q}(4)=2304\mathbb Z_{>0},\quad
\mathcal N_{1,\mathbb Q}(5)=4608\mathbb Z_{>0},\quad
\mathcal N_{1,\mathbb Q}(6)=\mathcal N_{1,\mathbb Q}(7)=9216\mathbb Z_{>0};
\]
see \cref{thm:circle-cardinality-spectra}.

Independently and contemporaneously, Cheng, Hu, and Yu~\cite{ChengHuYu2026} obtained the same minimum sizes at strengths $4,5,6$ and $7$, and constructed two nonisometric optimal designs at each of these strengths. They also show that the size of every rational spherical $t$-design on the circle is divisible by an explicit integer $D_t$ with $\log D_t=(\tfrac12+o(1))t\log t$ as $t\to\infty$ \cite[Theorem~2.14 and Corollary~2.15]{ChengHuYu2026}. Their divisibility theorem therefore gives a stronger asymptotic lower bound on the circle. At strengths $4$ through $7$, $D_t$ equals the minimum size. Taking unions of pairwise disjoint copies of an optimal design obtained by rational rotations, as in \cref{prop:conic-cardinality-additivity}, realizes every positive multiple of that size. Thus the full sets of possible sizes at these strengths also follow from their results.

Our general existence theorem for $(G,\mathcal S)$-PTE set solutions provides a common method for constructing rational designs of arbitrary strength on a variety of spaces and for obtaining the applications to algebraic number theory and discrete geometry described below.

On $\mathbb S^2$, we construct a $192$-point rational spherical $5$-design and prove that its size is minimal (\cref{prop:minimum-rational-sphere-two-five}).

As an application to algebraic number theory, we study representations as sums of units.
A unit of a ring is an element with a multiplicative inverse in that ring.
Narkiewicz~\cite[Problem~18]{Narkiewicz1990} asks for a classification of the number fields in which every element of the ring of integers is a finite sum of distinct units.
The same representation problem has also been studied in subrings of number fields.
For particular quadratic, cubic, and quartic number fields, every element of the ring of integers is known to be a finite sum of distinct units~\cite[Section~1]{KrennThuswaldnerZiegler2013}.
Thuswaldner--Ziegler proposed the following conjecture for rings arising from a family of cubic number fields.

\begin{tzconjecture}[{\cite[Conjecture~7.3]{ThuswaldnerZiegler2011}}]
Let $a\in\mathbb Z$, and let $\varrho$ be a root of
\[
X^3-(a-1)X^2-(a+2)X-1.
\]
Then every element of $\mathbb Z[\varrho]$ is a finite sum of distinct units of this ring.
\end{tzconjecture}

We prove this conjecture using $(G,\mathcal S)$-PTE set solutions (\cref{cor:thuswaldner-ziegler-conjecture}).
Two units $u,v$ are called multiplicatively independent if $u^{k_1}v^{k_2}=1$, with $k_1,k_2\in\mathbb Z$, implies $k_1=k_2=0$.
More generally, for a subring of a number field containing $1$ and two multiplicatively independent units, we prove that if every element is a finite sum of units with repetitions allowed, then every element is a finite sum of distinct units (\cref{prop:distinct-unit-sums}).
For number fields whose rings of integers contain two multiplicatively independent units, Narkiewicz's problem therefore reduces to determining whether every element of the ring of integers is a finite sum of units with repetitions allowed.

As an application to discrete geometry, we study the possible numbers of edges of equilateral polygons in a given planar lattice.
For a positive squarefree integer $h\equiv3\pmod4$, we consider the lattice $\mathbb Z\times\sqrt h\,\mathbb Z$.
These lattices contain convex equilateral $n$-gons for every even integer $n\geq4$~\cite[Theorem~4.1]{Maehara2019}.
Convex equilateral polygons have also been constructed for certain lattices and small odd numbers of edges~\cite[Section~3]{IinoSakiyama2026}.
Mahabaduge~\cite[Theorem~1]{Mahabaduge2025} proved that each lattice in this family contains equilateral polygons with every sufficiently large number of edges when convexity is not required.
He asked whether the same conclusion holds when convexity is required~\cite[Open Problem~1]{Mahabaduge2025}.
We prove that, for each such $h$, this lattice contains convex equilateral polygons with every sufficiently large number of edges, answering his problem affirmatively (\cref{cor:lattice-convex-equilateral-polygons}).
We obtain this conclusion by using $(G,\mathcal S)$-PTE set solutions to construct an odd number of distinct nonzero lattice vectors of equal length whose sum is zero.

Alongside these applications, we study minimum block sizes and explicit constructions for the $(G,\mathcal S)$-PTE problem itself.
For the classical PTE problem, we construct a family of ideal solutions of degree $5$ (\cref{thm:new-degree-five-family}).
The family contains infinitely many solutions that are affinely inequivalent to the famous solutions from either of the two constructions of Choudhry~\cite{Choudhry2003,Choudhry2017} considered in this paper (\cref{prop:degree-five-separation}).

\Cref{sec:pte-abelian} develops the $(G,\mathcal{S})$-PTE theory and its applications.
\Cref{sec:rational-design-existence} treats the existence of rational designs, and \cref{sec:local-arithmetic-minimum} studies their possible sizes on the circle.
\Cref{sec:problems} considers ellipsoids and related problems.

Throughout the remainder of the paper, $\mathbb{N} := \{1,2,\ldots\}$, and $\mathbf{e}_j$ denotes the $j$-th standard basis vector.

\section{Prouhet--Tarry--Escott problems on abelian groups} \label{sec:pte-abelian}

\begingroup
\setstretch{1.05}

In this section, we extend the classical PTE problem to the $(G,\mathcal S)$-PTE problem of equating character sums on abelian groups.
After introducing the basic notions, we prove a general existence theorem for set solutions under condition $(\mathrm D)$ (\cref{thm:infinite-existence}).
We adapt the counting argument to obtain set solutions of odd size and apply them to construct convex equilateral polygons in planar lattices (\cref{prop:odd-cardinality-pte,cor:lattice-convex-equilateral-polygons}).
We also use the counting argument to study representations of ring elements as sums of distinct units and prove the Thuswaldner--Ziegler conjecture (\cref{prop:distinct-unit-sums,cor:thuswaldner-ziegler-conjecture}).
We then study ideal solutions on subgroups of the unit circle and the corresponding problems on finite cyclic groups (\cref{thm:circle-minimum,cor:imaginary-quadratic-minimum,prop:prime-cyclic-rigidity,prop:cyclic-jump}).
Finally, we transfer character-sum identities to classical power sums, construct a parametric family of nonsymmetric rational ideal PTE solutions of degree five, and compare it with two constructions of Choudhry (\cref{thm:new-degree-five-family,prop:degree-five-separation}).

\subsection{The $(G,\mathcal S)$-PTE problem and basic notions}\label{subsec:pte-definitions}

Throughout this subsection, let $q\geq2$ be an integer.
Let $\mathbb K$ be a field of characteristic zero and $G$ be an abelian group with identity $1_G$.
We call a homomorphism $G\to \mathbb K^\times$ a character and write $\mathbf1$ for the trivial character, $\mathbf1(g)=1$ for all $g\in G$.
The set of finite multisets on $G$ is
\[
\operatorname{MSet}(G):=\{\xi:G\to\mathbb Z_{\geq0}:\operatorname{supp}(\xi)\text{ is finite}\}.
\]
Write $\xi_A:G\to\mathbb Z_{\geq0}$ for the multiplicity function of $A\in\operatorname{MSet}(G)$.
For a map $\varphi:G\to\mathbb K$, put
\begin{equation}\label{eq:multiset-sum}
\begin{gathered}
\operatorname{supp}(A):=\{g\in G:\xi_A(g)>0\},\qquad
|A|:=\sum_{g\in\operatorname{supp}(A)}\xi_A(g),\\
\sum_{g\in A}\varphi(g):=\sum_{g\in\operatorname{supp}(A)}\xi_A(g)\varphi(g),\qquad
\prod_{g\in A}\varphi(g):=\prod_{g\in\operatorname{supp}(A)}\varphi(g)^{\xi_A(g)}.
\end{gathered}
\end{equation}
We say that finite multisets on $G$ are pairwise disjoint if their supports are pairwise disjoint.

\begin{definition}[$(G,\mathcal S)$-PTE solution]\label{def:abelian-pte}
Let $\mathcal S$ be a finite subset of $\operatorname{Hom}(G,\mathbb K^\times)$ containing the trivial character $\mathbf1$ and satisfying $\bigcap_{\chi\in\mathcal S}\ker\chi=\{1_G\}$.
Let $n$ be a positive integer.
A $q$-tuple $(A_1,\ldots,A_q)$ of pairwise disjoint nonempty finite multisets on $G$, each of size $n$, is called a $(G,\mathcal S)$-PTE solution of size $n$ if
\begin{equation}\label{eq:multiple-pte-character-sums}
\sum_{g\in A_1}\chi(g)=\cdots=\sum_{g\in A_q}\chi(g)
\qquad(\chi\in\mathcal S).
\end{equation}
We write
\[
[A_1]=_{\mathcal S}^{n}\cdots=_{\mathcal S}^{n}[A_q]
\]
for the condition $|A_1|=\cdots=|A_q|=n$ together with \eqref{eq:multiple-pte-character-sums}.
The multisets $A_1,\ldots,A_q$ are called the blocks of the solution.
If every block is a set, the solution is called a $(G,\mathcal S)$-PTE set solution.
\end{definition}

\begin{remark}\label{rem:trivial-solutions}
The condition $\bigcap_{\chi\in\mathcal S}\ker\chi=\{1_G\}$ in \cref{def:abelian-pte} rules out pairs of distinct singleton blocks with equal character sums.
Indeed, for $g,h\in G$,
\[
\chi(g)=\chi(h)\quad(\chi\in\mathcal S)
\quad\Longleftrightarrow\quad
g^{-1}h\in\bigcap_{\chi\in\mathcal S}\ker\chi.
\]
\end{remark}

\begin{remark}\label{rem:classical-pte}
Let $t$ be a positive integer, let $\mathbb K=\mathbb Q$, and put $G=\mathbb Q^\times$.
Let $\mathcal S$ consist of the characters $z\mapsto z^k$ for $0\leq k\leq t$.
Then the defining equations are precisely
\[
\sum_{z\in A_1}z^k=\cdots=\sum_{z\in A_q}z^k
\qquad(0\leq k\leq t),
\]
that is, the classical multiple PTE equations for nonzero rational entries.
Conversely, any classical rational solution can first be translated by a common rational number so that all entries are nonzero and then scaled by a common positive integer to clear denominators.
Translation preserves the equations by the binomial theorem, and scaling preserves them by homogeneity.
See Frolov~\cite[Theorems~I and~III]{Frolov1889}.
\end{remark}

To distinguish multiset and set solutions, we often use the following notation.
\begin{definition}\label{def:min-cardinality}
Define
\[
\begin{aligned}
\kappa_{G,\mathcal S}(q)&:=\min\{|A_1|:(A_1,\ldots,A_q)\text{ is a }(G,\mathcal S)\text{-PTE solution}\},\\
\kappa^{\mathrm{set}}_{G,\mathcal S}(q)&:=\min\{|A_1|:(A_1,\ldots,A_q)\text{ is a }(G,\mathcal S)\text{-PTE set solution}\},
\end{aligned}
\]
where $\min\varnothing:=\infty$.
\end{definition}

\begin{proposition} \label{prop:PTEbound-newton}
Let $q,t$ be positive integers with $q\geq2$.
Let $\chi : G \to \mathbb K^\times$ be an injective character and $\mathcal S=\{\mathbf1,\chi,\ldots,\chi^t\}$.
Then
\begin{equation}\label{eq:PTEbound}
t+1\leq\kappa_{G,\mathcal S}(q)\leq\kappa^{\mathrm{set}}_{G,\mathcal S}(q).
\end{equation}
\end{proposition}
\begin{proof}
Suppose that $(A_1,\ldots,A_q)$ is a solution of size $n\leq t$.
For $1\leq k\leq n$, the defining equations give
\[
\sum_{g\in A_1}\chi(g)^k=\cdots=\sum_{g\in A_q}\chi(g)^k.
\]
Because $\mathbb K$ has characteristic zero, Newton's identities~\cite[Chapter~I, Equation~(2.11')]{Macdonald1995} determine all $n$ elementary symmetric functions from these sums, giving
\[
\prod_{g\in A_1}\bigl(X-\chi(g)\bigr)
=
\prod_{g\in A_2}\bigl(X-\chi(g)\bigr).
\]
Injectivity of $\chi$ then gives $A_1=A_2$ as multisets, contradicting their disjointness and nonemptiness.
The second inequality follows because every set solution is a multiset solution.
\end{proof}

\begin{definition}[Ideal $(G,\mathcal S)$-PTE solution]\label{def:ideal-solution}
Under the hypotheses of \cref{prop:PTEbound-newton}, a $(G,\mathcal S)$-PTE solution of size $t+1$ is called an ideal $(G,\mathcal S)$-PTE solution of degree $t$ \cite{CoppersmithEtAl2024}.
\end{definition}

\begin{definition}[Affine equivalence and symmetry]\label{def:rational-pte-equivalence-symmetry}
Let ${\alpha}:\mathbb Q\to\mathbb Q$ be an invertible rational affine map.
For a finite multiset $A$ on $\mathbb Q$, write ${\alpha_{*}}A$ for its image multiset under ${\alpha}$.
Following \cite[Definition~5.8]{MatsumuraSawa2025}, two rational PTE solutions $(A_1,A_2)$ and $(A'_1,A'_2)$ are \emph{affinely equivalent over $\mathbb Q$} if there is such a common map ${\alpha}$ with
\[
 (A'_1,A'_2)\in\{({\alpha_{*}}A_1,{\alpha_{*}}A_2),({\alpha_{*}}A_2,{\alpha_{*}}A_1)\}.
\]
A two-block rational PTE solution of odd degree is \emph{symmetric} if it is equivalent over $\mathbb Q$ to a solution $(A'_1,A'_2)$ satisfying $-A'_1=A'_1$ and $-A'_2=A'_2$; otherwise it is \emph{nonsymmetric}.
\end{definition}

\Needspace{5\baselineskip}
The following solution is affinely equivalent to the member at $u=-3$ of the family in \cref{thm:new-degree-five-family}, by \eqref{eq:degree-five-test-point-affine-map}, and is used in the proof of \cref{prop:degree-five-separation}.

\begin{example}\label{ex:degree-five-ideal-solution}
The two sets
\begin{equation}
A_1:=\{-2619,-1857,-475,59,1345,3547\}, \qquad A_2:=\{-2685,-1553,-1007,509,1181,3555\}
\label{eq:degree-five-test-point-concise}
\end{equation}
form a nonsymmetric ideal PTE solution of degree five.
The sets are disjoint and each has six elements.
Their power sums are
\[
\begin{array}{c|r}
k&\displaystyle\sum_{z\in A_1}z^k=\sum_{z\in A_2}z^k\\ \hline
1&0\\
2&24\,926\,950\\
3&22\,583\,808\,000\\
4&220\,550\,200\,239\,814\\
5&420\,518\,758\,589\,030\,400
\end{array}
\]
Both $A_1$ and $A_2$ have mean zero.
For either block $A\in\{A_1,A_2\}$, its image under an invertible rational affine map ${\alpha}$ has mean ${\alpha(0)}$ and third centered power sum
\[
\sum_{z\in A}\bigl({\alpha}(z)-{\alpha(0)}\bigr)^3
={\bigl(\alpha(1)-\alpha(0)\bigr)}^3\sum_{z\in A}z^3
={\bigl(\alpha(1)-\alpha(0)\bigr)}^3\cdot22\,583\,808\,000\ne0.
\]
Every finite multiset invariant under $z\mapsto-z$ has mean zero and third power sum zero, since its odd power sums cancel.
Hence no common rational affine map makes these blocks invariant under negation, proving nonsymmetry.
\end{example}

\subsection{Set solutions under condition $(\mathrm D)$} \label{subsec:pte-infinite-abelian}

We now prove the existence of set solutions that avoid any prescribed finite subset of $G$.

Let $\mathbb K$ be a number field with ring of integers $\mathcal O_{\mathbb K}$, let $G$ be an abelian group, and put $\mathbb Z^{(\mathbb N)}:=\bigoplus_{i\geq1}\mathbb Z\mathbf e_i$.
Here $\mathbf e_i$ is the $i$-th standard basis vector.
\begin{definition}[Condition $(\mathrm D)$]\label{def:condition-D}
Let $\nu:G\to\mathbb Z^{(\mathbb N)}$ be a homomorphism, and let $(g_i)_{i\geq1}$ be a sequence in $G$.
The pair $(\nu,(g_i)_{i\geq1})$ satisfies condition $(\mathrm D)$ if
\begin{equation}\label{eq:condition-D}
\nu(g_i)\in\mathbb Z_{>0}\mathbf e_i\qquad(i\geq1).
\tag{$\mathrm D$}
\end{equation}
\end{definition}

For $i\geq1$, let $\nu_i:G\to\mathbb Z$ be the $i$-th coordinate of $\nu$, so that
\begin{equation}\label{eq:nu-coordinates}
\nu(g)=(\nu_i(g))_{i\geq1}\qquad(g\in G).
\end{equation}
For $\mathbf k=(k_j)_{j\geq1}\in\mathbb Z^{(\mathbb N)}$, put
\[
\operatorname{supp}(\mathbf k):=\{i\in\mathbb N:k_i\neq0\}.
\]
By condition $(\mathrm D)$, for each finite $I\subseteq\mathbb N$, the products $\prod_{i\in I}g_i^{k_i}$ are distinct as $(k_i)_{i\in I}$ ranges over $\mathbb Z^I$.

The following is the main result of this subsection.

\begin{theorem}\label{thm:infinite-existence}
Suppose that $(\nu,(g_i)_{i\geq1})$ satisfies condition $(\mathrm D)$, and let $\mathcal S$ be a nonempty finite subset of $\operatorname{Hom}(G,\mathbb K^\times)$.
Let $q\geq2$ be an integer, and let $E\subseteq G$ and $I_0\subseteq\mathbb N$ be finite.
Then $G\setminus E$ contains pairwise disjoint nonempty finite subsets $A_1,\ldots,A_q$ of the same size such that
\[
\sum_{g\in A_1}\chi(g)=\cdots=\sum_{g\in A_q}\chi(g)
\qquad(\chi\in\mathcal S),
\]
and, for some finite set $I\subseteq\mathbb N\setminus I_0$,
\[
\bigcup_{j=1}^{q} A_j
\subseteq
\left\{\prod_{i\in I}g_i^{k_i}\ \middle|\ (k_i)_{i\in I}\in\mathbb Z_{>0}^{I}\right\}.
\]
If $\mathbf1\in\mathcal S$ and $\bigcap_{\chi\in\mathcal S}\ker\chi=\{1_G\}$, then $(A_1,\ldots,A_q)$ is a $(G,\mathcal S)$-PTE set solution.
\end{theorem}

The following lemma is useful in the proof of \cref{thm:infinite-existence}.

\begin{lemma}\label{lem:number-field-character-sum-tuples}
Fix a nonempty finite $I\subseteq\mathbb N$ and elements ${b_{j}},{b_{j}}^{\prime}\in\mathbb K$ for $j\in I$.
There exists a constant $c=c(\mathbb K,I,({b_{j}},{b_{j}}^{\prime})_{j\in I})>0$ such that, for every integer ${\mathsf{k}}\geq0$,
\[
\left|\left\{\left(|J|,\left(\sum_{(k_1,k_2)\in J}{b_{j}}^{k_1}({b_{j}}^{\prime})^{k_2}\right)_{j\in I}\right):J\subseteq\{0,\ldots,{\mathsf{k}}\}^2\right\}\right|
\leq\exp(c({\mathsf{k}}+1)).
\]
\end{lemma}

\begin{proof}
Choose ${\mathsf{b}}\in\mathbb Z_{>0}$ such that ${\mathsf{b}}{b_{j}},{\mathsf{b}}{b_{j}}^{\prime}\in\mathcal O_{\mathbb K}$ for every $j\in I$.
For ${\mathsf{k}}\geq0$ and $J\subseteq\{0,\ldots,{\mathsf{k}}\}^2$,
\[
{\mathsf{b}}^{2{\mathsf{k}}}\sum_{(k_1,k_2)\in J}{b_{j}}^{k_1}({b_{j}}^{\prime})^{k_2}\in\mathcal O_{\mathbb K},
\]
since every summand after multiplication by ${\mathsf{b}}^{2{\mathsf{k}}}$ is
\[
{\mathsf{b}}^{2{\mathsf{k}}-k_1-k_2}({\mathsf{b}}{b_{j}})^{k_1}({\mathsf{b}}{b_{j}}^{\prime})^{k_2},
\]
whose exponents are nonnegative.
Put
\[
c_0:=\max_{\substack{j\in I\\\sigma:\mathbb K\hookrightarrow\mathbb C}}
\max\{1,|\sigma({\mathsf{b}}{b_{j}})|,|\sigma({\mathsf{b}}{b_{j}}^{\prime})|,{\mathsf{b}}\}.
\]
Since $c_0\geq {\mathsf{b}}$, each embedded summand has absolute value at most $c_0^{2{\mathsf{k}}}$.
Thus, for every embedding $\sigma:\mathbb K\hookrightarrow\mathbb C$,
\[
\left|\sigma\!\left({\mathsf{b}}^{2{\mathsf{k}}}\sum_{(k_1,k_2)\in J}{b_{j}}^{k_1}({b_{j}}^{\prime})^{k_2}\right)\right|
\leq({\mathsf{k}}+1)^2c_0^{2{\mathsf{k}}}\leq\exp(c_1({\mathsf{k}}+1)),
\]
where $c_1>0$ is independent of ${\mathsf{k}}$, $J$, and $j$.
The Minkowski embedding identifies $\mathcal O_{\mathbb K}$ with a lattice of real dimension $[\mathbb K:\mathbb Q]$ \cite[Chapter~I, Sections~4--5]{Neukirch1999}.
The number of points of this fixed lattice satisfying these coordinate bounds is at most $\exp(c_2({\mathsf{k}}+1))$, with $c_2>0$ independent of ${\mathsf{k}}$, $J$, and $j$.
There are $({\mathsf{k}}+1)^2+1$ possibilities for $|J|$, and multiplication by ${\mathsf{b}}^{2{\mathsf{k}}}\ne0$ is injective.
Consequently the number of tuples in the statement is at most
\[
\bigl(({\mathsf{k}}+1)^2+1\bigr)\exp\bigl(|I|c_2({\mathsf{k}}+1)\bigr)\leq\exp(c({\mathsf{k}}+1))
\]
for a suitable constant $c>0$ with precisely the stated dependence.
\end{proof}

\Cref{lem:number-field-character-sum-tuples} bounds the number of character-sum tuples used in the pigeonhole argument for \cref{thm:infinite-existence}.

\begin{proof}[Proof of \cref{thm:infinite-existence}]
We first construct two blocks using two indices and then repeat the construction to obtain $q$ blocks.

For the construction of two blocks, put $[|\mathcal S|]:=\{1,\ldots,|\mathcal S|\}$ and enumerate the characters as $(\chi_j)_{j\in[|\mathcal S|]}$.
For distinct $i_1,i_2\in\mathbb N$ and an integer $\mathsf{k}\geq0$, the map
\[
\{0,\ldots,\mathsf{k}\}^2
\longrightarrow
\{g_{i_1}^{k_1}g_{i_2}^{k_2}:0\leq k_1,k_2\leq\mathsf{k}\},
\qquad
(k_1,k_2)\longmapsto g_{i_1}^{k_1}g_{i_2}^{k_2}
\]
is a bijection by condition $(\mathrm D)$.

Fix distinct indices $i_1,i_2$.
Apply \cref{lem:number-field-character-sum-tuples} with $I=[|\mathcal S|]$, $b_j=\chi_j(g_{i_1})$, and $b_j^{\prime}=\chi_j(g_{i_2})$.
The lemma and the bijectivity of the map give a constant $c>0$, independent of $\mathsf{k}$, such that for every integer $\mathsf{k}\geq0$,
\[
\#\left\{
\left(|A|,\left(\sum_{g\in A}\chi_j(g)\right)_{j\in[|\mathcal S|]}\right)
\ \middle|\
A\subseteq\{g_{i_1}^{k_1}g_{i_2}^{k_2}:0\leq k_1,k_2\leq\mathsf{k}\}
\right\}
\leq \exp\bigl(c(\mathsf{k}+1)\bigr).
\]

Choose $\mathsf{k}$ large enough that $2^{(\mathsf{k}+1)^2}>\exp\bigl(c(\mathsf{k}+1)\bigr)$.
Since the map is a bijection, its image has $2^{(\mathsf{k}+1)^2}$ subsets.
By the pigeonhole principle, there are distinct subsets $A,A'$ of the image such that
\[
|A|=|A'|,
\qquad
\sum_{g\in A}\chi_j(g)=\sum_{g\in A'}\chi_j(g)
\quad(j\in[|\mathcal S|]).
\]

Removing the common part $A\cap A'$ from both sets preserves the equalities of their sizes and character sums.
Since $A\ne A'$ and $|A|=|A'|$, both $A\setminus A'$ and $A'\setminus A$ are nonempty.
Put
\[
A'_0:=g_{i_1}g_{i_2}(A\setminus A'),
\qquad
A'_1:=g_{i_1}g_{i_2}(A'\setminus A).
\]
Then $A'_0,A'_1\subseteq\{g_{i_1}^{k_1}g_{i_2}^{k_2}:k_1,k_2\in\mathbb Z_{>0}\}$, $A'_0\cap A'_1=\varnothing$, and $|A'_0|=|A'_1|>0$.
Moreover,
\[
\sum_{g\in A'_0}\chi_j(g)=\sum_{g\in A'_1}\chi_j(g)
\qquad(j\in[|\mathcal S|]).
\]

We now repeat this construction to obtain $q$ blocks.

By hypothesis, $I_0$ and $E$ are finite and each $\nu(g)$ has finite support, so the set
\[
I_0\cup\bigcup_{g\in E}\operatorname{supp}(\nu(g))
\]
is finite.
For $1\leq j\leq\lceil\log_2q\rceil$, choose indices
\[
i_{j,1},i_{j,2}\in
\mathbb N\setminus\left(I_0\cup\bigcup_{g\in E}\operatorname{supp}(\nu(g))\right),
\]
all distinct, and put
\[
I:=\{i_{j,1},i_{j,2}:1\leq j\leq\lceil\log_2q\rceil\}.
\]
For each $j$, apply the two-block construction with $(i_1,i_2)=(i_{j,1},i_{j,2})$ to obtain disjoint nonempty sets
\[
A'_{j,0},A'_{j,1}
\subseteq
\{g_{i_{j,1}}^{k_1}g_{i_{j,2}}^{k_2}:k_1,k_2\in\mathbb Z_{>0}\}
\]
of the same size and with equal character sums.
For each $\boldsymbol{\varepsilon}=(\varepsilon_j)_{j=1}^{\lceil\log_2q\rceil}\in\{0,1\}^{\lceil\log_2q\rceil}$, put
\[
A_{\boldsymbol{\varepsilon}}
:=\prod_{j=1}^{\lceil\log_2q\rceil}A'_{j,\varepsilon_j}.
\]

The sets $A_{\boldsymbol{\varepsilon}}$ are pairwise disjoint.
Indeed, take $a_j\in A'_{j,\varepsilon_j}$ for each $j$.
Writing the $j$-th factor as $a_j=g_{i_{j,1}}^{k_1}g_{i_{j,2}}^{k_2}$, we obtain from condition $(\mathrm D)$
\[
\nu_{i_{j,\ell}}\bigl(a_1\cdots a_{\lceil\log_2q\rceil}\bigr)
=k_\ell\,\nu_{i_{j,\ell}}(g_{i_{j,\ell}})
\qquad(\ell=1,2).
\]
Since $\nu_{i_{j,\ell}}(g_{i_{j,\ell}})>0$, the product determines $k_1,k_2$ uniquely and hence determines $a_j$ uniquely.
Moreover, $A'_{j,0}\cap A'_{j,1}=\varnothing$, so $a_j$ also determines $\varepsilon_j$ uniquely.

Consequently,
\[
|A_{\boldsymbol{\varepsilon}}|
=\prod_{j=1}^{\lceil\log_2q\rceil}|A'_{j,\varepsilon_j}|.
\]
For each $\chi\in\mathcal S$, the multiplicativity of $\chi$ gives
\[
\sum_{g\in A_{\boldsymbol{\varepsilon}}}\chi(g)
=\prod_{j=1}^{\lceil\log_2q\rceil}
\sum_{g\in A'_{j,\varepsilon_j}}\chi(g).
\]
Since $A'_{j,0}$ and $A'_{j,1}$ have the same size and equal character sums for each $j$, the size and character sums of $A_{\boldsymbol{\varepsilon}}$ are independent of the choice of $\boldsymbol{\varepsilon}$.

The set $I$ is a finite subset of $\mathbb N\setminus I_0$, and each element of $A_{\boldsymbol{\varepsilon}}$ is a product of positive powers of all the $g_i$ with $i\in I$.
It follows from condition $(\mathrm D)$ and $I\cap\bigcup_{g\in E}\operatorname{supp}(\nu(g))=\varnothing$ that $A_{\boldsymbol{\varepsilon}}\cap E=\varnothing$ for every $\boldsymbol{\varepsilon}$.

Since $2^{\lceil\log_2q\rceil}\geq q$, choose $q$ distinct $\boldsymbol{\varepsilon}$ and denote the corresponding sets by $A_1,\ldots,A_q$.
\end{proof}

The following lemma verifies condition $(\mathrm D)$ for norm-one groups of imaginary quadratic fields; it will be used in \cref{cor:imaginary-quadratic-minimum} and \cref{sec:rational-design-existence}.

\begin{lemma}\label{lem:norm-one-condition-D}
Let $\mathbb K$ be an imaginary quadratic field with ring of integers $\mathcal O_{\mathbb K}$, and put
\[
T_{\mathbb K}:=\{z\in \mathbb K^\times:z\overline z=1\}.
\]
For a nonzero prime ideal $\mathfrak p$ of $\mathcal O_{\mathbb K}$, let $\operatorname{ord}_{\mathfrak p}:\mathbb K^\times\to\mathbb Z$ denote the normalized prime-ideal valuation.
There exist distinct nonzero prime ideals $\mathfrak p_i$ of $\mathcal O_{\mathbb K}$, a homomorphism $\nu:T_{\mathbb K}\to\mathbb Z^{(\mathbb N)}$, and elements $g_i\in T_{\mathbb K}$ such that
\[
\nu(z)=(\operatorname{ord}_{\mathfrak p_i}(z))_{i\geq1}, \qquad \operatorname{ord}_{\mathfrak p_i}(g_j)=\begin{cases}1,&i=j,\\0,&i\ne j\end{cases} \quad(i,j\geq1).
\]
Thus $\nu(g_i)=\mathbf e_i$ for every $i\geq1$, so condition $(\mathrm D)$ holds.
\end{lemma}

\begin{proof}
By the Chebotarev density theorem, infinitely many rational primes split completely in the Hilbert class field of $\mathbb K$.
A nonzero prime ideal of $\mathcal O_{\mathbb K}$ splits completely in the Hilbert class field of $\mathbb K$ exactly when it is principal; hence there are distinct rational primes ${\mathsf{p}}_i$ with ${\mathsf{p}}_i\mathcal O_{\mathbb K}=\mathfrak p_i\overline{\mathfrak p_i}$ and $\mathfrak p_i=({p_{\mathbb K,i}})$ principal; see for example~\cite[Chapter~VI, Sections~6--7 and Chapter~VII, Section~13]{Neukirch1999}.
Because the rational primes ${\mathsf{p}}_i$ are distinct, $\mathfrak p_i$ is different from both $\mathfrak p_j$ and $\overline{\mathfrak p_j}$ whenever $i\ne j$.
Put $g_i:={p_{\mathbb K,i}}/\overline{{p_{\mathbb K,i}}}\in T_{\mathbb K}$.
The identities $({p_{\mathbb K,i}})=\mathfrak p_i$ and $(\overline{{p_{\mathbb K,i}}})=\overline{\mathfrak p_i}$ give
\[
\operatorname{ord}_{\mathfrak p_i}(g_j)=\begin{cases}1,&i=j,\\0,&i\ne j\end{cases}.
\]
Each element of $\mathbb K^\times$ has nonzero valuation at only finitely many prime ideals, so $\nu(z):=(\operatorname{ord}_{\mathfrak p_i}(z))_{i\geq1}$ is a well-defined homomorphism to $\mathbb Z^{(\mathbb N)}$.
\end{proof}

\subsection{Applications to discrete geometry and algebraic number theory}\label{subsec:combinatorial-number-theoretic-applications}

The construction of sets with equal sums in \cref{subsec:pte-infinite-abelian} gives new results beyond design theory.
We first strengthen \cref{thm:infinite-existence} as follows.

\begin{proposition}\label{prop:odd-cardinality-pte}
In \cref{thm:infinite-existence}, assume also that $\mathbf1\in\mathcal S$ and $\bigcap_{\chi\in\mathcal S}\ker\chi=\{1_G\}$.
Then the $(G,\mathcal S)$-PTE set solution $(A_1,\ldots,A_q)$ in that theorem can be chosen to satisfy
\[
|A_i|\equiv1\pmod2\qquad(1\leq i\leq q).
\]
\end{proposition}

\begin{proof}
Use the notation and hypotheses of \cref{thm:infinite-existence}.
Choose two distinct elements $g_r,g_s$ of the sequence $(g_i)_{i\geq1}$ in condition $(\mathrm D)$, and put $\Omega_{\mathsf k}:=\{g_r^a g_s^b:0\leq a,b\leq\mathsf k\}$ for each nonnegative integer $\mathsf k$.
Using \cref{lem:number-field-character-sum-tuples} as in the proof of \cref{thm:infinite-existence}, for sufficiently large $\mathsf k$ we obtain a family $\mathcal B\subseteq2^{\Omega_{\mathsf k}}$ of subsets with a common size and common character sums such that
\[
|\mathcal B|>2^{|\Omega_{\mathsf k}|/2}.
\]
Suppose that $|A\setminus A'|$ is even for all $A,A'\in\mathcal B$.
Let $m$ be their common size, and let $\mathbf t_A\in\mathbb F_2^{\Omega_{\mathsf k}}$ be the indicator vector of $A$ over the two-element field $\mathbb F_2$.
Then $\mathbf t_A\cdot\mathbf t_{A'}=|A\cap A'|\equiv m\pmod2$.
Fixing $A_*\in\mathcal B$, we therefore have
\[
(\mathbf t_A+\mathbf t_{A_*})\cdot(\mathbf t_{A'}+\mathbf t_{A_*})=0
\qquad(A,A'\in\mathcal B).
\]
Hence
\[
\mathcal U:=\operatorname{span}_{\mathbb F_2}\{\mathbf t_A+\mathbf t_{A_*}:A\in\mathcal B\}\subseteq \mathcal U^\perp,
\]
where $\perp$ denotes the orthogonal complement with respect to the standard bilinear form.
This form is nondegenerate, and the map $A\mapsto\mathbf t_A+\mathbf t_{A_*}$ is injective, so
\[
|\mathcal B|\leq2^{\dim\mathcal U}\leq2^{|\Omega_{\mathsf k}|/2},
\]
a contradiction.
Thus $|A\setminus A'|$ is odd for some $A,A'\in\mathcal B$, and removing their intersection gives a two-block $(G,\mathcal S)$-PTE set solution $(A\setminus A',A'\setminus A)$ of a common odd size.
The rest of the proof proceeds as in the proof of \cref{thm:infinite-existence}.
\end{proof}

\subsubsection{Convex equilateral polygons in lattices}\label{subsubsec:lattice-equilateral-polygons}

A \emph{planar lattice} is the set of all integer linear combinations of two linearly independent vectors in $\mathbb R^2$.
A simple closed polygonal curve in the plane with exactly $n\geq3$ edges of equal length is called an \emph{equilateral $n$-gon}.
It is called a \emph{convex equilateral $n$-gon} if the convex hull of its vertex set has exactly $n$ vertices.
A lattice is said to contain a polygon if all its vertices belong to the lattice.
For a positive integer $h$, put $\Lambda_h:=\mathbb Z\times\sqrt h\,\mathbb Z$.

\begin{lemma}[Iino--Sakiyama~\cite{IinoSakiyama2026}, Lemma~2.1]\label{lem:lattice-equilateral-vectors}
For a planar lattice $L$ and an integer $n\geq3$, the lattice $L$ contains a convex equilateral $n$-gon if and only if there are distinct nonzero vectors $\mathbf v_1,\ldots,\mathbf v_n\in L$ satisfying
\[
\|\mathbf v_1\|=\cdots=\|\mathbf v_n\|,
\qquad
\sum_{i=1}^n\mathbf v_i=\mathbf0.
\]
\end{lemma}

It is known that $\Lambda_h$ contains a convex equilateral $n$-gon for every even $n\geq4$~\cite[Theorem~4.1]{Maehara2019}.
Mahabaduge~\cite[Open Problem~1]{Mahabaduge2025} asks whether, for every positive integer $m\equiv3\pmod4$, the lattice $\Lambda_m$ contains a convex equilateral $n$-gon for every sufficiently large integer $n$.
The following corollary settles the squarefree case.

\begin{corollary}\label{cor:lattice-convex-equilateral-polygons}
For every squarefree positive integer $h\equiv3\pmod4$, the lattice $\Lambda_h$ contains a convex equilateral $n$-gon for every sufficiently large integer $n$.
\end{corollary}

\begin{proof}
By Maehara~\cite[Theorem~4.1(a)]{Maehara2019}, $\Lambda_h$ contains a convex equilateral $n$-gon for every even $n\geq4$.
By Iino and Sakiyama~\cite[Lemma~2.3]{IinoSakiyama2026}, if $\Lambda_h$ contains a convex equilateral $n$-gon, then it also contains a convex equilateral $(n+2)$-gon.

Thus it suffices to construct a single convex equilateral polygon with an odd number of edges.
By \cref{lem:lattice-equilateral-vectors}, it is therefore enough to find, for some odd $n\geq3$, $n$ distinct nonzero vectors in $\Lambda_h$ of equal length whose sum is zero.

Identify $\Lambda_h$ with $\mathbb Z+\sqrt{-h}\,\mathbb Z\subset\mathbb C$.
By Maehara~\cite[Theorem~4.1(b)]{Maehara2019}, $\Lambda_h$ contains an equilateral $q$-gon for some odd $q$.
Let $v_1,\ldots,v_q$ be its edge vectors, and put
\[
G=\{z\in\mathbb Q(\sqrt{-h}):|z|=1\},
\qquad
w_i=\frac{v_i}{v_1}\quad(1\leq i\leq q).
\]
Then $(w_i)_{i=1}^q$ is a tuple of vectors in $G$, with repetitions allowed, satisfying $\sum_{i=1}^q w_i=0$.

Let $\mathcal S_1=\{z\mapsto1,\ z\mapsto z\}$ and suppose that there is a $(G,\mathcal S_1)$-PTE set solution $(A_i)_{i=1}^q$ such that $|A_i|\equiv1\pmod2$ for every $i$ and $w_iA_i\cap w_jA_j=\emptyset$ whenever $i\neq j$.
Then $X:=\bigsqcup_{i=1}^q w_iA_i\subseteq G$ is a set of distinct vectors with $|X|=q|A_1|\equiv1\pmod2$, and
\[
\sum_{x\in X}x
=
\left(\sum_{a\in A_1}a\right)\sum_{i=1}^q w_i
=0.
\]

To obtain such a solution, use \cref{lem:norm-one-condition-D} to choose $(\nu,(g_i)_{i\geq1})$ satisfying condition $(\mathrm D)$ for $G=T_{\mathbb Q(\sqrt{-h})}$.
Apply \cref{prop:odd-cardinality-pte} with $I_0=\bigcup_{i=1}^q\operatorname{supp}(\nu(w_i))$ and $E=\varnothing$.
For the resulting $I\subseteq\mathbb N\setminus I_0$, we have $\nu_r(w_i)=0$ for every $r\in I$ and every $i$.
Thus, by condition $(\mathrm D)$, if $w_i a=w_j b$ for $a\in A_i$ and $b\in A_j$, then $a=b$.
Since $A_i\cap A_j=\emptyset$ for $i\neq j$, it follows that $w_iA_i\cap w_jA_j=\emptyset$ for $i\neq j$.

Let $\mathsf b\in\mathbb Z_{>0}$ be a common multiple of the denominators of the coordinates of all elements of $X$ with respect to the basis $1,\sqrt{-h}$.
Then $\{\mathsf b x:x\in X\}\subseteq\Lambda_h$ is the desired set.
\end{proof}

For any positive integer $m\equiv3\pmod4$, write $m=s^2h$ with
$s,h\in\Z_{>0}$ and $h$ squarefree; then $s$ is odd and
$h\equiv3\pmod4$.
Since
\[
s\Lambda_h=s\Z\times\sqrt m\,\Z\subseteq\Lambda_m,
\]
scaling the polygons in \cref{cor:lattice-convex-equilateral-polygons}
by $s$ gives convex equilateral $n$-gons in $\Lambda_m$ for every
sufficiently large integer $n$.
Thus the corollary answers Mahabaduge's
question~\cite[Open Problem~1]{Mahabaduge2025} in full.

\subsubsection{Narkiewicz's problem}\label{subsubsec:narkiewicz-problem}

Narkiewicz~\cite[Problem~18]{Narkiewicz1990} asks for a classification of the number fields in which every element of the ring of integers is a finite sum of distinct units.
In this connection, Thuswaldner--Ziegler~\cite[Conjecture~7.3]{ThuswaldnerZiegler2011} conjectured that, for every integer $a$, if $\varrho$ is a root of
\[
X^3-(a-1)X^2-(a+2)X-1,
\]
then every element of $\mathbb Z[\varrho]$ is a finite sum of distinct units of this ring.
Representations using each unit at most twice are already known~\cite[Theorem~3.1]{KrennThuswaldnerZiegler2013}.
We give a general proposition for obtaining representations as sums of distinct units by using $(G,\mathcal S)$-PTE set solutions, and prove this conjecture as a consequence.

\begin{proposition}\label{prop:distinct-unit-sums}
Let $\mathbb K$ be a number field, and let $R$ be a subring of $\mathbb K$ containing $1$.
Suppose that $R^\times$ contains a subgroup $H$ isomorphic to $\mathbb Z^2$.
Then every element of $R$ is a finite sum of distinct units if and only if
\[
R=\operatorname{span}_{\mathbb Z}(R^\times).
\]
\end{proposition}

\begin{proof}
If every element of $R$ is a finite sum of units, then $R=\operatorname{span}_{\mathbb Z}(R^\times)$.
Conversely, suppose that $R=\operatorname{span}_{\mathbb Z}(R^\times)$.
First, consider the $(G,\mathcal S)$-PTE problem over the number field $\mathbb K$.
We verify that, if $G$ contains two multiplicatively independent elements $g_1,g_2$, then a two-block $(G,\mathcal S)$-PTE set solution exists in the subgroup $\langle g_1,g_2\rangle$ generated by $g_1,g_2$.
Indeed, write $\mathcal S=\{\chi_1,\ldots,\chi_s\}$ and, for each nonnegative integer $\mathsf k$, put
\[
\Omega_{\mathsf k}:=\{g_1^i g_2^j:0\leq i,j\leq\mathsf k\}.
\]
Multiplicative independence gives $|\Omega_{\mathsf k}|=(\mathsf k+1)^2$.
Applying \cref{lem:number-field-character-sum-tuples} with $I=\{1,\ldots,s\}$, $b_j=\chi_j(g_1)$, and $b'_j=\chi_j(g_2)$ and proceeding as in the proof of \cref{thm:infinite-existence} gives the required set solution.

Now choose generators $u,v$ of $H$.
For any finite set $E\subseteq R^\times$, choose a positive integer $\mathsf M$ larger than the absolute values of the exponents in every element of the form $\pm u^{k_1}v^{k_2}$ belonging to $E$, and put
\[
G_{\mathsf M}:=\langle u^{\mathsf M},v^{\mathsf M}\rangle.
\]
Multiplicative independence also excludes $u^{k_1}v^{k_2}=-1$, so
\[
E\cap\bigl(G_{\mathsf M}\cup(-G_{\mathsf M})\bigr)\subseteq\{1,-1\},
\qquad
G_{\mathsf M}\cap(-G_{\mathsf M})=\varnothing.
\]
Apply the preceding observation with $G=G_{\mathsf M}$ and $\mathcal S=\{\mathbf1,\chi\}$, where $\chi:G_{\mathsf M}\hookrightarrow\mathbb K^\times$ is the inclusion, and take a two-block $(G_{\mathsf M},\mathcal S)$-PTE set solution $(A_1,A_2)$.
For $w\in A_1$, the two parts of
\[
C:=w^{-1}A_2\cup\bigl(-w^{-1}(A_1\setminus\{w\})\bigr)
\]
are disjoint, contain neither $1$ nor $-1$, and have total sum $1$.
Thus $C\cap E=\varnothing$.

By assumption, every $x\in R$ is an integer linear combination of units.
Since $-\varepsilon$ is also a unit, there are a nonempty finite set $J\subseteq\mathbb N$ and units $\varepsilon_i\in R^\times$ such that, with repetitions allowed,
\[
x=\sum_{i\in J}\varepsilon_i.
\]
For zero, use $0=1+(-1)$.
For each $i\in J$, use the preceding construction to choose successively a finite set $C_i\subseteq R^\times$ whose elements sum to $1$ and such that
\[
C_i\cap\{\varepsilon_i^{-1}\varepsilon_j z:j\in J,\ j<i,\ z\in C_j\}=\varnothing.
\]
Then $\bigsqcup_{i\in J}\varepsilon_i C_i$ is a finite set of distinct units whose elements sum to
\[
\sum_{i\in J}\sum_{z\in C_i}\varepsilon_i z
=\sum_{i\in J}\varepsilon_i
=x.
\]
\end{proof}

\begin{corollary}[Thuswaldner--Ziegler's Conjecture~7.3]\label{cor:thuswaldner-ziegler-conjecture}
For $\varrho$ as above, every element of $\mathbb Z[\varrho]$ is a finite sum of distinct units of this ring.
\end{corollary}

\begin{proof}
By Thuswaldner--Ziegler~\cite[Lemma~7.1(3), (4)]{ThuswaldnerZiegler2011}, the elements $\varrho$ and $-1-\varrho^{-1}$ are multiplicatively independent units of $\mathbb Z[\varrho]$.
Moreover, the nonnegative integer powers of $\varrho$ are units and generate $\mathbb Z[\varrho]$ as a $\mathbb Z$-module.
Thus \cref{prop:distinct-unit-sums} applies.
\end{proof}

In \cref{sec:rational-design-existence}, we use $(G,\mathcal S)$-PTE set solutions to remove multiplicities from integer-weighted rational designs and construct equal-weight rational designs.


\subsection{Ideal solutions on the unit circle and finite cyclic groups} \label{subsec:pte-circle-ideal}

We first characterize ideal solutions on subgroups of the unit circle for the character family $\{z\mapsto z^k:0\leq k\leq t\}$ with $t\geq1$ and deduce a criterion for norm-one groups of imaginary quadratic fields (\cref{thm:circle-minimum,cor:imaginary-quadratic-minimum}).
We then consider finite cyclic groups.
For a cyclic group $G$ of prime order, we prove that no $(G,\mathcal S)$-PTE solution exists over a field of characteristic zero if $\mathcal S$ contains both the trivial character and a nontrivial one (\cref{prop:prime-cyclic-rigidity}).
For order $2\mathsf p$, with $\mathsf p$ an odd prime, we show that adjoining $\chi^2$ to $\{\mathbf 1,\chi\}$, where $\chi$ is a faithful character, raises the minimum size of a two-block solution from $2$ to $\mathsf p$ (\cref{prop:cyclic-jump}).
For each integer $k\geq1$, put
\begin{equation}\label{eq:roots-of-unity}
R_k:=\{z\in\mathbb C^\times:z^k=1\}.
\end{equation}
For a multiplicative subgroup $G\leq\mathbb C^\times$ and an integer $t\geq0$, define the power-character family
\begin{equation}\label{eq:power-character-family}
\mathcal T_t(G):=\{z\mapsto z^k:0\leq k\leq t\}\subseteq\operatorname{Hom}(G,\mathbb C^\times).
\end{equation}

\begin{theorem}\label{thm:circle-minimum}
Let $t,q$ be integers with $t\geq1$ and $q\geq2$, and put $n:=t+1$.
Let $G\leq\{z\in\mathbb C^\times:|z|=1\}$ and put $\mathcal S:=\mathcal T_t(G)$.
Then the following are equivalent.
\begin{enumerate}
\item[(i)] $\kappa_{G,\mathcal S}(q)=\kappa^{\mathrm{set}}_{G,\mathcal S}(q)=n$.
\item[(ii)] $R_n\subseteq G$ and $[G:R_n]\geq q$.
\end{enumerate}
Every solution of size $n$ consists of $q$ pairwise distinct cosets of $R_n$.
\end{theorem}
\begin{proof}
By \cref{prop:PTEbound-newton}, every solution has size at least $n$.
Write $\tau_k(A)$ for the $k$-th elementary symmetric function of the entries of $A$.
Let $(A_1,\ldots,A_q)$ be a solution of size $n$.
Newton's identities \cite[Chapter~I, Equation~(2.11')]{Macdonald1995} show that $\tau_k(A_i)$ is independent of $i$ for $1\leq k<n$.
For each $i$, write the entries of $A_i$ as $z_{i,1},\ldots,z_{i,n}$, counting multiplicities.
Since $|z_{i,j}|=1$, for $1\leq k<n$ we have
\[
\tau_n(A_i)\overline{\tau_k(A_i)}
=
\sum_{\substack{J\subseteq\{1,\ldots,n\}\\|J|=k}}
\prod_{j\in\{1,\ldots,n\}\setminus J}z_{i,j}
=
\tau_{n-k}(A_i).
\]
If the common value of some $\tau_k(A_i)$ with $1\leq k<n$ were nonzero, then
\[
\tau_n(A_i)=\frac{\tau_{n-k}(A_i)}{\overline{\tau_k(A_i)}}
\]
would also be independent of $i$.
Hence $A_1=A_2$, a contradiction.
Thus $\tau_k(A_i)=0$ for $1\leq k<n$.
For each $i$, choose ${v}\in\operatorname{supp}(A_i)$.
Then
\[
\prod_{z\in A_i}(X-z)=X^n-{v}^n.
\]
Since ${v}\ne0$, this polynomial has only simple roots, so $A_i={v}R_n$.
As ${v}\in G$ and $A_i\subseteq G$, it follows that $R_n={v}^{-1}A_i\subseteq G$.
The blocks give $q$ distinct cosets, and hence $[G:R_n]\geq q$.
Conversely, if $R_n\subseteq G$ and $[G:R_n]\geq q$, choose $q$ distinct cosets $A_1,\ldots,A_q$ of $R_n$.
For each $i$, write $A_i={v} R_n$ with ${v}\in G$. For $1\leq k<n$,
\[
\sum_{z\in {v} R_n}z^k
={v}^k\sum_{z'\in R_n}z'^k=0.
\]
They therefore form a set solution of size $n$, and the lower bound gives both asserted minimum values.
\end{proof}
For an imaginary quadratic subfield $\mathbb K\subseteq\mathbb C$, put
\[
R_{\mathbb K}:=\{z\in \mathbb K^\times:z^k=1\text{ for some }k\in\mathbb Z_{\geq1}\}.
\]
We have
\begin{equation}\label{eq:imaginary-quadratic-roots}
R_{\mathbb K}=
\begin{cases}
R_4,&\mathbb K=\mathbb Q(\mathrm i),\\
R_6,&\mathbb K=\mathbb Q(\sqrt{-3}),\\
R_2,&\text{otherwise},
\end{cases}
\end{equation}
see \cite[Proposition~13.1.5]{IrelandRosen1990}.
\Needspace{5\baselineskip}
Together with \cref{lem:norm-one-condition-D,thm:circle-minimum}, this gives the following corollary.

\begin{corollary}\label{cor:imaginary-quadratic-minimum}
Let $\mathbb K\subseteq\mathbb C$ be an imaginary quadratic field, and let $t,q$ be integers with $t\geq1$ and $q\geq2$.
Put $\mathcal S:=\mathcal T_t(T_{\mathbb K})$.
An ideal $(T_{\mathbb K},\mathcal S)$-PTE solution with $q$ blocks exists exactly in the following cases:
\[
\begin{cases}
t\in\{1,3\},&\mathbb K=\mathbb Q(\mathrm i),\\
t\in\{1,2,5\},&\mathbb K=\mathbb Q(\sqrt{-3}),\\
t=1,&\text{otherwise}.
\end{cases}
\]
\end{corollary}

\begin{proof}
Put $n:=t+1$.
By \cref{lem:norm-one-condition-D}, $T_{\mathbb K}$ is infinite, so $[T_{\mathbb K}:R_n]\geq q$ whenever $R_n\subseteq T_{\mathbb K}$.
Moreover, $R_n\subseteq T_{\mathbb K}$ if and only if $R_n\subseteq R_{\mathbb K}$, since every root of unity in $\mathbb K$ has norm one.
By \eqref{eq:imaginary-quadratic-roots}, this is equivalent to $n\mid4$, $n\mid6$, or $n\mid2$ in the three respective cases.
The result now follows from \cref{thm:circle-minimum}.
\end{proof}

Let $\mathbb K$ be a field of characteristic zero.
For each integer $k\geq1$, let $\Phi_k\in\mathbb Z[X]$ denote the $k$-th cyclotomic polynomial; it is irreducible over $\mathbb Q$ \cite[Chapter~13, Section~2, Theorem~1]{IrelandRosen1990}.
For each integer $m\geq2$, we use $0,\ldots,m-1$ as representatives for the elements of $\mathbb Z/m\mathbb Z$.

\begin{proposition}\label{prop:prime-cyclic-rigidity}
Let ${\mathsf{p}}$ be prime, put $G:=\mathbb Z/{\mathsf{p}}\mathbb Z$, and let $\mathcal S\subseteq\operatorname{Hom}(G,\mathbb K^\times)$ be finite, contain $\mathbf1$, and satisfy $\bigcap_{\chi\in\mathcal S}\ker\chi=\{0\}$.
Then $\kappa_{G,\mathcal S}(q)=\kappa^{\mathrm{set}}_{G,\mathcal S}(q)=\infty$ for every integer $q\geq2$.
\end{proposition}

\begin{proof}
Suppose that a solution exists, and let $A_1$ and $A_2$ be its first two blocks.
The condition $\bigcap_{\chi\in\mathcal S}\ker\chi=\{0\}$ ensures that there is a nontrivial character $\chi\in\mathcal S$.
Since $G$ has prime order, $\chi$ is injective, and hence $e:=\chi(1)$ is a primitive $\mathsf p$-th root of unity.

Under this assumption, we show that $A_1=A_2$.
This would contradict the defining conditions of a $(G,\mathcal S)$-PTE solution.
To this end, introduce the polynomial
\[
\phi(X):=\sum_{j=0}^{\mathsf p-1}
\bigl(\xi_{A_1}(j)-\xi_{A_2}(j)\bigr)X^j\in\mathbb Z[X].
\]
Recall that $\xi_A(j)$ is the multiplicity of $j$ in the multiset $A$.

By \eqref{eq:multiple-pte-character-sums}, we have $\phi(1)=\phi(e)=0$.
Since $\deg\phi\leq\mathsf p-1=\deg\Phi_{\mathsf p}$, the assumption that $\phi$ is not a rational multiple of $\Phi_{\mathsf p}$ would contradict the fact that $\Phi_{\mathsf p}$ is the minimal polynomial of $e$ over $\mathbb Q$.
Thus $\phi=c\Phi_{\mathsf p}$ for some $c\in\mathbb Q$.

Substituting $X=1$ on both sides gives $c=0$, since $\Phi_{\mathsf p}(1)=\mathsf p\neq0$.
Hence $\phi\equiv0$, so $\xi_{A_1}(j)=\xi_{A_2}(j)$ for every $j$.
Therefore $A_1=A_2$, which completes the proof.
\end{proof}

\begin{proposition}\label{prop:cyclic-jump}
Let ${\mathsf{p}}$ be an odd prime and put $G:=\mathbb Z/(2{\mathsf{p}})\mathbb Z$, ${e}:=\exp(2\pi\mathrm i/(2{\mathsf{p}}))$, and $\mathbb K:=\mathbb Q({e})$.
Define $\chi(j):={e}^j$, $\mathcal T_1:=\{\mathbf1,\chi\}$, and $\mathcal T_2:=\{\mathbf1,\chi,\chi^2\}$.
Then
\[
\kappa_{G,\mathcal T_1}(2)=\kappa^{\mathrm{set}}_{G,\mathcal T_1}(2)=2,
\qquad
\kappa_{G,\mathcal T_2}(2)=\kappa^{\mathrm{set}}_{G,\mathcal T_2}(2)={\mathsf{p}},
\]
and, up to interchanging the two blocks, the even and odd residue classes give the unique $\mathcal T_2$-solution of size ${\mathsf{p}}$.
\end{proposition}

\begin{proof}
First we show that $\kappa^{\mathrm{set}}_{G,\mathcal T_1}(2)=2$.
Indeed, $(\{0,\mathsf p\},\{1,\mathsf p+1\})$ is a $(G,\mathcal T_1)$-PTE set solution of size $2$.
Moreover, the injectivity of $\chi$ gives $\kappa^{\mathrm{set}}_{G,\mathcal T_1}(2)>1$.

Next we show that $\kappa^{\mathrm{set}}_{G,\mathcal T_2}(2)=\mathsf p$.
We first construct a $(G,\mathcal T_2)$-PTE set solution of size $\mathsf p$.
Put
\[
A_1:=\{0,2,\ldots,2\mathsf p-2\},
\qquad
A_2:=\{1,3,\ldots,2\mathsf p-1\},
\]
and, to verify that the character sums agree, define
\[
\phi_{A_1,A_2}(X):=\sum_{j\in A_1}X^j-\sum_{j\in A_2}X^j
=\sum_{j=0}^{2\mathsf p-1}(-1)^jX^j
=(1-X)(1+X^2+\cdots+X^{2\mathsf p-2}).
\]
Since this polynomial has $1,e,e^2$ as roots, $(A_1,A_2)$ is a $(G,\mathcal T_2)$-PTE set solution of size $\mathsf p$.

It therefore suffices to show that $\kappa^{\mathrm{set}}_{G,\mathcal T_2}(2)\geq\mathsf p$.
Let $(A'_1,A'_2)$ be a $(G,\mathcal T_2)$-PTE set solution.
To bound its size from below, introduce the polynomial
\[
\phi(X):=\sum_{j=0}^{2\mathsf p-1}
\bigl(\xi_{A'_1}(j)-\xi_{A'_2}(j)\bigr)X^j\in\mathbb Z[X].
\]
Here $\xi_A(j)$ denotes the multiplicity of $j$ in the multiset $A$.
Applying \eqref{eq:multiple-pte-character-sums} to the characters in $\mathcal T_2$ gives $\phi(1)=\phi(e)=\phi(e^2)=0$.

The minimal polynomials of $1,e,e^2$ over $\mathbb Q$ are $X-1,\Phi_{2\mathsf p},\Phi_{\mathsf p}$, respectively.
Since these are distinct irreducible polynomials and $\phi(1)=\phi(e)=\phi(e^2)=0$, we have
\[
(X-1)\Phi_{2\mathsf p}(X)\Phi_{\mathsf p}(X)\mid\phi(X).
\]
Since $\mathsf p$ is an odd prime, $\Phi_{2\mathsf p}(X)=\Phi_{\mathsf p}(-X)$, and hence $\deg\bigl((X-1)\Phi_{2\mathsf p}(X)\Phi_{\mathsf p}(X)\bigr)=2\mathsf p-1$.
The polynomial $\phi_{A_1,A_2}$ also has $1,e,e^2$ as roots and has degree $2\mathsf p-1$.
It follows that
\[
(X-1)\Phi_{2\mathsf p}(X)\Phi_{\mathsf p}(X)\mid\phi_{A_1,A_2}(X).
\]
Since $\deg\phi\leq2\mathsf p-1=\deg\phi_{A_1,A_2}$, there is a rational number $c$ such that $\phi=c\phi_{A_1,A_2}$.

Substituting $X=0$ in $\phi(X)=c\phi_{A_1,A_2}(X)$ and using $\phi_{A_1,A_2}(0)=1$, we obtain $c=\xi_{A'_1}(0)-\xi_{A'_2}(0)$.
Comparing the coefficients of $X^j$ on both sides, we obtain
\[
\xi_{A'_1}(j)-\xi_{A'_2}(j)
=(-1)^j\bigl(\xi_{A'_1}(0)-\xi_{A'_2}(0)\bigr)
\qquad(0\leq j\leq2\mathsf p-1).
\]
If $\xi_{A'_1}(0)-\xi_{A'_2}(0)=0$, then all multiplicities agree and $A'_1=A'_2$, contrary to the fact that the two blocks are nonempty and disjoint.

Since the two blocks are disjoint, $\lvert\xi_{A'_1}(j)-\xi_{A'_2}(j)\rvert=\xi_{A'_1}(j)+\xi_{A'_2}(j)$ for each $j$.
Hence
\[
|A'_1|=|A'_2|
=\frac12\sum_{j=0}^{2\mathsf p-1}
\bigl|\xi_{A'_1}(j)-\xi_{A'_2}(j)\bigr|
=\mathsf p\,\bigl|\xi_{A'_1}(0)-\xi_{A'_2}(0)\bigr|
\geq\mathsf p.
\]

Next we prove uniqueness for solutions of size $\mathsf p$.
The preceding argument applies also when the blocks are multisets.
Thus, for a solution of size $\mathsf p$, we have $\bigl|\xi_{A'_1}(0)-\xi_{A'_2}(0)\bigr|=1$.
Without loss of generality, we may take $A'_1,A'_2$ so that $\xi_{A'_1}(0)-\xi_{A'_2}(0)=1$.
Then for each $j$ we have
\[
\xi_{A'_1}(j)-\xi_{A'_2}(j)
=(-1)^j\bigl(\xi_{A'_1}(0)-\xi_{A'_2}(0)\bigr)
=(-1)^j.
\]
Since the two blocks are disjoint, it follows that
\[
\bigl(\xi_{A'_1}(j),\xi_{A'_2}(j)\bigr)
=
\begin{cases}
(1,0),&\text{if }j\equiv0\pmod{2},\\
(0,1),&\text{if }j\equiv1\pmod{2}.
\end{cases}
\]
Thus $(A'_1,A'_2)=(A_1,A_2)$.

The same arguments give $\kappa_{G,\mathcal T_1}(2)=2$ and $\kappa_{G,\mathcal T_2}(2)=\mathsf p$ when the blocks are multisets.
\end{proof}

\subsection{Classical ideal PTE solutions from character sums} \label{subsec:pte-block-transfer}

We first establish a criterion for transferring character sums to classical power sums (\cref{thm:ideal-block-transfer}).
We then use it to construct a nonsymmetric ideal family of degree five (\cref{thm:new-degree-five-family}) and show that, outside a finite set of parameters, its solutions are inequivalent to the two Choudhry constructions specified in \cref{prop:degree-five-separation}.
The degree-seven examples in \cref{ex:degree-seven-examples} give two further applications.

For finite multisets on $\mathbb C$, put
\begin{equation}\label{eq:Pk}
\psi_k:\operatorname{MSet}(\mathbb C)\longrightarrow\mathbb C \quad(k\in\mathbb Z_{\geq0}), \qquad \psi_0(A):=|A|,\qquad \psi_k(A):=\sum_{z\in A}z^k\quad(k\geq1).
\end{equation}
Let $G$ be an abelian group and let $\upsilon:G\to\operatorname{MSet}(\mathbb Q)$.
To replace each occurrence of $g$ in a finite multiset $A$ on $G$ by the multiset $\upsilon(g)$, define $\upsilon_*A$ by
\begin{equation}\label{eq:block-transfer-multiplicity}
\xi_{\upsilon_*A}(z) := \sum_{g\in\operatorname{supp}(A)} \xi_A(g)\xi_{\upsilon(g)}(z) \qquad (z\in\mathbb Q).
\end{equation}

\begin{definition}\label{def:classical-multiple-pte}
Let $q\geq2$, $n\geq1$, and $t\geq0$ be integers.
For finite multisets $A_1,\ldots,A_q$ on $\mathbb{Q}$, write
\[
[A_1]={}_t^n\cdots={}_t^n[A_q]
\]
if
\[
|A_i|=n\quad(1\leq i\leq q),\qquad
\psi_k(A_i)=\psi_k(A_1)
\quad(1\leq i\leq q,\ 1\leq k\leq t).
\]
The classical multiple PTE problem asks for such relations with integer entries and pairwise disjoint supports.
For $q=2$, $t\geq1$, and $n=t+1$, a rational relation with pairwise disjoint supports is called an ideal rational PTE solution of degree $t$.
\end{definition}

\begin{remark} \label{rem:multiple-PTE}
Equal-power-sum constructions go back to Prouhet \cite{Prouhet1851} and were studied further by Jacroux \cite{Jacroux1995}; the multiple PTE problem in \cref{def:classical-multiple-pte} is the one-dimensional case of \cite[Problem~7.2]{InagakiMatsumuraSawaUchida2026}.
\end{remark}

\begin{proposition}
\label{thm:ideal-block-transfer}
Let $t,q$ be integers with $t\geq1$ and $q\geq2$.
Let $\mathbb K,G,\mathcal S$ be as in \cref{def:abelian-pte}, and let $\upsilon:G\to\operatorname{MSet}(\mathbb Q)$.
Suppose that $(A_1,\ldots,A_q)$ is a $(G,\mathcal S)$-PTE solution and that $\psi_k\circ\upsilon\in\operatorname{span}_{\mathbb K}\mathcal S\subseteq\mathbb K^G$ for $0\leq k\leq t$.
Then the following hold.
\begin{enumerate}
\item[\rm(a)] The transferred multisets have equal power sums through degree $t$:
\begin{equation}
 \psi_k(\upsilon_*A_1)
 =\cdots=
 \psi_k(\upsilon_*A_q)
 \qquad(0\leq k\leq t).
 \label{eq:block-transfer}
\end{equation}
\item[\rm(b)] If \(q=2\) and $|\upsilon_*A_1|=|\upsilon_*A_2|=t+1$, then $\upsilon_*A_1\neq\upsilon_*A_2$ if and only if
\begin{equation}
 \prod_{g\in A_1}\left(\prod_{z\in\upsilon(g)}z\right)
 \neq
 \prod_{g\in A_2}\left(\prod_{z\in\upsilon(g)}z\right).
 \label{eq:block-product-test}
\end{equation}
\item[\rm(c)] Under the hypotheses of \textup{(b)}, if \eqref{eq:block-product-test} holds, then $\operatorname{supp}(\upsilon_*A_1)\cap\operatorname{supp}(\upsilon_*A_2)=\varnothing$.
For $c\in\mathbb Q$, choose $m\in\mathbb Z_{>0}$ such that $m(z+c)\in\mathbb Z$ for every $z\in\operatorname{supp}(\upsilon_*A_1)\cup\operatorname{supp}(\upsilon_*A_2)$, and put $\alpha(z):=m(z+c)$ for $z\in\mathbb Q$.
Then $(\alpha_*(\upsilon_*A_1),\alpha_*(\upsilon_*A_2))$ is an integer ideal solution of degree $t$.
\end{enumerate}
\end{proposition}

\begin{proof}
For \textup{(a)}, \eqref{eq:block-transfer-multiplicity} gives, for $1\leq i\leq q$ and $0\leq k\leq t$,
\[
\psi_k(\upsilon_*A_i)=\sum_{g\in A_i}\psi_k(\upsilon(g)).
\]
This identity and $\psi_k\circ\upsilon\in\operatorname{span}_{\mathbb K}\mathcal S$ give \eqref{eq:block-transfer}.

For \textup{(b)}, the definitions give
\[
\prod_{z\in\upsilon_*A_i}z=\prod_{g\in A_i}\left(\prod_{z\in\upsilon(g)}z\right)
\qquad(i=1,2).
\]
Newton's identities \cite[Chapter~I, Equation~(2.11')]{Macdonald1995} show that the polynomials $\prod_{z\in\upsilon_*A_i}(X-z)$, $i=1,2$, agree except possibly in their constant terms, whose difference is
\[
(-1)^{t+1}\left(\prod_{g\in A_1}\left(\prod_{z\in\upsilon(g)}z\right)-\prod_{g\in A_2}\left(\prod_{z\in\upsilon(g)}z\right)\right).
\]
This proves the equivalence.

For \textup{(c)}, the polynomials in \textup{(b)} differ by a nonzero constant, so they have no common root.
Thus the transferred multisets have disjoint supports.
The binomial theorem gives $\psi_k(\alpha_*(\upsilon_*A_1))=\psi_k(\alpha_*(\upsilon_*A_2))$ for $0\leq k\leq t$.
Since $\alpha$ is injective and takes every entry to an integer, the image blocks form an integer ideal solution of degree $t$.
\end{proof}

For a finite set $I\subseteq\mathbb N$ and $\mathbf k=(k_i)_{i\in I}\in\mathbb Z_{\geq0}^I$, write
\[
 |\mathbf k|:=\sum_{i\in I}k_i,\qquad
 \binom{k}{\mathbf k}:=\frac{k!}{\prod_{i\in I}k_i!}\quad(k=|\mathbf k|),
\]
and $\mathbf a^{\mathbf k}:=\prod_{i\in I}a_i^{k_i}$.
For a field $\mathbb K$, an abelian group $G$, and characters $\chi_i:G\to\mathbb K^\times$ indexed by $I$, define the character
\[
 \chi_{\mathbf k}(g):=\prod_{i\in I}\chi_i(g)^{k_i}.
\]

\begin{lemma}\label{lem:block-character-support}
Let $\mathbb K$ be a field of characteristic zero, let $G$ be an abelian group, and let $I\subseteq\mathbb N$ be finite.
Let $\chi_i:G\to\mathbb K^\times$ be characters indexed by $i\in I$.
Let $C\in\operatorname{MSet}(\mathbb K^I)$ satisfy $|C|>0$ and
\[
 \sum_{i\in I}a_i\chi_i(g)\in\mathbb Q
 \qquad\bigl(\mathbf a=(a_i)_{i\in I}\in\operatorname{supp}(C),\ g\in G\bigr).
\]
Define $\upsilon_C:G\to\operatorname{MSet}(\mathbb Q)$ by
\[
\xi_{\upsilon_C(g)}(z)
:=
\sum_{\substack{\mathbf a\in\operatorname{supp}(C)\\
\sum_{i\in I}a_i\chi_i(g)=z}}
\xi_C(\mathbf a)
\qquad(g\in G,\ z\in\mathbb Q).
\]
Then, for every integer $k\geq0$ and every $g\in G$,
\begin{equation}
 \psi_k(\upsilon_C(g))=\sum_{\substack{\mathbf k\in\mathbb Z_{\geq0}^I\\|\mathbf k|=k}}\binom{k}{\mathbf k}
 \left(\sum_{\mathbf a\in C}\mathbf a^{\mathbf k}\right)\chi_{\mathbf k}(g).
 \label{eq:block-character-support}
\end{equation}
\end{lemma}

\begin{proof}
Expand each $\left(\sum_{i\in I} a_i\chi_i(g)\right)^k$ by the multinomial theorem and sum over $\mathbf a\in C$.
\end{proof}

\begin{example}\label{ex:block-field-trace}
Let $\mathbb K/\mathbb Q$ be a finite separable extension, and let $\widetilde{\mathbb K}$ be its normal closure over $\mathbb Q$ inside a fixed algebraic closure.
Let $\sigma_i:\mathbb K\hookrightarrow\widetilde{\mathbb K}$, $1\leq i\leq[\mathbb K:\mathbb Q]$, be all the distinct $\mathbb Q$-embeddings.
Let $G\leq\mathbb K^\times$, and let $\varnothing\ne C'\subseteq\mathbb K$ be finite.
In \cref{lem:block-character-support}, take the coefficient field to be $\widetilde{\mathbb K}$ and the index set to be $I=\{1,\ldots,[\mathbb K:\mathbb Q]\}$, and put
\[
 \chi_i=\sigma_i|_G:G\to\widetilde{\mathbb K}^{\times}\quad(i\in I),\qquad
 C=\{(\sigma_i(z))_{i\in I}:z\in C'\}\subseteq\widetilde{\mathbb K}^{I}.
\]
Then
\[
\xi_{\upsilon_C(g)}(r)
=
\left|\left\{z\in C':
\operatorname{Tr}_{\mathbb K/\mathbb Q}(zg)=r
\right\}\right|
\qquad(g\in G,\ r\in\mathbb Q).
\]
Thus \cref{lem:block-character-support} expresses their power sums as linear combinations of products of the characters $\sigma_i|_G$.
\end{example}

For the degree-five construction, we apply \cref{thm:ideal-block-transfer} to two blocks on $\{\pm1\}^2$, replacing each element by three rational numbers.

\begin{theorem} \label{thm:new-degree-five-family}
Let $u\in\mathbb{Q}^\times$, and define
\[
 y_1:=u-u^{-1},\qquad y_2:=u+u^{-1},\qquad y_2^2=y_1^2+4.
\]
Define $h_i,{a'_{i,1}},{a'_{i,2}}\in\mathbb Q$ $(1\leq i\leq3)$ by {\small
\begin{equation}
\begin{array}{@{}ll@{}}
 h_1=-y_2(y_1+4)^2(y_1^2+2),
 &h_2=-2y_2(y_1+1)^2(y_1^2+8),\qquad h_3=0,\\
 {a'_{1,1}}=y_1(y_1-2)(y_1+1)(y_1^2+2y_1+4),
 &{a'_{1,2}}=(y_1+1)(3y_1^2-2y_1+16)(y_1^2+2y_1+4),\\
 {a'_{2,1}}=-(y_1+4)(2y_1^2-y_1+6)(y_1^2+2y_1+4),
 &{a'_{2,2}}=(y_1-2)(y_1+4)(y_1^2+2y_1+4),\\
 {a'_{3,1}}=y_1(y_1-2)(2y_1^3+3y_1^2+12y_1+20),
 &{a'_{3,2}}=(y_1-2)(5y_1^3+12y_1^2+12y_1+32).
\end{array}
\label{eq:degree-five-data-concise}
\end{equation}
}
and set
\begin{equation}
A_1(u)=\{h_1\pm {a'_{1,1}},h_2\pm {a'_{2,1}},\pm {a'_{3,1}}\}, \qquad A_2(u)=\{h_1\pm {a'_{1,2}},h_2\pm {a'_{2,2}},\pm {a'_{3,2}}\}.
\label{eq:degree-five-family-concise}
\end{equation}
The braces denote six-element multisets, with each indicated occurrence counted once.
\begin{enumerate}
\item[\rm(a)] If \(u\neq\pm1\), then \([A_1(u)]={}_5^6[A_2(u)]\) is a nonsymmetric rational ideal solution.
\item[\rm(b)] There exists a sequence $(u_i)_{i\geq1}\in(\mathbb Z\setminus\{0,\pm1\})^{\mathbb N}$ of distinct parameters such that the blocks in each relation
\[
 [u_i^5A_1(u_i)]={}_5^6[u_i^5A_2(u_i)]
\]
are disjoint six-element sets of integers, and the resulting nonsymmetric ideal solutions are pairwise affinely inequivalent over $\mathbb Q$.
\end{enumerate}
\end{theorem}

\begin{proof}
See Appendix~\ref{app:degree-five-family}.
\end{proof}

Borwein's degree-five family is symmetric \cite[Eq.~(7)]{MishimaLuSawaUchida2026}.
We compare our nonsymmetric family with Choudhry's constructions in \cite[Eq.~(3.8)]{Choudhry2003} and \cite[Eqs.~(3.27)--(3.30)]{Choudhry2017}.
The latter gives equal power sums through degree four, and also through degree five when its parameters satisfy the condition in \cite[Eq.~(4.28)]{Choudhry2017}.

\begin{proposition}
\label{prop:degree-five-separation}
Except for finitely many $u\in\mathbb{Q}^\times \setminus \{\pm1\}$, the solution in \eqref{eq:degree-five-family-concise} is not equivalent over $\mathbb{Q}$ to either of the following:
\begin{enumerate}
\item[(i)] A member of the family in \cite[Eq.~(3.8)]{Choudhry2003} obtained at a rational parameter for which the formula is defined.
\item[(ii)] A pair of six-element multisets obtained from \cite[Eq.~(3.30)]{Choudhry2017} with all twelve entries rational.
\end{enumerate}
In \textup{(ii)}, we allow any parameter values in \cite[Eqs.~(3.27)--(3.30)]{Choudhry2017} for which all twelve entries are rational.
\end{proposition}

\begin{proof}
See Appendix~\ref{app:degree-five-separation}.
\end{proof}

\begin{example} \label{ex:degree-seven-examples}
Let $\mathbb{K} = \mathbb{Q}$, $G = \mathbb{Q}_{>0}^{\times}$, and $\mathcal{S} = \{\mathbf{1},g\mapsto g^2,g\mapsto g^4,g\mapsto g^6\}$.
For $\upsilon(g) = \{g,-g\}$, the moments satisfy, for $g\in G$,
\[
 \psi_0(\upsilon(g))=2,\qquad
 \psi_{2j}(\upsilon(g))=2g^{2j}\quad(1\leq j\leq3),\qquad
 \psi_{2j+1}(\upsilon(g))=0\quad(0\leq j\leq3).
\]
The corresponding reduction of degree seven to the even degrees $2,4,6$ is recorded in \cite[Eqs.~(2.1)--(2.2)]{Choudhry2022}.
Define
\[
\begin{aligned}
 A_{1,1}&=\{985,6611,7204,9993\},
 \qquad A_{1,2}=\{2243,5340,8129,9839\},\\
 A_{2,1}&=\{2242551,2371260367,2639169308,3841724555\},\\
 A_{2,2}&=\{1142625215,1575093099,3045557287,3779551972\}.
\end{aligned}
\]
For $i=1,2$, the sets $A_{i,1}$ and $A_{i,2}$ are disjoint, and the displayed integers satisfy
\[
 \psi_{2j}(A_{i,1})=\psi_{2j}(A_{i,2}) \qquad (1\leq j\leq3).
\]
Thus $(A_{i,1},A_{i,2})$ is a $(G,\mathcal{S})$-PTE solution of size $4$.
Moreover,
\[
\prod_{z\in A_{1,1}}z^2-\prod_{z\in A_{1,2}}z^2\equiv41\pmod{59}, \qquad \prod_{z\in A_{2,1}}z^2-\prod_{z\in A_{2,2}}z^2\equiv30\pmod{59}.
\]
Since $\prod_{z\in\upsilon(g)}z=-g^2$, these two nonzero congruences and \cref{thm:ideal-block-transfer} give the following two integer ideal solutions:
\[
 [\pm985,\pm6611,\pm7204,\pm9993]
 ={}_7^8
 [\pm2243,\pm5340,\pm8129,\pm9839]
\]
and
\[
\begin{aligned}
 &[\pm2242551,\pm2371260367,\pm2639169308,\pm3841724555]\\
 &\qquad={}_7^8
 [\pm1142625215,\pm1575093099,\pm3045557287,\pm3779551972],
\end{aligned}
\]
where $\pm z$ denotes the two occurrences $z$ and $-z$.
\end{example}

\endgroup

\section{From Integer-Weighted to Equal-Weight Rational Designs}\label{sec:rational-design-existence}
\begingroup
\setstretch{1.05}

This section uses pairwise disjoint sets with equal character sums to replace integer weights by equal weights on distinct rational points.
The resulting theorem is applied to spheres, orthogonal and unitary groups, and real and complex Grassmannians.
Recall that, for $Z\subseteq\R^d$, we write $Z(\Q):=Z\cap\Q^d$.

\begin{definition}\label{def:integer-weighted-design}
Let $Z\subseteq\R^d$ be equipped with a probability measure $\mu$, and let $\mathcal V$ be a finite-dimensional real vector space of $\mu$-integrable real-valued functions on $Z$ containing $\mathbf1$.
A nonempty finite set $\mathfrak X\subseteq Z(\Q)$ equipped with positive integer weights $m_{\mathbf x}$, $\mathbf x\in\mathfrak X$, is called an \emph{integer-weighted rational $\mathcal V$-design} if, with $m:=\sum_{\mathbf x\in\mathfrak X}m_{\mathbf x}$,
\begin{equation}\label{eq:integer-weighted-design}
\sum_{\mathbf x\in\mathfrak X}m_{\mathbf x}\varphi(\mathbf x)
=m\int_Z\varphi\,\mathrm d\mu
\qquad(\varphi\in\mathcal V).
\end{equation}
\end{definition}

The following consequence of Cui--Xia--Xiang~\cite[Lemmas~2.17, 3.4 and 3.7]{CuiXiaXiang2019} supplies the integer-weighted designs used in \cref{thm:integer-to-equal-weight}.

\begin{proposition}\label{prop:cxx-integer-weighted}
Let $d\geq1$ be an integer.
Let $Z\subseteq \R^d$ be compact and $\mu$ be a probability measure on $Z$.
Assume that every nonempty relatively open subset of $Z$ has positive measure.
Suppose that $Z_0 \subseteq Z(\Q)$ is dense in $Z$ and that
\[
\int_Z(f|_Z)\,\d\mu\in\Q
\qquad
\bigl(f\in\Q[X_1,\ldots,X_d]\bigr).
\]
Then, for every integer $t\geq0$, there exists an integer-weighted rational $\mathcal P_t(Z;\R)$-design supported on $Z_0$.
\end{proposition}
\begin{proof}
For $t=0$, any point of $Z_0$ with weight $1$ suffices; suppose henceforth that $t\geq1$.
In the proof of the cited Lemma~3.4, the approximating sequences may be chosen in $Z_0$ by density, and Lemma~2.17 applies because polynomial evaluations converge along these sequences.
The results of Cui, Xia, and Xiang~\cite[Lemmas~2.17, 3.4 and 3.7]{CuiXiaXiang2019} yield a rationally weighted $\mathcal P_t(Z;\R)$-design on a nonempty finite set $\mathfrak X\subseteq Z_0$, with positive rational weights $(m_{\mathbf{x}}^{\prime})_{\mathbf{x}\in\mathfrak X}$.
Write these weights as
\[
m_{\mathbf{x}}^{\prime}=\frac{m_{\mathbf{x}}}{\sum_{\mathbf{x}^{\prime}\in\mathfrak X}m_{\mathbf{x}^{\prime}}},
\qquad m_{\mathbf{x}}\in\Z_{>0}\quad(\mathbf{x}\in\mathfrak X).
\]
Substituting into \eqref{eq:weighted-design} gives \eqref{eq:integer-weighted-design}.
\end{proof}

\Needspace{5\baselineskip}

\subsection{Removal of Multiplicities Using $(G,\mathcal S)$-PTE Solutions}\label{subsec:remove-integer-weights}

\begin{lemma}\label{lem:nontrivial-character-integral}
Suppose that an abelian group \(G\) acts on \(Z\) by measurable bijections and preserves a probability measure \(\mu\).
Let \(\varphi:Z\to\C\) be integrable, and let \(\chi:G\to\C^\times\) be a character such that \(\varphi(g\cdot\mathbf{x})=\chi(g)\varphi(\mathbf{x})\) for all \(g\in G\) and \(\mathbf{x}\in Z\).
If \(\chi\neq\mathbf1\), then
\[
\int_Z\varphi\,\d\mu=0.
\]
\end{lemma}

\begin{proof}
Choose \(g\in G\) such that \(\chi(g)\neq1\).
Measure invariance gives
\[
\int_Z\varphi\,\d\mu
=
\int_Z\varphi(g\cdot\mathbf{x})\,\d\mu(\mathbf{x})
=
\chi(g)\int_Z\varphi\,\d\mu.
\]
Therefore, the integral vanishes.
\end{proof}

In what follows, let $\mathbb K\subseteq\C$ be a number field and $G$ an abelian group, and assume that condition $(\mathrm D)$ holds for $\nu$ and $(g_i)_{i\geq1}$.
Recall the coordinate maps $\nu_i$ from \eqref{eq:nu-coordinates}.
Real and imaginary parts of functions are taken pointwise.
The following theorem uses \cref{thm:infinite-existence} to replace integer multiplicities by distinct rational points of equal weight, while avoiding any prescribed finite set.

\begin{theorem}\label{thm:integer-to-equal-weight}
Let $d \geq 1$ and ${\mathsf{s}} \geq 0$ be integers, and suppose that the following conditions hold.
\begin{enumerate}
\item[\rm(a)] Let $Z \subseteq \R^d$ be compact.
The group $G$ acts on $Z$ by measurable bijections and preserves $Z(\Q)$ and a probability measure $\mu$.
A subset $W\subseteq Z$ is $G$-invariant, and we write $W(\Q)=W\cap Z(\Q)$.
\item[\rm(b)] There are integrable functions $\varphi_1,\ldots,\varphi_{\mathsf{s}}:Z\to\C$ and characters $\chi_1,\ldots,\chi_{\mathsf{s}}:G\to \mathbb K^\times$ such that
\[
\mathcal V=\operatorname{span}_{\R}\{\mathbf1,\Re\varphi_1,\Im\varphi_1,\ldots,\Re\varphi_{\mathsf{s}},\Im\varphi_{\mathsf{s}}\}
\]
and $\varphi_j(g\cdot\mathbf{x}) = \chi_j(g)\varphi_j(\mathbf{x})$ for $g \in G$, $\mathbf{x} \in Z$, and $1 \leq j \leq {\mathsf{s}}$.
\item[\rm(c)] There are a character $\chi_0:G\to \mathbb K^\times$, homomorphisms $\theta_i : \mathbb K^\times \to \Z$ indexed by $i\geq1$, and a map $\eta:W(\Q)\to \mathbb K^\times$ such that, for every $i\geq1$, the map $\theta_i\circ\chi_0$ is a positive integer multiple of $\nu_i$, and
\[
\eta(g\cdot\mathbf{x}) = \chi_0(g)\eta(\mathbf{x})
\qquad(g\in G,\ \mathbf{x}\in W(\Q)).
\]
For every $\mathbf{x} \in W(\Q)$, only finitely many $i$ satisfy $\theta_i(\eta(\mathbf{x})) \neq 0$.
\end{enumerate}
Then, for every integer-weighted rational $\mathcal V$-design supported on $W(\Q)$
and every finite set $E \subseteq Z(\Q)$, there exists a rational $\mathcal V$-design $\mathfrak X \subseteq W(\Q)$ disjoint from $E$.
\end{theorem}

\begin{proof}
Expand the integer weights into a finite sequence $(\mathbf{x}^{(1)},\ldots,\mathbf{x}^{(m)})$ of points in $W(\Q)$, allowing repetitions, where $m$ is the sum of the weights. Put
\[
I_0:=
\bigcup_{i=1}^{m}\{j\geq1:\theta_j(\eta(\mathbf{x}^{(i)}))\neq0\}
\cup
\bigcup_{\mathbf{x}^{\prime}\in E\cap W}\{j\geq1:\theta_j(\eta(\mathbf{x}^{\prime}))\neq0\}.
\]
This set is finite by \textup{(c)}.
Put
\[
G':=\langle g_j:j\geq1\rangle.
\]
Condition $(\mathrm D)$ implies that $\nu|_{G'}$ is injective.
If $g\in G'$ satisfies $\chi_0(g)=1$, then \textup{(c)} gives $\nu_j(g)=0$ for every $j$, so $g=1_G$.
Thus $\chi_0|_{G'}$ is injective, and the character set
\[
\mathcal S:=\{\mathbf1,\chi_0|_{G'},\chi_1|_{G'},\ldots,\chi_{\mathsf{s}}|_{G'}\}
\]
has common kernel $\{1_G\}$.
Apply \cref{thm:infinite-existence} to $G'$, with $\nu|_{G'}$, this character set, $q=\max\{2,m\}$, the index set $I_0$, and $E=\varnothing$ in that theorem.
It gives a finite set $I\subseteq\N\setminus I_0$ and a $q$-block $(G',\mathcal S)$-PTE set solution $(A_1,\ldots,A_q)$ of nonempty finite sets.
In particular, $I\ne\varnothing$: if $I$ were empty, all these pairwise disjoint nonempty sets would be contained in $\{1_G\}$, contradicting $q\geq2$.
We use the first $m$ blocks.
Every element of every $A_i$ is a product of positive integer powers of the $g_j$, $j\in I$, and the PTE equalities give
\[
|A_i|=|A_1|\qquad(1\leq i\leq m), \qquad \sum_{g\in A_i}\chi_j(g)=\sum_{g\in A_1}\chi_j(g) \qquad(1\leq i\leq m,\ 1\leq j\leq {\mathsf{s}}).
\]
Put
\[
\mathfrak X:=\{g\cdot\mathbf{x}^{(i)}:1\leq i\leq m,\ g\in A_i\}.
\]
Since $W$ is $G$-invariant and the action preserves
$Z(\Q)$, we have $\mathfrak X\subseteq W(\Q)$.

We next verify that all points in $\mathfrak X$ are distinct.
If $g\cdot\mathbf{x}^{(i)}=g'\cdot\mathbf{x}^{(i^{\prime})}$ with $g\in A_i$ and $g'\in A_{i'}$, then condition \textup{(c)} gives, for every $j\in I$,
\[
\theta_j(\chi_0(g))
 =\theta_j(\eta(g\cdot\mathbf{x}^{(i)}))-\theta_j(\eta(\mathbf{x}^{(i)}))
 =\theta_j(\eta(g'\cdot\mathbf{x}^{(i^{\prime})}))-\theta_j(\eta(\mathbf{x}^{(i^{\prime})}))
 =\theta_j(\chi_0(g')),
\]
because $\theta_j(\eta(\mathbf{x}^{(i)}))=\theta_j(\eta(\mathbf{x}^{(i^{\prime})}))=0$ for $j\in I$.
Since $\theta_j\circ\chi_0$ is a positive integer multiple of $\nu_j$, we have $\nu_j(g)=\nu_j(g')$ for every $j\in I$.
By condition $(\mathrm D)$, $\nu_j(g_j)>0$, and $\nu_j(g)$ is $\nu_j(g_j)$ times the exponent of $g_j$ in $g$; the same holds for $g'$.
Thus the exponents agree and $g=g'$.
As $g\in A_i\cap A_{i'}$, the disjointness of the sets $A_i$ gives $i=i'$.
Thus $\mathfrak X$ consists of $m|A_1|$ distinct points.

We next show that $\mathfrak X\cap E=\varnothing$.
If $\mathbf{x}^{\prime}=g\cdot\mathbf{x}^{(i)}\in E$, then, for every $j\in I$,
\[
 0=\theta_j(\eta(\mathbf{x}^{\prime}))-\theta_j(\eta(\mathbf{x}^{(i)}))=\theta_j(\chi_0(g))>0.
\]
Here the first equality uses $\mathbf{x}^{\prime}\in E\cap W$ and $I\cap I_0=\varnothing$, the second uses \textup{(c)}, and the strict inequality follows from \textup{(c)}, $(\mathrm D)$, and the positive exponents of $g$.
Since $I\ne\varnothing$, this is a contradiction.
Thus $\mathfrak X\cap E=\varnothing$.

Finally, we verify that $\mathfrak X$ is a $\mathcal V$-design.
For $1\leq j\leq {\mathsf{s}}$, applying \eqref{eq:integer-weighted-design} to $\Re\varphi_j$ and $\Im\varphi_j$, and using the PTE equalities above, gives
\begin{align*}
 \sum_{\mathbf{x}^{\prime}\in \mathfrak X}\varphi_j(\mathbf{x}^{\prime})
 &=\sum_{i=1}^{m}\sum_{g\in A_i}
   \chi_j(g)\varphi_j(\mathbf{x}^{(i)})\\
 &=\left(\sum_{g\in A_1}\chi_j(g)\right)\sum_{i=1}^{m}\varphi_j(\mathbf{x}^{(i)})\\
 &=m\left(\sum_{g\in A_1}\chi_j(g)\right)\int_Z\varphi_j\,\d\mu.
\end{align*}
If $\chi_j=\mathbf1$ on $G$, then $\sum_{g\in A_1}\chi_j(g)=|A_1|$.
If $\chi_j\neq\mathbf1$ on $G$, then $\int_Z\varphi_j\,\d\mu=0$ by \cref{lem:nontrivial-character-integral}.
Hence, in either case,
\[
 \sum_{\mathbf{x}^{\prime}\in \mathfrak X}\varphi_j(\mathbf{x}^{\prime})
 =m|A_1|\int_Z\varphi_j\,\d\mu.
\]
Taking real and imaginary parts and dividing by $|\mathfrak X|=m|A_1|$ gives \eqref{eq:equal-weight-design} for each $\Re\varphi_j$ and $\Im\varphi_j$.
Equation \eqref{eq:equal-weight-design} also holds for the constant function because $\mu(Z)=1$, and hence for every $\varphi\in \mathcal V$ by \textup{(b)}.
\end{proof}

The following is a specialization of \cref{thm:integer-to-equal-weight} to norm-one groups of imaginary quadratic fields.

\begin{corollary}\label{cor:norm-one-equal-weight}
Let $d\geq1$ and ${\mathsf{s}}\geq0$ be integers.
Let $\mathbb{K} \subseteq \C$ be an imaginary quadratic field and
\[
T_{\mathbb K} :=\{{v}\in \mathbb K^\times:{v}\overline {v}=1\}.
\]
Suppose that $T_{\mathbb K}$ acts on a compact set $Z\subseteq\R^d$ by measurable bijections preserving $Z(\Q)$, a probability measure $\mu$, and a subset $W\subseteq Z$.
Let $\varphi_1,\ldots,\varphi_{\mathsf{s}}:Z\to\C$ be integrable, let $k_1,\ldots,k_{\mathsf{s}}\in\Z$, let $\eta:W(\Q)\to\mathbb K^\times$, and let $k_0\in\Z_{>0}$.
Suppose that
\[
\begin{gathered}
\mathcal V=\operatorname{span}_{\R}\{\mathbf1,\Re\varphi_1,\Im\varphi_1,\ldots,\Re\varphi_{\mathsf{s}},\Im\varphi_{\mathsf{s}}\},\\
\varphi_j({v}\cdot\mathbf{x})={v}^{k_j}\varphi_j(\mathbf{x})
\quad({v}\in T_{\mathbb K},\ \mathbf{x}\in Z,\ 1\leq j\leq {\mathsf{s}}),\\
\eta({v}\cdot\mathbf{x})={v}^{k_0}\eta(\mathbf{x})
\quad({v}\in T_{\mathbb K},\ \mathbf{x}\in W(\Q)).
\end{gathered}
\]
%
Then, for every integer-weighted rational $\mathcal{V}$-design supported on $W(\Q)$ and every finite set $E\subseteq Z(\Q)$, there exists a rational $\mathcal V$-design $\mathfrak{X} \subseteq W(\Q)$ disjoint from $E$.
For every integer $k\geq1$, put
\[
T_{\mathbb K}^{(k)}:=\{{v}^k:{v}\in T_{\mathbb K}\}.
\]
The same conclusion holds for the restricted action of this subgroup.
\end{corollary}

\begin{proof}
Take $\mathfrak p_i,\nu,g_i$ as in \cref{lem:norm-one-condition-D}, and set $\theta_i:=\operatorname{ord}_{\mathfrak p_i}:\mathbb K^\times\to\Z$, so that $\nu_i=\theta_i|_{T_{\mathbb K}}$ for every $i\geq1$.
For every $\mathbf{x}\in W(\Q)$, only finitely many $i$ satisfy $\theta_i(\eta(\mathbf{x}))\ne0$, because a nonzero element of $\mathbb K$ has nonzero prime-ideal valuation at only finitely many prime ideals.
Also, $\theta_i({v}^{k_0})=k_0\nu_i({v})$ and $k_0>0$, so $\theta_i\circ\chi_0$ is a positive integer multiple of $\nu_i$ for $\chi_0({v})={v}^{k_0}$.
Thus the theorem applies with coefficient field \(\mathbb K\) and \(\chi_j({v})={v}^{k_j}\).
For the subgroup $T_{\mathbb K}^{(k)}$, restrict $\nu$ and the given action to this subgroup and use the elements $g_i^k$.
Since $\nu(g_i^k)=k\mathbf e_i\in\Z_{>0}\mathbf e_i$, condition $(\mathrm D)$ holds for these elements.
The identities for $\varphi_j$ and $\eta$ remain valid on this subgroup.
Thus \cref{thm:integer-to-equal-weight} applies to the restricted action.
\end{proof}

\subsection{Rational Designs on Spheres, Matrix Groups, and Grassmannians}\label{subsec:tq-applications}

We verify the density and integral conditions in \cref{prop:cxx-integer-weighted}, then construct the group actions required by \cref{cor:norm-one-equal-weight}.

\Needspace{5\baselineskip}

\subsubsection{Rationality, density, and group actions}\label{subsubsec:rational-homogeneous-input}

\begin{lemma}\label{lem:topological-inputs}
The following hold:
\begin{enumerate}
\item[\rm(i)] Let a compact group $\Gamma$ act continuously and transitively on a compact Hausdorff space $Z$, and let $\mu$ be a $\Gamma$-invariant probability measure.
Then every nonempty open subset of $Z$ has positive measure.
\item[\rm(ii)] Let $Z$ be a real-analytic manifold with finitely many connected components, and let $\varphi$ be a complex-valued real-analytic function on $Z$.
If $\varphi$ is not identically zero on any connected component, then $\{\mathbf{x} \in Z : \varphi(\mathbf{x}) \neq 0 \}$ is open and dense.
\end{enumerate}
\end{lemma}
\begin{proof}
For \textup{(i)}, suppose that a nonempty open subset of $Z$ has measure zero.
By transitivity its $\Gamma$-translates cover $Z$, and compactness gives a finite subcover.
Measure invariance would make every member of this subcover have measure zero, contradicting $\mu(Z)=1$; see also \cite[Proposition~2.60, p.~67]{Folland2016}.

For \textup{(ii)}, the set $\{\mathbf{x}\in Z:\varphi(\mathbf{x})\ne0\}$ is open by continuity.
On each connected component, the identity principle implies that $\lvert \varphi\rvert^2$ is not identically zero in any connected analytic coordinate chart.
In each such chart, the zero set of $\lvert \varphi\rvert^2$ has empty interior by \cite{Mityagin2020}.
Thus $\{\mathbf{x}\in Z:\varphi(\mathbf{x})\ne0\}$ is dense in $Z$.
\end{proof}

\begin{proposition}\label{prop:rational-invariant-average}
Let a compact matrix group $\Gamma$ act continuously and transitively on a compact set $Z$, and let $\mu$ be a $\Gamma$-invariant probability measure.
For $g \in \Gamma$ and a continuous real-valued function $\varphi$ on $Z$, define $(g\cdot \varphi)(\mathbf{x}):=\varphi(g^{-1}\cdot\mathbf{x})$.
Let $\Gamma_0\subseteq \Gamma$ be a dense subgroup, and let $Z_0 \subseteq Z$ be a $\Gamma_0$-invariant dense subset.
Let $\mathcal V_{\Q}$ be a finite-dimensional $\Q$-vector space of continuous real-valued functions on $Z$, containing the constant function $\mathbf1$ and invariant under $\Gamma_0$.
Assume that $\varphi(\mathbf{x}) \in \Q$ for every $\varphi \in \mathcal V_{\Q}$ and $\mathbf{x} \in Z_0$.
Then
\[
\int_Z\varphi\,\d\mu\in\Q \qquad (\varphi\in \mathcal V_{\Q}).
\]
\end{proposition}

\begin{proof}
Put $\mathcal V:=\operatorname{span}_{\R}\mathcal V_{\Q}$ and choose a $\Q$-basis $({\epsilon_{i}})_{i\in I}$ with finite $I\subseteq\Z_{>0}$, $1\in I$, and ${\epsilon_{1}}=\mathbf1$.
The evaluation vectors $\bigl({\epsilon_{i}}(\mathbf{x})\bigr)_{i\in I}$, $\mathbf{x}\in Z_0$, span $\Q^I$; otherwise a nonzero $\Q$-linear combination of the ${\epsilon_{i}}$ would vanish on the dense set $Z_0$, and hence on $Z$.
Choose $(\mathbf{x}^{(j)})_{j\in I}$ in $Z_0$ so that $\det\!\bigl({\epsilon_{i}}(\mathbf{x}^{(j)})\bigr)_{i,j\in I}\ne0$.
This rational matrix is also invertible over $\R$, so the ${\epsilon_{i}}$ form an $\R$-basis of $\mathcal V$.

Since $\Gamma_0$ is dense, continuity of the action and finite-dimensionality of $\mathcal V$ imply that $\mathcal V$ is $\Gamma$-invariant, and that every $\Gamma_0$-invariant real-linear functional on $\mathcal V$ is $\Gamma$-invariant.
Indeed, continuity of the action on the compact space $Z$ makes $g\mapsto g\cdot\varphi$ continuous in the uniform norm; the conclusions follow because $\mathcal V$ is closed in the uniform norm and every real-linear functional on $\mathcal V$ is continuous.
Transitivity makes the $\Gamma$-invariant functions in $\mathcal V$ constant.
Averaging any inner product on $\mathcal V$ over $\Gamma$ gives a $\Gamma$-invariant inner product, which identifies $\mathcal V$ with its real dual and identifies their $\Gamma$-invariant subspaces.
Thus the space of $\Gamma$-invariant real-linear functionals on $\mathcal V$ is one-dimensional, and the map
\[
\mathcal V\longrightarrow\R,\qquad \varphi\longmapsto\int_Z\varphi\,\d\mu
\]
is the unique $\Gamma_0$-invariant real-linear functional taking the value $1$ on $\mathbf1$.

Let ${\omega}:\mathcal V\to\R$ be a real-linear functional.
The conditions
\[
 {\omega}(g\cdot \varphi-\varphi)=0\quad(g\in \Gamma_0,\ \varphi\in \mathcal V_{\Q}),
 \qquad {\omega}(\mathbf1)=1
\]
form a system of linear equations with rational coefficients for $\bigl({\omega}({\epsilon_{i}})\bigr)_{i\in I}$.
The differences $g\cdot \varphi-\varphi$ span a finite-dimensional $\Q$-subspace of $\mathcal V_{\Q}$, so finitely many of these equations suffice.
By the uniqueness established above, this system has a unique real solution.
Gaussian elimination over $\Q$ therefore gives $\int_Z {\epsilon_{i}}\,\d\mu\in\Q$ for every $i\in I$, and hence $\int_Z \varphi\,\d\mu\in\Q$ for every $\varphi\in \mathcal V_{\Q}$.
\end{proof}

For a field $\mathbb K$ and an integer ${\mathsf{m}}\geq1$, let $\operatorname{Mat}_{\mathsf{m}}(\mathbb K)$ and $\mathrm{GL}_{\mathsf{m}}(\mathbb K)$ denote the ${\mathsf{m}}\times {\mathsf{m}}$ matrices and invertible matrices over $\mathbb K$.
We indicate matrix entries by $(\cdot)_{ij}$; the superscript $*$ denotes conjugate transpose.
The identity and zero matrices are
\[
\operatorname{Id}_{\mathsf{m}}:=\operatorname{diag}(1,\ldots,1),\qquad
0_{{\mathsf{m}}\times {\mathsf{m}}}:=(0)_{1\leq i,j\leq {\mathsf{m}}}.
\]
We also set $\operatorname{Id}_0=0_{0\times0}$ equal to the empty $0\times0$ matrix; these zero-dimensional blocks are omitted in direct sums and diagonal matrices.

For an integer ${\mathsf{m}} \geq 1$, define
\[
\begin{aligned}
\mathrm{O}({\mathsf{m}})&:=\{g\in \operatorname{Mat}_{\mathsf{m}}(\R):g^{\mathsf T}g=\operatorname{Id}_{\mathsf{m}}\},
& \mathrm{O}({\mathsf{m}};\Q)&:=\mathrm{O}({\mathsf{m}})\cap \operatorname{Mat}_{\mathsf{m}}(\Q),\\
\mathrm{SO}({\mathsf{m}})&:=\{g\in \mathrm{O}({\mathsf{m}}):\det g=1\},
& \mathrm{SO}({\mathsf{m}};\Q)&:=\mathrm{SO}({\mathsf{m}})\cap \operatorname{Mat}_{\mathsf{m}}(\Q),\\
\mathrm{U}({\mathsf{m}})&:=\{g\in \operatorname{Mat}_{\mathsf{m}}(\C):g^*g=\operatorname{Id}_{\mathsf{m}}\},
& \mathrm{U}({\mathsf{m}};\Q(\mathrm i))&:=\mathrm{U}({\mathsf{m}})\cap \operatorname{Mat}_{\mathsf{m}}(\Q(\mathrm i)).
\end{aligned}
\]
Identify $\operatorname{Mat}_{\mathsf{m}}(\C)$ with $\R^{2{\mathsf{m}}^2}$ by real and imaginary parts, so that $\operatorname{Mat}_{\mathsf{m}}(\Q(\mathrm i))$ is its rational point set.
Equip the three groups with normalized Haar measure.
For the remaining applications, let ${\mathsf{m}}\geq2$ and $1 \leq {\mathsf{r}} < {\mathsf{m}}$ be integers, and use the projection-matrix models
\[
\begin{aligned}
\operatorname{Gr}_{\R}({\mathsf{r}},{\mathsf{m}})
&:=\{M_2\in \operatorname{Mat}_{\mathsf{m}}(\R):M_2^{\mathsf T}=M_2,\ M_2^2=M_2,\ \operatorname{tr}M_2={\mathsf{r}}\},\\
\operatorname{Gr}_{\C}({\mathsf{r}},{\mathsf{m}})
&:=\{M_2\in \operatorname{Mat}_{\mathsf{m}}(\C):M_2^*=M_2,\ M_2^2=M_2,\ \operatorname{tr}M_2={\mathsf{r}}\}.
\end{aligned}
\]
Their rational point sets are
\[
\operatorname{Gr}_{\R}({\mathsf{r}},{\mathsf{m}})(\Q):=\operatorname{Gr}_{\R}({\mathsf{r}},{\mathsf{m}})\cap \operatorname{Mat}_{\mathsf{m}}(\Q),
\qquad \operatorname{Gr}_{\C}({\mathsf{r}},{\mathsf{m}})(\Q(\mathrm i)):=\operatorname{Gr}_{\C}({\mathsf{r}},{\mathsf{m}})\cap \operatorname{Mat}_{\mathsf{m}}(\Q(\mathrm i)).
\]
We equip them with the pushforwards of normalized Haar measure under
\[
g\longmapsto g\operatorname{diag}(\operatorname{Id}_{{\mathsf{r}}},0_{({\mathsf{m}}-{\mathsf{r}})\times({\mathsf{m}}-{\mathsf{r}})})g^{\mathsf T}\quad(g\in \mathrm{SO}({\mathsf{m}})),
\qquad g\longmapsto g\operatorname{diag}(\operatorname{Id}_{{\mathsf{r}}},0_{({\mathsf{m}}-{\mathsf{r}})\times({\mathsf{m}}-{\mathsf{r}})})g^*\quad(g\in \mathrm{U}({\mathsf{m}})),
\]
respectively; see \cite{BachocCoulangeonNebe2002,BachocBannaiCoulangeon2004,Roy2008}.

\begin{proposition}\label{prop:classical-homogeneous-rational-input}
Let $d,{\mathsf{m}}\geq2$ and $1\leq {\mathsf{r}}<{\mathsf{m}}$ be integers.
Equip the sphere $\mathbb S^{d-1}$ with normalized invariant measure and the matrix groups and Grassmannians with the measures specified above.
Then the corresponding rational point sets
\[
\mathbb S^{d-1}(\Q),\quad \mathrm{O}({\mathsf{m}};\Q),\quad \mathrm{SO}({\mathsf{m}};\Q),\quad \mathrm{U}({\mathsf{m}};\Q(\mathrm{i})),\quad \operatorname{Gr}_{\R}({\mathsf{r}},{\mathsf{m}})(\Q),\quad \operatorname{Gr}_{\C}({\mathsf{r}},{\mathsf{m}})(\Q(\mathrm{i}))
\]
are dense, every nonempty open set has positive measure, and the integral of every real polynomial with rational coefficients is rational.
\end{proposition}
\begin{proof}
Schmutz proves that $\mathrm{O}({\mathsf{m}};\Q)$ is dense in $\mathrm{O}({\mathsf{m}})$ \cite[Theorem~3.1]{Schmutz2008}.
Since $\mathrm{SO}({\mathsf{m}})$ is an open component of $\mathrm{O}({\mathsf{m}})$, this implies that $\mathrm{SO}({\mathsf{m}};\Q)=\mathrm{O}({\mathsf{m}};\Q)\cap \mathrm{SO}({\mathsf{m}})$ is dense in $\mathrm{SO}({\mathsf{m}})$.
\Needspace{5\baselineskip}
For the unitary group, the Cayley map on the skew-Hermitian matrices in $\operatorname{Mat}_{\mathsf{m}}(\C)$,
\[
M_3\longmapsto(\operatorname{Id}_{\mathsf{m}}-M_3)(\operatorname{Id}_{\mathsf{m}}+M_3)^{-1}
\]
is a homeomorphism onto $\{g\in \mathrm{U}({\mathsf{m}}):\det(\operatorname{Id}_{\mathsf{m}}+g)\ne0\}$, with inverse
\[
g\longmapsto(\operatorname{Id}_{\mathsf{m}}-g)(\operatorname{Id}_{\mathsf{m}}+g)^{-1}.
\]
The formulas are rational in the real and imaginary parts of the entries; $\operatorname{Id}_{\mathsf{m}}+M_3$ is invertible because the eigenvalues of a skew-Hermitian matrix are purely imaginary.
This image is dense by \cref{lem:topological-inputs}(ii), because $\mathrm{U}({\mathsf{m}})$ is connected and $g\mapsto\det(\operatorname{Id}_{\mathsf{m}}+g)$ is nonzero at $\operatorname{Id}_{\mathsf{m}}$.
Skew-Hermitian matrices with entries in $\Q(\mathrm i)$ are dense in the skew-Hermitian matrices, and the Cayley map sends them into $\mathrm{U}({\mathsf{m}};\Q(\mathrm i))$; hence $\mathrm{U}({\mathsf{m}};\Q(\mathrm i))$ is dense in $\mathrm{U}({\mathsf{m}})$.

Every unit vector extends to an oriented orthonormal basis, so the orbit map $g\mapsto g(1,0,\ldots,0)^{\mathsf T}$ from $\mathrm{SO}(d)$ onto $\mathbb S^{d-1}$ is surjective.
Likewise, the spectral theorem gives the surjective orbit maps
\[
g\mapsto g\operatorname{diag}(\operatorname{Id}_{{\mathsf{r}}},0_{({\mathsf{m}}-{\mathsf{r}})\times({\mathsf{m}}-{\mathsf{r}})})g^{\mathsf T}\quad(g\in\mathrm{SO}({\mathsf{m}})),
\qquad g\mapsto g\operatorname{diag}(\operatorname{Id}_{{\mathsf{r}}},0_{({\mathsf{m}}-{\mathsf{r}})\times({\mathsf{m}}-{\mathsf{r}})})g^*\quad(g\in\mathrm{U}({\mathsf{m}})),
\]
onto the real and complex Grassmannians, respectively.
In the real case, an orthogonal diagonalizer can be made to have determinant $1$ by changing the sign of one column, without changing the projection.
The images of the dense rational subgroups under these continuous maps are rational points and are dense in the corresponding spaces.

Use left multiplication for the matrix groups, the usual linear action on the sphere, and conjugation for the Grassmannians.
Apply \cref{prop:rational-invariant-average} with $(\Gamma,\Gamma_0,Z_0)$ equal to
\[
\begin{gathered}
(\mathrm{SO}(d),\mathrm{SO}(d;\Q),\mathbb S^{d-1}(\Q)),\\
(\mathrm{O}({\mathsf{m}}),\mathrm{O}({\mathsf{m}};\Q),\mathrm{O}({\mathsf{m}};\Q)),\quad (\mathrm{SO}({\mathsf{m}}),\mathrm{SO}({\mathsf{m}};\Q),\mathrm{SO}({\mathsf{m}};\Q)),\\
(\mathrm{U}({\mathsf{m}}),\mathrm{U}({\mathsf{m}};\Q(\mathrm i)),\mathrm{U}({\mathsf{m}};\Q(\mathrm i))),\\
(\mathrm{SO}({\mathsf{m}}),\mathrm{SO}({\mathsf{m}};\Q),\operatorname{Gr}_{\R}({\mathsf{r}},{\mathsf{m}})(\Q)),\quad (\mathrm{U}({\mathsf{m}}),\mathrm{U}({\mathsf{m}};\Q(\mathrm i)),\operatorname{Gr}_{\C}({\mathsf{r}},{\mathsf{m}})(\Q(\mathrm i)))
\end{gathered}
\]
for the sphere, the two orthogonal groups, the unitary group, and the two Grassmannians, respectively.
For each integer $t\geq0$, let $\mathcal V_\Q$ consist of the restrictions of polynomials of degree at most $t$ with rational coefficients in the ambient real coordinates.
The indicated rational subgroups act by linear maps with rational coefficients in these coordinates, so $\mathcal V_\Q$ is invariant under each subgroup and its functions take rational values on the corresponding rational point set.
This proves that $\int_Z \varphi\,\d\mu\in\Q$ for every restriction $\varphi$ of a polynomial with rational coefficients, and \cref{lem:topological-inputs}(i) gives positive measure to every nonempty open set.
\end{proof}

\begin{proposition}\label{prop:tq-action-data}
Let $d,{\mathsf{m}}\geq2$ and $1\leq {\mathsf{r}}<{\mathsf{m}}$ be integers.
The group $\TQ$ is $T_{\Q(\mathrm i)}$ in the notation of \cref{lem:norm-one-condition-D}.
Each of the five actions below, on the corresponding space $Z$, preserves its rational points and invariant probability measure and is real-linear in the ambient coordinates, hence preserves polynomial degree.
Moreover, for every integer $t \geq 0$, the space $\mathcal P_t(Z;\R)$ is generated by the real and imaginary parts of finitely many complex-valued functions $\varphi$ satisfying
\[
\varphi({v}\cdot\mathbf{x})={v}^k\varphi(\mathbf{x}) \qquad ({v} \in \TQ,\ \mathbf{x} \in Z)
\]
for some $k \in \mathbb{Z}$, with the constant function corresponding to $k = 0$.
\begin{enumerate}
\item For the sphere, take $\mathbf{x}=(x_1,\ldots,x_d)\in\mathbb S^{d-1}$ and define
\[
{v}\cdot(x_1,x_2,x_3,\ldots,x_d) := (\Re({v}(x_1+\mathrm i x_2)),\Im({v}(x_1+\mathrm i x_2)),x_3,\ldots,x_d).
\]
For nonnegative integers $k_1,\ldots,k_d$ with $\sum_{j=1}^d k_j\leq t$, we have
\[
({v}(x_1+\mathrm i x_2))^{k_1}\overline{{v}(x_1+\mathrm i x_2)}^{\,k_2}\prod_{j=3}^d x_j^{k_j} ={v}^{k_1-k_2}(x_1+\mathrm i x_2)^{k_1}\overline{(x_1+\mathrm i x_2)}^{\,k_2}\prod_{j=3}^d x_j^{k_j}.
\]
\item For $Z=\mathrm{O}({\mathsf{m}})$ or $\mathrm{SO}({\mathsf{m}})$, put
\[
\operatorname{rot}({v}):=\begin{pmatrix}\Re {v}&-\Im {v}\\ \Im {v}&\Re {v}\end{pmatrix}\in \mathrm{SO}(2;\Q)\qquad({v}\in\TQ).
\]
Define
\[
{v}\cdot g':=(\operatorname{rot}({v})\oplus \operatorname{Id}_{{\mathsf{m}}-2})g'\qquad(g'\in Z).
\]
Then
\[
({v}\cdot g')_{1j}+\mathrm i({v}\cdot g')_{2j}
={v}\bigl((g')_{1j}+\mathrm i(g')_{2j}\bigr)
\qquad(1\leq j\leq{\mathsf{m}}),
\]
and the entries in rows $3,\ldots,{\mathsf{m}}$ are invariant.
\item For the unitary group, define ${v}\cdot g' := \operatorname{diag}({v},1,\ldots,1)g'$ for $g'\in \mathrm{U}({\mathsf{m}})$.
The entries in the first row are multiplied by ${v}$, their complex conjugates are multiplied by ${v}^{-1}$, and the entries in the other rows are invariant.
\item For the real Grassmannian, define
\[
{v}\cdot M_2 := (\operatorname{rot}({v})\oplus \operatorname{Id}_{{\mathsf{m}}-2})M_2(\operatorname{rot}({v}) \oplus \operatorname{Id}_{{\mathsf{m}}-2})^{\mathsf T}\qquad(M_2\in\operatorname{Gr}_{\R}({\mathsf{r}},{\mathsf{m}})).
\]
Then
\[
({v}\cdot M_2)_{11}-({v}\cdot M_2)_{22}+2\mathrm i({v}\cdot M_2)_{12}
={v}^2((M_2)_{11}-(M_2)_{22}+2\mathrm i(M_2)_{12})
\]
and
\[
({v}\cdot M_2)_{1j}+\mathrm i({v}\cdot M_2)_{2j}
={v}((M_2)_{1j}+\mathrm i(M_2)_{2j})
\qquad(3\leq j\leq {\mathsf{m}}),
\]
while the entries $(M_2)_{ij}$ with $i,j\geq3$ and the sum $(M_2)_{11}+(M_2)_{22}$ are invariant.
\item For the complex Grassmannian, define
\[
{v}\cdot M_2 := \operatorname{diag}({v},1,\ldots,1)M_2\operatorname{diag}(\overline {v},1,\ldots,1)\qquad(M_2\in\operatorname{Gr}_{\C}({\mathsf{r}},{\mathsf{m}})).
\]
Then the entry $(M_2)_{1j}$ is multiplied by ${v}$, and $(M_2)_{j1}$ is multiplied by ${v}^{-1}$ for $2\leq j\leq {\mathsf{m}}$, and the other entries are invariant.
\end{enumerate}
\end{proposition}

\begin{proof}
The acting matrices have entries in \(\Q\) or \(\Q(\mathrm{i})\), and each action is the restriction of the corresponding orthogonal or unitary action.
Thus rational points and the invariant measures are preserved.
Direct matrix multiplication verifies the displayed formulas.

For the real Grassmannian, the invariant sum $(M_2)_{11}+(M_2)_{22}$ together with the real and imaginary parts of $(M_2)_{11}-(M_2)_{22}+2\mathrm i(M_2)_{12}$ and $(M_2)_{1j}+\mathrm i(M_2)_{2j}$, $3\leq j\leq {\mathsf{m}}$, recovers all entries involving the first two rows or columns.
Thus, in every case, each real coordinate function is a real linear combination of real and imaginary parts of the functions listed above, each of which transforms by an integer power of ${v}$.
Expanding in these functions and their conjugates expresses each real coordinate monomial of degree at most $t$ as a complex linear combination of products of at most $t$ linear functions of this form.
Since $\overline {v}={v}^{-1}$ for ${v}\in\TQ$, these products transform by integer powers of ${v}$ and have degree at most $t$.
Taking real parts of the expansion puts each real coordinate monomial in the real linear span of the real and imaginary parts of these products; conversely, all these parts belong to $\mathcal P_t(Z;\R)$.
\end{proof}

\begin{proposition}\label{prop:tq-nonzero-loci}
For each space $Z$ and action in \cref{prop:tq-action-data},
take $\eta:Z\to\C$ and $k_0\in\Z_{>0}$ as follows:
\[
\begin{array}{c|c|c}
Z & \eta & k_0\\ \hline
\mathbb S^{d-1} & \mathbf{x}\mapsto x_1+\mathrm i x_2 & 1\\
\mathrm{O}({\mathsf{m}}),\,\mathrm{SO}({\mathsf{m}}) & g'\mapsto((g')_{11}+\mathrm i(g')_{21}) & 1\\
\mathrm{U}({\mathsf{m}}) & g'\mapsto(g')_{11} & 1\\
\operatorname{Gr}_{\R}({\mathsf{r}},{\mathsf{m}}) & M_2\mapsto (M_2)_{11}-(M_2)_{22}+2\mathrm i (M_2)_{12} & 2\\
\operatorname{Gr}_{\C}({\mathsf{r}},{\mathsf{m}}) & M_2\mapsto (M_2)_{12} & 1
\end{array}
\]
Then $\eta({v}\cdot\mathbf{x})={v}^{k_0}\eta(\mathbf{x})$ for all ${v}\in\TQ$ and $\mathbf{x}\in Z$.
Moreover, the set
\[
W:=\{\mathbf{x}\in Z:\eta(\mathbf{x})\ne0\}
\]
is $\TQ$-invariant, open, and dense.
Its rational point set $W(\Q)=W\cap Z(\Q)$ is dense in $Z$, and $\eta(W(\Q))\subseteq\Q(\mathrm i)^\times$.
\end{proposition}
\begin{proof}
By \cref{prop:tq-action-data}, $\eta({v}\cdot\mathbf{x})={v}^{k_0}\eta(\mathbf{x})$; since ${v}\ne0$, this makes $W$ invariant.
The definitions of $\eta$ give $\eta(W(\Q))\subseteq\Q(\mathrm i)^\times$, and continuity of $\eta$ makes $W$ open.

It remains to prove density of $W$.
Each $\eta$ is real analytic.
The groups $\mathrm{SO}({\mathsf{m}})$ and $\mathrm{U}({\mathsf{m}})$ are connected, whereas $\mathrm{O}({\mathsf{m}})$ has two connected components distinguished by the determinant \cite[Chapter~1, pp.~3--30]{Hall2015}.
The sphere and the real and complex Grassmannians are connected because the connected groups $\mathrm{SO}(d)$, $\mathrm{SO}({\mathsf{m}})$, and $\mathrm{U}({\mathsf{m}})$ act transitively on them, respectively.

In the sphere case, $\eta((1,0,\ldots,0))=1$.
For $\mathrm{SO}({\mathsf{m}})$ and $\mathrm{U}({\mathsf{m}})$, $\eta(\operatorname{Id}_{\mathsf{m}})=1$.
In the orthogonal case, the two connected components of $\mathrm{O}({\mathsf{m}})$ contain, respectively,
\[
\operatorname{Id}_{\mathsf{m}} \qquad \text{and} \qquad \operatorname{diag}(-1,1,\ldots,1),
\]
and $\eta$ takes the values $1$ and $-1$ at these points, respectively.

For the real Grassmannian, $\eta$ takes the value $1$ at the rank-${\mathsf{r}}$ orthogonal projection
\[
\operatorname{diag}(1,0,\operatorname{Id}_{{\mathsf{r}}-1},0_{({\mathsf{m}}-{\mathsf{r}}-1)\times({\mathsf{m}}-{\mathsf{r}}-1)}).
\]

For the complex Grassmannian, the matrix
\[
\frac12 \begin{pmatrix} 1&1\\ 1&1 \end{pmatrix} \oplus \operatorname{Id}_{{\mathsf{r}}-1}\oplus 0_{({\mathsf{m}}-{\mathsf{r}}-1)\times({\mathsf{m}}-{\mathsf{r}}-1)}
\]
belongs to $\operatorname{Gr}_{\C}({\mathsf{r}},{\mathsf{m}})$, since its $2\times2$ block is a rank-one orthogonal projection.
Its $(1,2)$-entry is $\frac12$, so $\eta$ is nonzero at this point.

Hence, in every case, $\eta$ is not identically zero on any connected component of $Z$.
By \cref{lem:topological-inputs}(ii), it follows that $W$ is dense in $Z$.
Since $W$ is open and $Z(\Q)$ is dense by \cref{prop:classical-homogeneous-rational-input}, $W(\Q)$ is dense in $Z$.
\end{proof}

\subsubsection{Applications to classical homogeneous spaces}

In each application, the restriction to $W(\Q)$ of the map $\eta$ in \cref{prop:tq-nonzero-loci} is the map required in condition~\textup{(c)} of \cref{thm:integer-to-equal-weight} and used in \cref{cor:norm-one-equal-weight}.
Proposition~\ref{prop:cxx-integer-weighted} supplies an integer-weighted rational design supported on $W(\Q)$, and \cref{cor:norm-one-equal-weight} transforms it into
an equi-weighted rational design
supported there and disjoint from the prescribed finite set.

The following theorem gives rational spherical designs in every dimension $d\geq2$ and every strength, resolving Xiang's conjecture~\cite[Conjecture~1.4]{Xiang2022}.

\begin{theorem}\label{thm:rational-spherical-existence}
For every pair of integers $d\geq2$ and $t\geq1$ and every finite set $E\subseteq \mathbb{S}^{d-1}(\Q)$, there exists a rational spherical $t$-design $\mathfrak{X}$ such that
\[
\mathfrak{X} \subseteq \{\mathbf{x} \in \mathbb{S}^{d-1}(\Q) : x_1 + \mathrm{i}x_2 \neq 0\} \setminus E.
\]
\end{theorem}

\begin{proof}
By \cref{prop:classical-homogeneous-rational-input}, the integral of every polynomial with rational coefficients on the sphere is rational, and every nonempty open set has positive measure.
The rational points in $W=\{\mathbf{x}\in \mathbb S^{d-1}:x_1+\mathrm i x_2\ne0\}$ are dense by \cref{prop:tq-nonzero-loci}, so \cref{prop:cxx-integer-weighted} supplies an integer-weighted rational spherical $t$-design supported there.
Use the sphere action in \cref{prop:tq-action-data} and apply \cref{cor:norm-one-equal-weight} with $\mathcal V=\mathcal P_t(\mathbb S^{d-1};\R)$, $\mathbb K=\Q(\mathrm i)$, $\eta(\mathbf{x})=x_1+\mathrm i x_2$, $k_0=1$, and the prescribed set $E$.
\end{proof}

\begin{remark}\label{rem:cxx-alternative-proof}
\begin{enumerate}
\item[\rm(i)]
The existence of a rational $t$-design on a space $Z$ satisfying the hypotheses of \cite[Theorem~7.1]{CuiXiaXiang2019} can also be proved using \cref{thm:rational-spherical-existence}.
The rational interval designs in \textup{(ii)} can be incorporated into the construction in \cite[Section~7.1]{CuiXiaXiang2019}, replacing the point-separation argument used there.

\item[\rm(ii)]
For integers $d\geq2$ and $r\geq0$, a nonempty finite set $\mathfrak{X} \subseteq (-1,1) \cap \mathbb{Q}$ is called a rational $r$-design for
the Gegenbauer measure
\[
\frac{(1-u^2)^{(d-3)/2}\,\mathrm du}
{\int_{-1}^{1}(1-u^2)^{(d-3)/2}\,\mathrm du}
\]
on $(-1,1)$
if
\[
\frac{1}{|\mathfrak{X}|}\sum_{u \in \mathfrak{X}}\varphi(u)
=
\frac{1}{\int_{-1}^{1}(1-u^2)^{(d-3)/2}\,\mathrm du} \int_{-1}^1 \varphi(u)\, (1-u^2)^{(d-3)/2}\,\mathrm du
\quad \bigl(\varphi\in\mathbb R[X],\ \deg\varphi\leq r\bigr);
\]
see for example Xiang~\cite{Xiang2022} and Mishima et al.~\cite{MishimaLuSawaUchida2026}.
For $d=2$, this measure is the Chebyshev measure.
Taking inner products of the points of a rational spherical $r$-design on $\mathbb S^{d-1}$, $d\geq2$, with a rational unit vector yields a rational $r$-design for the Gegenbauer measure, provided that the projected values are distinct and lie in $(-1,1)$.
The underlying method, which also applies in the real setting, is described in \cite[Section~2.1, Eq.~(2.1)]{BilykGlazyrinMatzkeParkVlasiuk2019}.
\end{enumerate}
\end{remark}

For designs on orthogonal and special orthogonal groups, see Hashagen--Flammia--Gross--Wallman \cite{HashagenFlammiaGrossWallman2018} and Okuda \cite{Okuda2014}, respectively.
We work with the real vector spaces
\[
\mathcal P_t(\mathrm O(\mathsf m);\mathbb R), \qquad \mathcal P_t(\mathrm{SO}(\mathsf m);\mathbb R),
\]
consisting of restrictions of
real polynomials of total degree at most \(t\) in the matrix entries.
We consider designs that reproduce integrals of these functions with respect to normalized Haar measure by equal-weight averages over finite sets.
We use the action rotating the first two rows
specified in \cref{prop:tq-action-data}.

\begin{proposition}\label{prop:rational-orthogonal-designs}
For every pair of integers \({\mathsf{m}}\geq2\) and \(t\geq1\) and every finite set \(E\subseteq \mathrm{O}({\mathsf{m}};\Q)\), there exists a rational \(\mathcal P_t(\mathrm{O}({\mathsf{m}});\R)\)-design disjoint from \(E\).
Likewise, for every finite set \(E\subseteq \mathrm{SO}({\mathsf{m}};\Q)\), there exists a rational \(\mathcal P_t(\mathrm{SO}({\mathsf{m}});\R)\)-design disjoint from \(E\).
\end{proposition}

\begin{proof}
For $Z=\mathrm{O}({\mathsf{m}})$ and $Z=\mathrm{SO}({\mathsf{m}})$, \cref{prop:classical-homogeneous-rational-input} gives rational integrals for polynomials with rational coefficients and positive measure on every nonempty open set, and \cref{prop:tq-nonzero-loci} gives density of $\{g\in Z(\Q):(g_{11}+\mathrm i g_{21})\ne0\}$.
Thus \cref{prop:cxx-integer-weighted} supplies integer-weighted rational $\mathcal P_t(Z;\R)$-designs on these subsets.
Use the respective actions in \cref{prop:tq-action-data} and apply \cref{cor:norm-one-equal-weight} with $\mathcal V=\mathcal P_t(Z;\R)$, $\mathbb K=\Q(\mathrm i)$, $\eta(g)=g_{11}+\mathrm i g_{21}$, $k_0=1$, and the prescribed finite set $E$ in each group.
\end{proof}

For designs on unitary groups, see Roy--Scott \cite{RoyScott2009}.
For integers $\mathsf{m} \ge 1$ and $t\ge0$, let $\mathcal{P}_{t,t}(\mathrm U(\mathsf m);\mathbb{C})$ denote the complex vector space
consisting of restrictions of polynomials homogeneous of degree $t$ in the matrix entries and of degree $t$ in their complex conjugates.
Unitary $t$-designs are designs for the real vector space
\[
\operatorname{span}_{\mathbb R} \left\{ \Re\varphi,\Im\varphi: \varphi\in\mathcal P_{t,t}(\mathrm U(\mathsf m);\mathbb C) \right\}.
\]
This real vector space contains the constant functions and is contained in $\mathcal{P}_{2t}(\mathrm U(\mathsf m);\mathbb{R})$, so it suffices to construct a $\mathcal{P}_{2t}$-design.
We use the action multiplying the first row by a unit complex number specified in \cref{prop:tq-action-data}.

\begin{proposition}\label{prop:rational-unitary-designs}
For every pair of integers ${\mathsf{m}},t\geq1$ and every finite set $E\subseteq \mathrm{U}({\mathsf{m}};\Q(\mathrm{i}))$, there exists a unitary $t$-design $\mathfrak{X} \subseteq \mathrm{U}({\mathsf{m}};\Q(\mathrm{i}))$ disjoint from $E$.
\end{proposition}
\begin{proof}
For ${\mathsf{m}}=1$, the group $\mathrm{U}(1;\Q(\mathrm i))=\TQ$ is infinite by \cref{lem:norm-one-condition-D}, and $\mathcal P_{t,t}(\mathrm{U}(1);\C)=\C\mathbf1$ because $|(g)_{11}|^{2t}=1$ for every $g\in \mathrm{U}(1)$.
Thus a singleton outside $E$ is a unitary $t$-design.

Let ${\mathsf{m}} \geq2$.
By \cref{prop:classical-homogeneous-rational-input,prop:tq-nonzero-loci,prop:cxx-integer-weighted}, there exists an integer-weighted rational $\mathcal P_{2t}(\mathrm{U}({\mathsf{m}});\R)$-design supported on
\[
\{g\in \mathrm{U}({\mathsf{m}};\Q(\mathrm i)):(g)_{11}\ne0\}.
\]
Use the unitary action in \cref{prop:tq-action-data} and apply \cref{cor:norm-one-equal-weight} to $\mathcal P_{2t}(\mathrm{U}({\mathsf{m}});\R)$ with $\mathbb K=\Q(\mathrm i)$, $\eta(g)=(g)_{11}$, $k_0=1$, and the prescribed finite set $E$.
Thus the resulting set is a unitary $t$-design.
\end{proof}

\Needspace{5\baselineskip}
For designs on Grassmannians, see Bachoc--Coulangeon--Nebe~\cite{BachocCoulangeonNebe2002} for the real case and Roy~\cite{Roy2008} for the complex case.
We work with the real vector spaces
\[
\mathcal{P}_t\bigl(\operatorname{Gr}_{\mathbb R}(\mathsf r,\mathsf m);\mathbb R\bigr),
\qquad \mathcal P_t\bigl(\operatorname{Gr}_{\mathbb C}(\mathsf r,\mathsf m);\mathbb R\bigr),
\]
consisting of restrictions of real polynomials of total degree at most $t$ in the entries of the projection matrices in the real case, and in their real and imaginary parts in the complex case.
We use the conjugation actions specified in \cref{prop:tq-action-data}.

\begin{proposition}\label{prop:rational-grassmann-designs}
For all integers ${\mathsf{r}},{\mathsf{m}},t$ with $1\leq {\mathsf{r}}<{\mathsf{m}}$ and $t\geq1$ and every finite set $E\subseteq\operatorname{Gr}_{\R}({\mathsf{r}},{\mathsf{m}})(\Q)$, there exists a $\mathcal P_t(\operatorname{Gr}_{\R}({\mathsf{r}},{\mathsf{m}});\R)$-design disjoint from $E$ whose projection matrices have rational entries.
Likewise, for every finite set $E\subseteq\operatorname{Gr}_{\C}({\mathsf{r}},{\mathsf{m}})(\Q(\mathrm{i}))$, there exists a $\mathcal P_t(\operatorname{Gr}_{\C}({\mathsf{r}},{\mathsf{m}});\R)$-design disjoint from $E$ whose projection matrices have entries in $\Q(\mathrm{i})$.
\end{proposition}
\begin{proof}
By \cref{prop:classical-homogeneous-rational-input}, both Grassmannians have rational integrals for polynomials with rational coefficients and positive measure on every nonempty open set.
By \cref{prop:tq-nonzero-loci}, the rational points with $(M_2)_{11}-(M_2)_{22}+2\mathrm i(M_2)_{12}\ne0$ in the real case and $(M_2)_{12}\ne0$ in the complex case are dense.
Apply \cref{prop:cxx-integer-weighted} to these subsets.
Then apply \cref{cor:norm-one-equal-weight} to these integer-weighted designs with $\mathbb K=\Q(\mathrm i)$, the actions in \cref{prop:tq-action-data}, and the prescribed finite set $E$ in each space.
In the real case, use $\eta(M_2)=(M_2)_{11}-(M_2)_{22}+2\mathrm i(M_2)_{12}$ and $k_0=2$; in the complex case, use $\eta(M_2)=(M_2)_{12}$ and $k_0=1$.
\end{proof}

For designs on real and complex projective spaces, see Hoggar \cite{Hoggar1982}.
We identify a real or complex line $[\mathbf x]$ with its orthogonal projection
\[
[\mathbf x]\longmapsto \frac{\mathbf x\mathbf x^*}{\mathbf x^*\mathbf x} \qquad(\mathbf x\ne\mathbf0).
\]
Real and complex projective \(t\)-designs are designs for the respective real vector spaces
\[
\mathcal P_t\bigl(\operatorname{Gr}_{\mathbb R}(1,\mathsf m);\mathbb R\bigr),
\qquad \mathcal P_t\bigl(\operatorname{Gr}_{\mathbb C}(1,\mathsf m);\mathbb R\bigr),
\]
consisting of restrictions of real polynomials of total degree at most $t$ in the entries of the projection matrices, or in their real and imaginary parts in the complex case.
Taking $\mathsf r=1$ in \cref{prop:rational-grassmann-designs} yields the following corollary.

\begin{corollary}\label{cor:rational-projective-designs}
For every pair of integers ${\mathsf{m}}\geq2$ and $t\geq1$, $\mathbb RP^{{\mathsf{m}}-1}\cong\operatorname{Gr}_{\R}(1,{\mathsf{m}})$ has a design for $\mathcal P_t(\operatorname{Gr}_{\R}(1,{\mathsf{m}});\R)$ consisting of rational projection matrices.
Likewise, $\mathbb CP^{{\mathsf{m}}-1}\cong\operatorname{Gr}_{\C}(1,{\mathsf{m}})$ has a design for $\mathcal{P}_t(\operatorname{Gr}_{\C}(1,{\mathsf{m}});\R)$ consisting of projection matrices over $\Q(\mathrm{i})$.
\end{corollary}

\begin{remark}\label{rem:rational-designs-products}
\Cref{thm:integer-to-equal-weight} applies to a wide range of spaces, of which the spheres, matrix groups, and Grassmannians considered in this section are representative examples.
The theorem requires no assumptions of connectedness, smoothness, or a global product decomposition.
\end{remark}

\begin{remark}\label{rem:circle-arithmetic}
The case $d=2$ of \cref{thm:rational-spherical-existence} gives rational designs of every strength on $\mathbb{S}^1$.
In \cref{sec:local-arithmetic-minimum}, we study divisibility conditions on their sizes and determine all possible sizes for strengths $4$ through $7$.
\end{remark}

\endgroup

\section{Rational Designs on the Unit Circle}\label{sec:local-arithmetic-minimum}

We construct PTE set solutions on the rational unit circle using the sunflower lemma, then derive divisibility conditions for rational designs on positive-definite rational conics.
On the unit circle, these conditions and the constructions in \cref{thm:minimum-circle-constructions} determine all possible sizes for strengths $4$ through $7$, as stated in \cref{thm:circle-cardinality-spectra}.
We also give an explicit upper bound for the size of a rational spherical $11$-design.
We use $\psi_k$ and $R_4$ as defined in \eqref{eq:Pk} and \eqref{eq:roots-of-unity}, respectively.

\subsection{Rational Designs and PTE Solutions on the Unit Circle}\label{subsec:circle-designs-pte}\label{subsubsec:circle-pte-constructions}

We apply \cref{thm:integer-to-equal-weight} to obtain $R_4$-invariant rational spherical designs, then construct PTE set solutions by lifting from the fourth-power subgroup.

Identify $\mathbf x=(x_1,x_2)\in\mathbb S^1$ with $z=x_1+\mathrm i x_2$.
\begin{lemma}\label{lem:circle-moment-criterion}
Let $t\geq1$ be an integer.
A nonempty finite set \(\mathfrak X\subseteq \mathbb S^1\) is a spherical \(t\)-design if and only if
\[
\psi_k(\mathfrak X)=0\qquad(1\leq k\leq t).
\]
\end{lemma}
\begin{proof}
This is the circle case of the harmonic characterization of spherical designs; see, for example, \cite[Cor.~3.8 and Thm.~5.5]{DelsarteGoethalsSeidel1977}.
\end{proof}

The fourth-power subgroup of $\TQ$ is
\begin{equation} \label{eq:TQ4}
\TQ^{\,4} := \{z^4:z\in\TQ\}.
\end{equation}
The map
\begin{equation}\label{eq:circle-fourth-power-map}
\rho:\TQ\longrightarrow \TQ^{\,4},\qquad \rho(z):=z^4
\end{equation}
is a surjective group homomorphism with kernel $R_4$.

For ${v}\in\TQ$ and $\mathfrak X\subseteq\TQ$, write
\[
{v}\mathfrak X:=\{{v}z:z\in\mathfrak X\}.
\]

\begin{lemma} \label{lem:circle-coordinate-separating-rotation}
For every finite set $\mathfrak X \subseteq \TQ$, there exists ${v} \in \TQ$ such that the first coordinate $z \mapsto \Re z$ is injective on ${v}\mathfrak X$.
For every integer $t\geq1$, if $\mathfrak X$ is a spherical $t$-design, then so is ${v}\mathfrak X$.
If $\mathrm i\mathfrak X=\mathfrak X$, then $\mathrm i({v}\mathfrak X)={v}\mathfrak X$.
\end{lemma}

\begin{proof}
For distinct $z,z'\in\mathfrak X$, the identity
\[
2\bigl(\Re({v} z)-\Re({v} z')\bigr)
=(z-z')\bigl({v}-{v}^{-1}(zz')^{-1}\bigr)
\]
shows that $\Re({v} z)=\Re({v} z')$ if and only if ${v}^2=(zz')^{-1}$.
Each unordered pair in $\mathfrak X$ forbids at most two values of ${v}$.
The group $\TQ$ is infinite by \cref{lem:norm-one-condition-D}, so a value outside this finite forbidden set exists.
Rotation multiplies the $k$-th complex moment by ${v}^k$, so \cref{lem:circle-moment-criterion} shows that it preserves the strength of a spherical design.
The identity $\mathrm i({v}\mathfrak X)={v}(\mathrm i\mathfrak X)$ proves the last assertion.
\end{proof}

\begin{proposition}\label{prop:c4-invariant-rational-circle-design}
For every positive integer $t$, there exists a spherical $t$-design $\mathfrak X\subseteq\TQ$ satisfying
\[
\mathrm i\mathfrak X=\mathfrak X.
\]
Moreover, the map $z\mapsto\Re z$ is injective on $\mathfrak X$.
\end{proposition}

\begin{proof}
If $t < 4$, apply \cref{lem:circle-coordinate-separating-rotation} to $R_4$ and take $\mathfrak X:={v} R_4$ for the resulting ${v}$.
By \cref{lem:circle-moment-criterion}, $\mathfrak X$ is a spherical $t$-design.

Suppose $t\geq4$, and put $t':=\lfloor t/4\rfloor$.
By \cref{prop:classical-homogeneous-rational-input}, every nonempty open subset of $\mathbb S^1$ has positive measure and every real polynomial with rational coefficients has rational integral.
The subgroup $\TQ^{\,4}$ is dense in $\mathbb S^1$, since the fourth-power map on $\mathbb S^1$ is continuous and surjective and $\TQ$ is dense by \cref{prop:classical-homogeneous-rational-input}.
Thus \cref{prop:cxx-integer-weighted}, with $Z=\mathbb S^1$ and $Z_0=\TQ^{\,4}$, supplies an integer-weighted rational $\mathcal P_{t'}(\mathbb S^1;\R)$-design supported on $\TQ^{\,4}$.
Apply \cref{thm:integer-to-equal-weight} to this action of $\TQ^{\,4}$ with $\mathbb K=\Q(\mathrm i)$, $W=\TQ^{\,4}$, $\mathcal V=\mathcal P_{t'}(\mathbb S^1;\R)$, and $\eta(z)=z$.
Conditions \textup{(a)} and \textup{(b)} of that theorem follow by restricting the action in \cref{prop:tq-action-data} to $\TQ^{\,4}$.
Take $\mathfrak p_i,g_i$ as in \cref{lem:norm-one-condition-D} for $\mathbb K=\Q(\mathrm i)$, and restrict its valuation homomorphism $\nu$ to $\TQ^{\,4}$. Then the identity $\nu(g_i^4)=4\mathbf e_i$ verifies condition $(\mathrm D)$ on $\TQ^{\,4}$.
Condition \textup{(c)} holds with $\theta_i=\operatorname{ord}_{\mathfrak p_i}$ and $\chi_0({v})={v}$: one has $\theta_i\circ\chi_0=\nu_i$ on $\TQ^{\,4}$, and every $z\in \TQ^{\,4}$ has nonzero Gaussian-prime valuation at only finitely many primes.
Taking $E=\varnothing$ gives a spherical $t'$-design $\mathfrak X'\subseteq \TQ^{\,4}$.
The finite $R_4$-invariant set $\mathfrak X'' := \rho^{-1}(\mathfrak X')$ satisfies
\[
\psi_k(\mathfrak X'')=
\begin{cases}
0,&4\nmid k,\\
4\psi_{k/4}(\mathfrak X'),&4\mid k
\end{cases}
\qquad(1\leq k\leq t).
\]
Since $k/4\leq t'$ when $4\mid k$, these moments vanish; hence $\mathfrak X''$ is a spherical $t$-design by \cref{lem:circle-moment-criterion}.
Apply \cref{lem:circle-coordinate-separating-rotation} to $\mathfrak X''$ and put $\mathfrak X := {v} \mathfrak X''$ for the resulting ${v}$.
\end{proof}

For the following PTE constructions, fix an integer $t\geq4$ and put $t':=\lfloor t/4\rfloor$.
Recall the consecutive-power character families
\begin{equation}\label{eq:circle-character-family}
\mathcal T_{t'}(\TQ^{\,4})=\{\mathbf1,z\mapsto z,\ldots,z\mapsto z^{t'}\},\qquad
\mathcal T_t(\TQ)=\{\mathbf1,z\mapsto z,\ldots,z\mapsto z^t\},
\end{equation}
with values in $\Q(\mathrm i)^\times$.
The inclusion character is injective on each group, so each family has trivial common kernel, as required in \cref{def:abelian-pte}.

Let ${\mathsf{p}}_1<{\mathsf{p}}_2<\cdots$ be the rational primes congruent to $1\pmod4$.
For each such prime, choose a Gaussian prime ${p_{{\mathsf{p}}}}\in\Z[\mathrm i]$ with ${p_{{\mathsf{p}}}}\overline{{p_{{\mathsf{p}}}}}={\mathsf{p}}$, fixing ${p_{5}}=2+\mathrm i$ and ${p_{17}}=4+\mathrm i$.
The explicit constructions in \cref{app:certificate-2304,app:certificate-4608-9216} use these same choices.
Define
\begin{equation}\label{eq:circle-prime-ratios}
g_j:=\frac{{p_{{\mathsf{p}}_j}}}{\overline{{p_{{\mathsf{p}}_j}}}}\in\TQ,\qquad g_j^4\in \TQ^{\,4}\qquad(j\geq1).
\end{equation}
For an integer ${\mathsf{n}}\geq1$, put
\begin{equation}\label{eq:circle-common-denominators}
{\mathsf{b}}_{\mathsf{n}}:=\prod_{j=1}^{{\mathsf{n}}}{\mathsf{p}}_j.
\end{equation}
For an integer vector $\mathbf k=(k_1,\ldots,k_{\mathsf{n}})\in\Z^{\mathsf{n}}$, define
\begin{equation}\label{eq:circle-products}
{g_{\mathbf k}}:=\prod_{j=1}^{{\mathsf{n}}}g_j^{k_j}\in\TQ.
\end{equation}
The homomorphism $\Z^{\mathsf{n}}\to \TQ^{\,4}$ given by $\mathbf k\mapsto {g_{\mathbf k}}^4$ is injective: the Gaussian-prime valuation $\operatorname{ord}_{({p_{{\mathsf{p}}_j}})}$ takes the value $4k_j$ on ${g_{\mathbf k}}^4$.

Put
\begin{equation}\label{eq:circle-exponent-set}
J_{\mathsf{n}}:=\left\{\mathbf k\in\{-1,0,1\}^{\mathsf{n}}\setminus\{\mathbf0\}:\exists j\in\{1,\ldots,{\mathsf{n}}\}:k_j=1,\ k_i=0\ (1\leq i<j)\right\}.
\end{equation}
For each nonzero $\mathbf k\in\{-1,0,1\}^{\mathsf{n}}$, exactly one of $\mathbf k$ and $-\mathbf k$ belongs to $J_{\mathsf{n}}$, so $|J_{\mathsf{n}}|=(3^{\mathsf{n}}-1)/2$.
Injectivity shows that the two-element sets $\{{g_{\mathbf k}}^4,{g_{\mathbf k}}^{-4}\}$, $\mathbf k\in J_{\mathsf{n}}$, are pairwise disjoint.
Since $g_j^{\pm4k}\in {\mathsf{p}}_j^{-4k}\Z[\mathrm i]$, one has ${g_{\mathbf k}}^{4k}\in {\mathsf{b}}_{\mathsf{n}}^{-4k}\Z[\mathrm i]$ for every integer $k\geq1$ and every $\mathbf k\in J_{\mathsf{n}}$.

The following notion plays a key role in our subsequent arguments.
\begin{definition}[Sunflower]\label{def:circle-sunflower}
Following \cite{ErdosRado1960}, let $q\geq2$ be an integer and let $\mathfrak F_1,\ldots,\mathfrak F_q$ be pairwise distinct finite sets.
The family is a \emph{sunflower with core $\mathfrak K$} if
\[
\mathfrak F_i\cap \mathfrak F_j=\mathfrak K\qquad(1\leq i<j\leq q).
\]
The sets $\mathfrak F_i\setminus \mathfrak K$, $1\leq i\leq q$, are its \emph{petals}.
\end{definition}
Write
\[
\binom{J_{\mathsf{n}}}{{\mathsf{n}}}:=\{\mathfrak F\subseteq J_{\mathsf{n}}:|\mathfrak F|={\mathsf{n}}\}.
\]
For $\mathfrak F\subseteq J_{\mathsf{n}}$, define
\[
{\lambda_{t'}}:\{\mathfrak F:\mathfrak F\subseteq J_{\mathsf{n}}\}\to\Q(\mathrm i)^{t'},\qquad{\lambda_{t'}}(\mathfrak F):=
\left(\sum_{\mathbf k\in\mathfrak F}{g_{\mathbf k}}^{4k}\right)_{k=1}^{t'}.
\]
Since ${g_{\mathbf k}}^{-4k}=\overline{{g_{\mathbf k}}^{4k}}$, we have
\[
2\Re{\lambda_{t'}}(\mathfrak F)
=\left(\sum_{\mathbf k\in\mathfrak F}({g_{\mathbf k}}^{4k}+{g_{\mathbf k}}^{-4k})\right)_{k=1}^{t'}\in\Q^{t'}.
\]

\begin{lemma}\label{lem:circle-moment-sunflowers}
Fix integers $t'\geq1$ and $q\geq2$.
For all sufficiently large integers ${\mathsf{n}}$, there exist two sunflowers, each with $q$ members,
\[
\mathfrak F_1,\ldots,\mathfrak F_q\in\binom{J_{\mathsf{n}}}{{\mathsf{n}}},\qquad
\mathfrak F'_1,\ldots,\mathfrak F'_q\in\binom{J_{\mathsf{n}}}{{\mathsf{n}}},
\]
such that
\[
{\lambda_{t'}}(\mathfrak F_1)=\cdots={\lambda_{t'}}(\mathfrak F_q),\qquad
2\Re{\lambda_{t'}}(\mathfrak F'_1)=\cdots=2\Re{\lambda_{t'}}(\mathfrak F'_q).
\]
\end{lemma}
\begin{proof}
For each $\mathfrak F\in\binom{J_{\mathsf{n}}}{{\mathsf{n}}}$ and $1\leq k\leq t'$, we have
\[
\sum_{\mathbf k\in\mathfrak F}{g_{\mathbf k}}^{4k}\in {\mathsf{b}}_{\mathsf{n}}^{-4k}\Z[\mathrm i],
\qquad
\left|\sum_{\mathbf k\in\mathfrak F}{g_{\mathbf k}}^{4k}\right|\leq {\mathsf{n}}.
\]
The real and imaginary parts lie in ${\mathsf{b}}_{\mathsf{n}}^{-4k}\Z\cap[-{\mathsf{n}},{\mathsf{n}}]$, so the $k$-th coordinate of ${\lambda_{t'}}(\mathfrak F)$ takes $O({\mathsf{n}}^2{\mathsf{b}}_{\mathsf{n}}^{8k})$ possible values.
The $k$-th coordinate of $2\Re{\lambda_{t'}}(\mathfrak F)$ lies in ${\mathsf{b}}_{\mathsf{n}}^{-4k}\Z\cap[-2{\mathsf{n}},2{\mathsf{n}}]$, so it takes $O({\mathsf{n}}{\mathsf{b}}_{\mathsf{n}}^{4k})$ possible values.
By the prime number theorem for arithmetic progressions (see, e.g., \cite[Chapter~20]{Davenport2000}),
\[
\log {\mathsf{b}}_{\mathsf{n}}=\sum_{j=1}^{\mathsf{n}}\log {\mathsf{p}}_j=O({\mathsf{n}}\log {\mathsf{n}}).
\]
Since $t'$ is fixed, each vector takes at most $\exp(O_{t'}({\mathsf{n}}\log {\mathsf{n}}))$ distinct values as $\mathfrak F$ varies.

On the other hand,
\[
\log\binom{(3^{\mathsf{n}}-1)/2}{{\mathsf{n}}}=(\log3){\mathsf{n}}^2+O({\mathsf{n}}\log {\mathsf{n}}),
\qquad
\log\bigl({\mathsf{n}}!(q-1)^{\mathsf{n}}\bigr)=O_q({\mathsf{n}}\log {\mathsf{n}}).
\]
For all sufficiently large ${\mathsf{n}}$, the pigeonhole principle gives more than ${\mathsf{n}}!(q-1)^{\mathsf{n}}$ members of $\binom{J_{\mathsf{n}}}{{\mathsf{n}}}$ with the same value of ${\lambda_{t'}}$, and more than ${\mathsf{n}}!(q-1)^{\mathsf{n}}$ members with the same value of $2\Re{\lambda_{t'}}$.
Applying the Erd\H{o}s--Rado sunflower lemma \cite{ErdosRado1960} separately to these two collections gives the required sunflowers.
\end{proof}

\begin{proposition} \label{prop:circle-moment-pte}
Let $q\geq2$ be an integer.
Fix a sufficiently large integer ${\mathsf{n}}$, choose sunflowers $\mathfrak F_1,\ldots,\mathfrak F_q$ and $\mathfrak F'_1,\ldots,\mathfrak F'_q$ as in \cref{lem:circle-moment-sunflowers}, and let $\mathfrak K$ and $\mathfrak K'$ be their respective cores.
Define
\[
A_i:=\{{g_{\mathbf k}}^4:\mathbf k\in \mathfrak F_i\setminus \mathfrak K\},
\qquad
A'_i:=\bigcup_{\mathbf k\in \mathfrak F'_i\setminus \mathfrak K'}\{{g_{\mathbf k}}^4,{g_{\mathbf k}}^{-4}\}
\qquad(1\leq i\leq q).
\]
Then $(A_1,\ldots,A_q)$ and $(A'_1,\ldots,A'_q)$ are $(\TQ^{\,4},\mathcal T_{t'}(\TQ^{\,4}))$-PTE set solutions.
\end{proposition}
\begin{proof}
The petals in each sunflower are pairwise disjoint and have common positive size, respectively ${\mathsf{n}}-|\mathfrak K|$ and ${\mathsf{n}}-|\mathfrak K'|$, because its members are distinct ${\mathsf{n}}$-element sets.
Injectivity of $\mathbf k\mapsto {g_{\mathbf k}}^4$ and disjointness of the inverse pairs give the same properties for each family of blocks.
Subtracting the corresponding core sums in \cref{lem:circle-moment-sunflowers} gives, for $1\leq k\leq t'$ and $1\leq i\leq q$,
\[
\psi_k(A_i)=\sum_{\mathbf k\in \mathfrak F_i\setminus \mathfrak K}{g_{\mathbf k}}^{4k}=\psi_k(A_1),
\qquad \psi_k(A'_i)=\sum_{\mathbf k\in \mathfrak F'_i\setminus \mathfrak K'}({g_{\mathbf k}}^{4k}+{g_{\mathbf k}}^{-4k})=\psi_k(A'_1).\qedhere
\]
\end{proof}
\begin{proposition}\label{prop:circle-pte-lifting}
Recall the fourth-power map $\rho(z)=z^4$ in \eqref{eq:circle-fourth-power-map}.
Let $q\geq2$ be an integer.
If $(A_1,\ldots,A_q)$ is a $(\TQ^{\,4},\mathcal T_{t'}(\TQ^{\,4}))$-PTE set solution, then
\[
\bigl(\rho^{-1}(A_1),\ldots,\rho^{-1}(A_q)\bigr)
\]
is a $(\TQ,\mathcal T_t(\TQ))$-PTE set solution whose blocks are $R_4$-invariant.
\end{proposition}

\begin{proof}
Since every fiber of $\rho$ has four points, the inverse images are nonempty, pairwise disjoint, and have common size $4|A_i|$ for $1\leq i\leq q$.
Using $\psi_k$ as defined in \eqref{eq:Pk}, summation over the fibers of $\rho$ gives
\[
\psi_k\bigl(\rho^{-1}(A_i)\bigr)=
\begin{cases}
0,&4\nmid k,\\
4\psi_{k/4}(A_i),&4\mid k
\end{cases}
\qquad(1\leq k\leq t).
\]
For $4\mid k\leq t$, one has $k/4\leq t'$, so these moments are independent of $i$.
Each inverse image is $R_4$-invariant because $\ker\rho=R_4$.
\end{proof}

See~\cite[Sections~4 and~7]{SpiranecZenodo2026} for an earlier use of the sunflower lemma to construct rational spherical designs on the circle and for an exact-arithmetic verification of a rational spherical $7$-design with $9216$ points.

\subsection{Local Divisibility on Rational Ellipses and the Unit Circle}\label{subsec:norm-one-local-divisibility}\label{subsec:circle-local-divisibility}

Let $\ell$ be a positive squarefree integer and put
\[
\mathcal E_\ell:=\{(x_1,x_2)\in\R^2:x_1^2+\ell x_2^2=1\}.
\]
We first derive divisibility conditions for the size of rational designs on $\mathcal E_\ell$, obtaining the exponential lower bound in \cref{thm:conic-exponential-lower-bound}. We also show that the sum of two possible sizes is again a possible size. For $\ell=1$, additional $2$-adic congruences strengthen the divisibility result in \cref{thm:low-strength-divisibility}.

Put $\mathbb K_\ell:=\Q(\sqrt{-\ell})$, let \(\overline{\phantom z}\) denote its nontrivial $\Q$-automorphism, and use the embedding $\mathbb K_\ell\hookrightarrow\C$ determined by $\sqrt{-\ell}\mapsto\mathrm{i}\sqrt \ell$.
The group $T_{\mathbb K_\ell}$ defined in \cref{lem:norm-one-condition-D} is the group of $\Q$-points of the norm-one torus associated with $\mathbb K_\ell/\Q$ (see~\cite{HoshiKanaiYamasaki2022}).
Define a probability measure $\mu_{\mathcal E_\ell}$ on $\mathcal E_\ell$ by requiring, for every continuous $\varphi:\mathcal E_\ell\to\R$, that
\[
\int_{\mathcal E_\ell}\varphi\,\d\mu_{\mathcal E_\ell}
:=\int_{\mathbb S^1}\varphi(x_1,x_2/\sqrt\ell)\,\d\mu_{\mathbb S^1}(x_1,x_2),
\]
where $\mu_{\mathbb S^1}$ is normalized arc-length measure on $\mathbb S^1$.
If $\ell\neq1$, then $\mu_{\mathcal E_\ell}$ is not normalized arc-length measure.
Define
\[
\eta_\ell:\mathcal E_\ell\to\C,\qquad
\eta_\ell(x_1,x_2):=x_1+x_2\sqrt{-\ell}.
\]
Its restriction to $\mathcal E_\ell(\Q):=\mathcal E_\ell\cap\Q^2$ is a bijection onto $T_{\mathbb K_\ell}$; we identify these sets via this restriction.

\begin{lemma}\label{lem:conic-moment-criterion}
Let $t\geq1$ be an integer.
A nonempty finite set $\mathfrak X\subseteq\mathcal E_\ell(\Q)$ is a $\mathcal P_t(\mathcal E_\ell;\R)$-design on $(\mathcal E_\ell,\mu_{\mathcal E_\ell})$ if and only if
\[
\sum_{\mathbf x\in\mathfrak X}\eta_\ell(\mathbf x)^k=0\qquad(1\leq k\leq t).
\]
\end{lemma}
\begin{proof}
On $\mathcal E_\ell$, one has $\eta_\ell\overline{\eta_\ell}=\mathbf1$ and
\[
x_1=\frac{\eta_\ell(\mathbf x)+\overline{\eta_\ell}(\mathbf x)}2,
\qquad
x_2=\frac{\eta_\ell(\mathbf x)-\overline{\eta_\ell}(\mathbf x)}{2\sqrt{-\ell}}.
\]
These identities show that the $\C$-linear span of $\mathcal P_t(\mathcal E_\ell;\R)$ is generated by $\mathbf1$, $\eta_\ell^k$, and $\overline{\eta_\ell}^{\,k}$ for $1\leq k\leq t$.
The defining integral identity and rotation invariance of $\mu_{\mathbb S^1}$ give
\[
\int_{\mathcal E_\ell}\eta_\ell^k\,\d\mu_{\mathcal E_\ell}
=\int_{\mathbb S^1}(x_1+\mathrm i x_2)^k\,\d\mu_{\mathbb S^1}(x_1,x_2)=0
\qquad(k\geq1).
\]
The conjugate identities give the criterion.
\end{proof}

The linear map $(x_1,x_2)\mapsto(x_1,x_2/\sqrt\ell)$ from $\mathbb S^1$ to $\mathcal E_\ell$ conjugates circle isometries to linear transformations preserving $x_1^2+\ell x_2^2$.
Hence $\mu_{\mathcal E_\ell}$ is the invariant measure used in \cref{thm:rational-ellipsoidal-existence}.
Since $(1,0)\in \mathcal E_\ell(\Q)$, \cref{thm:rational-ellipsoidal-existence} with $M=\operatorname{diag}(1,\ell)$ and $r=1$ gives rational designs of every strength on $(\mathcal E_\ell,\mu_{\mathcal E_\ell})$, avoiding any prescribed finite subset of $\mathcal E_\ell(\Q)$.
\begin{definition}
For an integer $t\geq1$, define
\begin{equation}\label{eq:conic-cardinalities}
\mathcal N_{\ell,\Q}(t)
:=
\left\{\,|\mathfrak X|:
\begin{array}{l}
\varnothing\ne\mathfrak X\subseteq\mathcal E_\ell(\Q)\text{ is finite,}\\
\mathfrak X\text{ is a }\mathcal P_t(\mathcal E_\ell;\R)\text{-design}\\
\text{on }(\mathcal E_\ell,\mu_{\mathcal E_\ell})
\end{array}
\,\right\}.
\end{equation}
\end{definition}

Write $\mathcal O_{\mathbb K_\ell}$ for the ring of integers of $\mathbb K_\ell$.
For every prime ${\mathsf p}$, put
\[
\Z_{({\mathsf p})}:=\{z/{\mathsf b}:z\in\Z,\ {\mathsf b}\in\Z_{>0},\ {\mathsf p}\nmid{\mathsf b}\}.
\]
An odd prime ${\mathsf{p}}$ is \emph{inert in $\mathbb K_\ell$} if ${\mathsf{p}}\mathcal O_{\mathbb K_\ell}$ is a prime ideal of $\mathcal O_{\mathbb K_\ell}$.
\Needspace{5\baselineskip}
The next lemma gives a divisibility condition for the vanishing of a single moment.

\begin{lemma}\label{lem:inert-prime-moment-divisibility}
Let ${\mathsf{p}}$ be an odd prime inert in $\mathbb K_\ell$, let ${\mathsf{v}}\geq1$ be an integer, and put $k:=({\mathsf{p}}+1){\mathsf{p}}^{{\mathsf{v}}-1}$.
Let $I\subseteq\mathbb N$ be a nonempty finite position set and let $(\mathbf x^{(j)})_{j\in I}\in\mathcal E_\ell(\Q)^I$ be a finite sequence, with repetitions allowed.
If
\[
\sum_{j\in I}\eta_\ell(\mathbf x^{(j)})^k=0,
\]
then ${\mathsf{p}}^{2{\mathsf{v}}}\mid|I|$.
\end{lemma}
\begin{proof}
For $\mathbf x=(x_1,x_2)\in\mathcal E_\ell(\Q)$, choose ${\mathsf{b}}\in\Z_{>0}$ with ${\mathsf{b}}x_1,{\mathsf{b}}x_2\in\Z$ and $\gcd({\mathsf{b}}x_1,{\mathsf{b}}x_2,{\mathsf{b}})=1$.
Then $({\mathsf{b}}x_1)^2+\ell({\mathsf{b}}x_2)^2={\mathsf{b}}^2$.
Since ${\mathsf{p}}$ is inert, $-\ell$ is a quadratic nonresidue modulo ${\mathsf{p}}$.
If ${\mathsf{p}}\mid {\mathsf{b}}$, the identity would force ${\mathsf{p}}\mid {\mathsf{b}}x_1,{\mathsf{b}}x_2$, contrary to primitivity.
Thus ${\mathsf{p}}\nmid {\mathsf{b}}$.
Put
\[
\mathcal O_{\mathsf{p}}:=\Z_{({\mathsf{p}})}[\sqrt{-\ell}]=\{x_1+x_2\sqrt{-\ell}:x_1,x_2\in\Z_{({\mathsf{p}})}\}.
\]
Because ${\mathsf{p}}$ is odd and $[\mathcal O_{\mathbb K_\ell}:\Z[\sqrt{-\ell}]]\in\{1,2\}$, multiplication gives the identification
\[
\mathcal O_{\mathbb K_\ell}\otimes_\Z\Z_{({\mathsf{p}})}\xrightarrow{\sim}\mathcal O_{\mathsf{p}},
\qquad z\otimes a\longmapsto az.
\]
Both $z:=\eta_\ell(\mathbf x)$ and $z^{-1}=\overline z$ belong to $\mathcal O_{\mathsf{p}}$.
Since ${\mathsf{p}}$ is inert, $\mathcal O_{\mathsf{p}}/{\mathsf{p}}\mathcal O_{\mathsf{p}}\cong\mathbb F_{{\mathsf{p}}^2}$, and its nontrivial $\mathbb F_{\mathsf{p}}$-automorphism is $z\mapsto z^{\mathsf{p}}$ \cite[Chapter~12, Section~3, pp.~183--184]{IrelandRosen1990}.
Hence $z^{\mathsf{p}}\equiv\overline z\pmod{{\mathsf{p}}\mathcal O_{\mathsf{p}}}$ and $z^{{\mathsf{p}}+1}\equiv1\pmod{{\mathsf{p}}\mathcal O_{\mathsf{p}}}$.
For ${\mathsf{v}}'\geq1$ and $z'\in\mathcal O_{\mathsf{p}}$ with $z'\equiv1\pmod{{\mathsf{p}}^{{\mathsf{v}}'}\mathcal O_{\mathsf{p}}}$, the binomial theorem gives $(z')^{\mathsf{p}}\equiv1\pmod{{\mathsf{p}}^{{\mathsf{v}}'+1}\mathcal O_{\mathsf{p}}}$.
Applying this implication ${\mathsf{v}}-1$ times gives $z^k\equiv1\pmod{{\mathsf{p}}^{\mathsf{v}}\mathcal O_{\mathsf{p}}}$.
Consequently,
\[
\frac{z^k+z^{-k}}2-1=\frac{(z^k-1)^2}{2z^k}\in {\mathsf{p}}^{2{\mathsf{v}}}\mathcal O_{\mathsf{p}}.
\]
The left-hand side is fixed by conjugation and $\mathcal O_{\mathsf{p}}\cap\Q=\Z_{({\mathsf{p}})}$, so $(z^k+z^{-k})/2\equiv1\pmod{{\mathsf{p}}^{2{\mathsf{v}}}\Z_{({\mathsf{p}})}}$.
Taking the real part of the assumed identity therefore yields
\[
0=\sum_{j\in I}\frac{\eta_\ell(\mathbf x^{(j)})^k+\eta_\ell(\mathbf x^{(j)})^{-k}}2
\equiv|I|\pmod{{\mathsf{p}}^{2{\mathsf{v}}}\Z_{({\mathsf{p}})}}.
\]
Since ${\mathsf{p}}^{2{\mathsf{v}}}\Z_{({\mathsf{p}})}\cap\Z={\mathsf{p}}^{2{\mathsf{v}}}\Z$, this proves ${\mathsf{p}}^{2{\mathsf{v}}}\mid|I|$.
\end{proof}

For an integer $t\geq1$, put
\begin{equation}\label{eq:conic-divisor}
{\delta_{\ell}}(t):=\prod_{\substack{{\mathsf{p}}\text{ odd inert in }\mathbb K_\ell\\{\mathsf{p}}+1\leq t}}
 {\mathsf{p}}^{2\max\{{\mathsf{v}}\in\Z_{\geq1}:({\mathsf{p}}+1){\mathsf{p}}^{{\mathsf{v}}-1}\leq t\}}.
\end{equation}
Applying \cref{lem:inert-prime-moment-divisibility} to the moment conditions in \cref{lem:conic-moment-criterion} gives the following divisibility condition.

\begin{corollary}\label{cor:conic-local-divisibility}
For every integer $t\geq1$, every rational $\mathcal P_t(\mathcal E_\ell;\R)$-design $\mathfrak X\subseteq \mathcal E_\ell(\Q)$ satisfies
\[
{\delta_{\ell}}(t)\mid|\mathfrak X|.
\]
In particular, ${\delta_{\ell}}(t)\mid \min\mathcal N_{\ell,\Q}(t)$.
\end{corollary}
\begin{proof}
For each odd prime ${\mathsf{p}}$ inert in $\mathbb K_\ell$ with ${\mathsf{p}}+1\leq t$, put
\[
{\mathsf{v}}:=\max\{{\mathsf{v}}'\in\Z_{\geq1}:({\mathsf{p}}+1){\mathsf{p}}^{{\mathsf{v}}'-1}\leq t\}.
\]
The moment of degree $({\mathsf{p}}+1){\mathsf{p}}^{{\mathsf{v}}-1}$ vanishes by \cref{lem:conic-moment-criterion}, so \cref{lem:inert-prime-moment-divisibility} gives ${\mathsf{p}}^{2{\mathsf{v}}}\mid|\mathfrak X|$.
The factors in ${\delta_{\ell}}(t)$ are powers of distinct primes, so their product divides $|\mathfrak X|$.
\end{proof}

\Needspace{5\baselineskip}
\begin{theorem}\label{thm:conic-exponential-lower-bound}
For each fixed positive squarefree integer $\ell$, as $t\to\infty$,
\[
\log\bigl(\min\mathcal N_{\ell,\Q}(t)\bigr)\geq t+\operatorname{o}_\ell(t),
\]
where the error term may depend on $\ell$ and satisfies $\operatorname{o}_\ell(t)/t\to0$.
\end{theorem}
\begin{proof}
The prime number theorem for arithmetic progressions, 
applied to the quadratic Dirichlet character associated with $\mathbb K_\ell$,
gives
\[
\sum_{\substack{{\mathsf{p}}\leq t-1\\{\mathsf{p}}\text{ odd inert in }\mathbb K_\ell}}\log {\mathsf{p}}
=\frac{t}{2}+\operatorname{o}_\ell(t).
\]
Every prime in this sum contributes a factor divisible by ${\mathsf{p}}^2$ to ${\delta_{\ell}}(t)$, so \cref{cor:conic-local-divisibility} gives
\[
\log\bigl(\min\mathcal N_{\ell,\Q}(t)\bigr)\geq\log{\delta_{\ell}}(t)
\geq2\sum_{\substack{{\mathsf{p}}\leq t-1\\{\mathsf{p}}\text{ odd inert in }\mathbb K_\ell}}\log {\mathsf{p}}
=t+\operatorname{o}_\ell(t).\qedhere
\]
\end{proof}

\Needspace{5\baselineskip}
\begin{proposition}\label{prop:conic-cardinality-additivity}
For every integer $t\geq1$, the set $\mathcal N_{\ell,\Q}(t)\cup\{0\}$ is closed under addition.
In particular, every rational $t$-design $\mathfrak X\subseteq\mathcal E_\ell(\Q)$ satisfies
\[
|\mathfrak X|\Z_{>0}\subseteq\mathcal N_{\ell,\Q}(t).
\]
\end{proposition}
\begin{proof}
The elements $(1+j\sqrt{-\ell})/(1-j\sqrt{-\ell})$, $j\in\Z_{>0}$, belong to $T_{\mathbb K_\ell}$ and are pairwise distinct, so $T_{\mathbb K_\ell}$ is infinite.
Addition by $0$ is immediate.
For nonempty finite $\mathfrak X,\mathfrak X'\subseteq T_{\mathbb K_\ell}$, choose
\[
{v}\in T_{\mathbb K_\ell}\setminus\{z'z^{-1}:z\in\mathfrak X,\ z'\in\mathfrak X'\}.
\]
Then ${v}\mathfrak X\cap\mathfrak X'=\varnothing$.
If $\psi_k(\mathfrak X)=\psi_k(\mathfrak X')=0$ for $1\leq k\leq t$, then
\[
\psi_k({v}\mathfrak X\sqcup \mathfrak X')
={v}^k\psi_k(\mathfrak X)+\psi_k(\mathfrak X')=0
\qquad(1\leq k\leq t).
\]
Thus ${v}\mathfrak X\sqcup \mathfrak X'$ satisfies the moment criterion in \cref{lem:conic-moment-criterion}, so the sizes add.
Induction gives every positive multiple.
\end{proof}

Now, for $\ell=1$, we derive sharper divisibility information from $2$-adic congruences, complementing \cref{cor:conic-local-divisibility}.
Below $\Z_2$ denotes the ring of $2$-adic integers.

\Needspace{5\baselineskip}
\begin{lemma}\label{lem:primitive-pythagorean-local-form}
Let $z=x_1+\mathrm i x_2\in\TQ$, and choose ${\mathsf{b}}\in\Z_{>0}$ such that
\[
{\mathsf{b}}x_1,{\mathsf{b}}x_2\in\Z,\qquad\gcd({\mathsf{b}}x_1,{\mathsf{b}}x_2,{\mathsf{b}})=1,
\qquad({\mathsf{b}}x_1)^2+({\mathsf{b}}x_2)^2={\mathsf{b}}^2.
\]
Then ${\mathsf{b}}$ is odd, and one of ${\mathsf{b}}x_1,{\mathsf{b}}x_2$ is odd while the other is divisible by $4$.
\end{lemma}
\begin{proof}
Reduce $({\mathsf{b}}x_1)^2+({\mathsf{b}}x_2)^2={\mathsf{b}}^2$ modulo $4$ and $8$; see also \cite[Chapter~17, Section~1, p.~271]{IrelandRosen1990}.
If ${\mathsf{b}}$ were even, reduction modulo $4$ would force both ${\mathsf{b}}x_1$ and ${\mathsf{b}}x_2$ to be even, contradicting primitivity.
Thus ${\mathsf{b}}$ is odd, and exactly one numerator is odd.
Every odd square is $1$ modulo $8$, so the even numerator has square $0$ modulo $8$ and is divisible by $4$.
\end{proof}

For $z=x_1+\mathrm i x_2\in\TQ$, define the functions $\gamma_4,\gamma_5,\gamma_6:\TQ\to\Q$ by
\begin{equation}\label{eq:circle-congruence-functions}
\begin{aligned}
\gamma_4(z):={}&1-\Re(z^4)-8\Im(z^4),\\
\gamma_5(z):={}&2+2x_1-2x_2-9\Re(z^3)+7\Im(z^3)-2\Re(z^4)+16\Im(z^4)+7\Re(z^5)+9\Im(z^5),\\
\gamma_6(z):={}&2-352\Re(z)+452\Im(z)+511\Re(z^2)-928\Im(z^2)-672\Re(z^3)+450\Im(z^3)\\
&+510\Re(z^4)+16\Im(z^4)-2\Im(z^5)+\Re(z^6).
\end{aligned}
\end{equation}

\begin{lemma}\label{lem:circle-two-adic-congruences}
For every $z\in\TQ$,
\[
\gamma_4(z)\in2^8\Z_2,\qquad
\gamma_5(z)\in2^{10}\Z_2,\qquad
\gamma_6(z)\in2^{11}\Z_2.
\]
\end{lemma}
\begin{proof}
Write $z=x_1+\mathrm i x_2$ and choose the primitive denominator ${\mathsf{b}}$ as in \cref{lem:primitive-pythagorean-local-form}.
Then $x_1,x_2\in\Z_2$, with one odd and the other in $4\Z_2$.
Choose $i\in\{1,2\}$ so that $x_i\in4\Z_2$.
Direct expansion gives
\[
\Re(z^4)=1-128(x_i/4)^2x_{3-i}^2,\qquad
\frac{\Im(z^4)}{16}=(x_i/4)x_{3-i}(x_1^2-x_2^2)
\equiv x_i/4\equiv(x_i/4)^2x_{3-i}^2\pmod2.
\]
Consequently,
\[
\gamma_4(z)=128\left((x_i/4)^2x_{3-i}^2-\frac{\Im(z^4)}{16}\right)\in2^8\Z_2.
\]
Introduce the following formal polynomials in $\Q[X_1,X_2]$:
\[
\begin{aligned}
{\vartheta_{5,2}}:={}&144(X_2/4)^3+28X_1(X_2/4)^2+32(X_2/4)^2-X_1X_2/4-5X_2/4-X_1-1,\\
{\vartheta_{6,2}}:={}&128(X_2/4)^6+32(X_2/4)^5-1032(X_2/4)^4+8X_1(X_2/4)^3\\
&+110(X_2/4)^3-42X_1(X_2/4)^2+80(X_2/4)^2+7X_1X_2/4-7X_2/4+X_1-1,\\
{\vartheta_{5,1}}:={}&112(X_1/4)^3+36X_2(X_1/4)^2-32(X_1/4)^2-X_1X_2/4-3X_1/4-X_2+1,\\
{\vartheta_{6,1}}:={}&-64(X_1/4)^4-16(X_1/4)^3+4X_2(X_1/4)^2-508(X_1/4)^2
-3X_1X_2/4-43X_1/4-15X_2+13.
\end{aligned}
\]
If $x_2\in4\Z_2$, expansion using $x_1^2+x_2^2=1$ gives
\[
\gamma_5(z)=256(x_2/4)(x_2-1){\vartheta_{5,2}}(x_1,x_2),\qquad
\gamma_6(z)=-1024{\vartheta_{6,2}}(x_1,x_2),
\]
and
\[
{\vartheta_{5,2}}(x_1,x_2)\equiv-(x_2/4+1)(x_1+1)\pmod4, \qquad {\vartheta_{6,2}}(x_1,x_2)\equiv(x_2/4+1)(x_1+1)\equiv0\pmod2.
\]
If $x_1\in4\Z_2$, the corresponding identities are
\[
\gamma_5(z)=256(x_1/4)(x_1+1){\vartheta_{5,1}}(x_1,x_2),\qquad
\gamma_6(z)=-512(x_1/4)(x_1-1){\vartheta_{6,1}}(x_1,x_2),
\]
and
\[
{\vartheta_{5,1}}(x_1,x_2)\equiv(x_1/4+1)(1-x_2)\pmod4, \qquad {\vartheta_{6,1}}(x_1,x_2)\equiv(x_1/4+1)(x_2+1)\pmod4.
\]
In the first case $x_1$ is odd, so $(x_2/4)(x_2/4+1)\in2\Z_2$ and $x_1+1\in2\Z_2$ imply
$(x_2/4){\vartheta_{5,2}}(x_1,x_2)\in4\Z_2$ and ${\vartheta_{6,2}}(x_1,x_2)\in2\Z_2$.
In the second case $x_2$ is odd, so $(x_1/4)(x_1/4+1)\in2\Z_2$ and $1\pm x_2\in2\Z_2$ imply
$(x_1/4){\vartheta_{5,1}}(x_1,x_2)\in4\Z_2$ and $(x_1/4){\vartheta_{6,1}}(x_1,x_2)\in4\Z_2$.
The displayed factorizations prove the assertions for $\gamma_5$ and $\gamma_6$.
\end{proof}

The congruence for $\gamma_4$ yields the following divisibility result from a vanishing fourth moment.

\begin{proposition}\label{prop:fourth-moment-divisibility}
Let $I\subseteq\mathbb N$ be a nonempty finite position set, and suppose $z_j\in\TQ$ $(j\in I)$ satisfy $\sum_{j\in I}z_j^4=0$.
Then $2304\mid|I|$.
\end{proposition}
\begin{proof}
The assumption and \cref{lem:circle-two-adic-congruences} give
\[
|I|=\sum_{j\in I}\gamma_4(z_j)\in2^8\Z_2.
\]
Since $-1$ is not a square modulo $3$, the prime $3$ is inert in $\Q(\mathrm i)$.
Apply \cref{lem:inert-prime-moment-divisibility} with $\ell=1$, ${\mathsf{p}}=3$, ${\mathsf{v}}=1$, and the corresponding point sequence $(\Re z_j,\Im z_j)_{j\in I}$ to obtain $9\mid|I|$.
Together these give $2304\mid|I|$.
\end{proof}

For $t\geq4$, define
\begin{equation}\label{eq:circle-divisor}
{\delta_{\TQ}}(t)
:=
\begin{cases}
2^8{\delta_{1}}(t),&t=4,\\
2^9{\delta_{1}}(t),&t=5,\\
2^{10}{\delta_{1}}(t),&t\geq6.
\end{cases}
\end{equation}

\begin{theorem}\label{thm:low-strength-divisibility}
Let $t\geq4$ be an integer, and let $\mathfrak X\subseteq\TQ$ be a spherical $t$-design.
Then
\[
{\delta_{\TQ}}(t)\mid|\mathfrak X|.
\]
For \(4\leq t\leq11\), this gives explicitly
\[
{\delta_{\TQ}}(t)
=
\begin{cases}
2304=2^8\cdot3^2,&t=4,\\
4608=2^9\cdot3^2,&t=5,\\
9216=2^{10}\cdot3^2,&6\leq t\leq7,\\
451\,584=2^{10}\cdot3^2\cdot7^2,&8\leq t\leq11.
\end{cases}
\]
\end{theorem}

\begin{proof}
Applying \cref{cor:conic-local-divisibility} with $\ell=1$ gives ${\delta_{1}}(t)\mid |\mathfrak X|$.
Moreover, \cref{lem:circle-moment-criterion} and \eqref{eq:circle-congruence-functions} give
\[
\sum_{z\in \mathfrak X}\gamma_4(z)=|\mathfrak X|\quad(t\geq4),
\qquad
\sum_{z\in \mathfrak X}\gamma_5(z)=2|\mathfrak X|\quad(t\geq5),
\qquad
\sum_{z\in \mathfrak X}\gamma_6(z)=2|\mathfrak X|\quad(t\geq6).
\]
By \cref{lem:circle-two-adic-congruences}, we have $2^8\mid |\mathfrak X|$ for $t=4$, $2^9\mid |\mathfrak X|$ for $t=5$, and $2^{10}\mid |\mathfrak X|$ for $t\geq6$.
Since ${\delta_{1}}(t)$ is odd, \eqref{eq:circle-divisor} gives ${\delta_{\TQ}}(t)\mid |\mathfrak X|$.
The explicit values follow from \eqref{eq:conic-divisor}: ${\delta_{1}}(t)=3^2$ for $4\leq t\leq7$, and ${\delta_{1}}(t)=3^2\cdot7^2$ for $8\leq t\leq11$.
\end{proof}

\Needspace{5\baselineskip}
For spherical $4$-designs, the conditions $-\mathfrak X=\mathfrak X$ and $\mathrm i\mathfrak X=\mathfrak X$ impose the stronger divisibility conditions in the following corollary.

\begin{corollary}\label{cor:fourth-moment-symmetric-divisibility}
Let $\mathfrak X\subseteq\TQ$ be a spherical $4$-design.
\begin{enumerate}
    \item[(i)] If $-\mathfrak X=\mathfrak X$, then $4608\mid|\mathfrak X|$.
    \item[(ii)] If $\mathrm{i}\mathfrak X=\mathfrak X$, then $9216\mid|\mathfrak X|$.
\end{enumerate}
\end{corollary}
\begin{proof}
By \cref{lem:circle-moment-criterion}, $\psi_4(\mathfrak X)=0$, with $\psi_4$ as in \eqref{eq:Pk}.
Assume first that $-\mathfrak X=\mathfrak X$.
Choose one representative from each pair $\{z,-z\}$ and index these representatives by a finite position set $I$, writing them as $(z_j)_{j\in I}$.
Then
\[
\psi_4(\mathfrak X)=2\sum_{j\in I}z_j^4=0,\qquad|\mathfrak X|=2|I|.
\]
By \cref{prop:fourth-moment-divisibility}, $2304\mid|I|$, hence $4608\mid|\mathfrak X|$.
If $\mathrm i\mathfrak X=\mathfrak X$, choose instead one representative from each $R_4$-orbit, indexed by $I$.
Every orbit has four points with the same fourth power, so
\[
\psi_4(\mathfrak X)=4\sum_{j\in I}z_j^4=0,\qquad|\mathfrak X|=4|I|.
\]
Again $2304\mid|I|$, giving $9216\mid|\mathfrak X|$.
\end{proof}

For the PTE problem on the rational unit circle, we use the localization
\[
\Z_{(2)}\cap\Z_{(3)}=\{z/{\mathsf{b}}:z\in\Z,\ {\mathsf{b}}\in\Z_{>0},\ \gcd({\mathsf{b}},6)=1\},
\]
(see \cite[p.~18]{Milne2020}).
Write $\operatorname{Tr}:=\operatorname{Tr}_{\Q(\mathrm i)/\Q}$ for the field trace, that is,
\begin{equation}\label{eq:gaussian-field-trace}
\operatorname{Tr}(z):=z+\overline z\qquad(z\in\Q(\mathrm i)).
\end{equation}
For $z\in\TQ$, this equals $z+z^{-1}$.
Congruences below are taken in $\Z_{(2)}\cap\Z_{(3)}$.

\Needspace{5\baselineskip}
\begin{proposition}\label{prop:circle-fourth-image-trace-congruence}
The following hold:
\begin{enumerate}
\item[(a)] For every $z\in\TQ$,
\[
\Re(z^4)+24\Im(z^4)\equiv1\pmod{2304}.
\]
Consequently, if finite multisets $A_1,A_2$ in $\TQ$ satisfy $\psi_4(A_1)=\psi_4(A_2)$, then $|A_1|\equiv|A_2|\pmod{2304}$, where multiplicities are counted.
\item[(b)] For every $z\in \TQ^{\,4}$,
\[
\operatorname{Tr}(z)\equiv2\pmod{2304}.
\]
Consequently, for any nonempty finite position set \(I\subseteq\mathbb N\) and any \(z_j\in \TQ^{\,4}\) \((j\in I)\),
\[
\sum_{j\in I}\operatorname{Tr}(z_j)=0\quad\Longrightarrow\quad1152\mid|I|.
\]
\end{enumerate}
\end{proposition}
\begin{proof}
Write $z=x_1+\mathrm i x_2$ and take ${\mathsf{b}}$ as in \cref{lem:primitive-pythagorean-local-form}.
If $3\mid {\mathsf{b}}$, reduction of $({\mathsf{b}}x_1)^2+({\mathsf{b}}x_2)^2={\mathsf{b}}^2$ modulo $3$ would force both numerators to be divisible by $3$, contrary to primitivity.
Thus $3\nmid {\mathsf{b}}$, and exactly one numerator is divisible by $3$.
The same lemma gives $2\nmid {\mathsf{b}}$ and $4\mid({\mathsf{b}}x_1)({\mathsf{b}}x_2)$.
Therefore $({\mathsf{b}}x_1)({\mathsf{b}}x_2)/12$ is an integer, whereas $({\mathsf{b}}x_1)^2-({\mathsf{b}}x_2)^2$ is odd.
Expansion yields
\[
\Re(z^4)=1-\frac{1152}{{\mathsf{b}}^4}\left(\frac{({\mathsf{b}}x_1)({\mathsf{b}}x_2)}{12}\right)^2,
\qquad \Im(z^4)=\frac{48}{{\mathsf{b}}^4}\frac{({\mathsf{b}}x_1)({\mathsf{b}}x_2)}{12}\bigl(({\mathsf{b}}x_1)^2-({\mathsf{b}}x_2)^2\bigr).
\]
The product of an integer and its difference from an odd integer is even.
Combining these two formulas thus gives $\Re(z^4)+24\Im(z^4)-1\in2304(\Z_{(2)}\cap\Z_{(3)})$.
Since $2304(\Z_{(2)}\cap\Z_{(3)})\cap\Z=2304\Z$, summing over $A_1$ and $A_2$ proves part~\textup{(a)}.

For part~\textup{(b)}, write $z=(z')^4$ with $z'\in\TQ$ and apply the formula for the real part to $z'$.
It gives $\operatorname{Tr}(z)=2\Re((z')^4)\equiv2\pmod{2304}$.
Summation of this congruence gives $2|I|\in2304(\Z_{(2)}\cap\Z_{(3)})$, and hence $1152\mid|I|$.
\end{proof}

\begin{corollary}\label{cor:circle-integral-trace-relation}
Let $I\subseteq\mathbb N$ be a finite position set, let ${\mathsf{c}_{j}}\in\Z$, and let $z_j\in \TQ^{\,4}$ for $j\in I$.
If $\sum_{j\in I}{\mathsf{c}_{j}}\operatorname{Tr}(z_j)\in\Z$, then
\[
\sum_{j\in I}{\mathsf{c}_{j}}\bigl(\operatorname{Tr}(z_j)-2\bigr)\in2304\Z.
\]
\end{corollary}
\begin{proof}
Multiply the congruence in \cref{prop:circle-fourth-image-trace-congruence}(b) by ${\mathsf{c}_{j}}$ and sum.
The resulting sum belongs both to $2304(\Z_{(2)}\cap\Z_{(3)})$ and to $\Z$, whose intersection is $2304\Z$.
\end{proof}
\subsection{Explicit Constructions and Quantitative Bounds}\label{subsec:explicit-circle-constructions}

The constructions below attain the divisibility bounds in \cref{thm:low-strength-divisibility} for strengths $4$ through $7$ (see \cref{thm:minimum-circle-constructions}).
The explicit point sets, index lists, and proofs are given together in
\cref{app:certificate-2304,app:certificate-4608-9216}.
We then give an explicit upper bound for the size of a rational spherical $11$-design (\cref{prop:circle-strength-eleven-quantitative}).

For $z=x_1+\mathrm i x_2\in\TQ$ with $x_1,x_2>0$, consider the sets
\[
\{z,\overline z\},\qquad
\{z,-z,\overline z,-\overline z\},\qquad
\{{v} z,{v}\overline z:{v}\in R_4\}.
\]
These sets have two, four, and eight points, respectively, since $x_1=x_2$ would give the impossible rational identity $2x_1^2=1$.
The strength-$4$ construction uses the finite index sets
\[
L_4^{(2)}\subseteq\{-2185,\ldots,-1,1,\ldots,2185\},
\qquad L_4^{(4)},L_4^{(8)}\subseteq\{0,\ldots,2185\},
\]
specified by the lists in \cref{app:certificate-2304}.
They select $192$ conjugate pairs, $96$ antipodal four-point sets, and $192$ $R_4$-invariant eight-point sets, respectively.
The signs in $L_4^{(2)}$ select one of the two conjugate pairs in each corresponding antipodal four-point set.
The selected pairs contribute zero to the first and third power sums.
The four-point sets selected by $L_4^{(4)}$ have vanishing odd power sums, and their total second sum cancels that of the selected pairs.
The eight-point sets selected by $L_4^{(8)}$ contribute zero to the first, second, and third power sums, and their total fourth power sum cancels that of the selected pairs and four-point sets.
The points lie in ${\mathsf{b}}_7^{-1}\Z[\mathrm i]$ (see \eqref{eq:circle-seven-prime-denominator}); multiplying their $k$-th power sums by ${\mathsf{b}}_7^k$ and taking real and imaginary parts gives equations in $\Z$ for $1\leq k\leq4$.
\Cref{prop:certificate-four} verifies these equalities and shows that no point is selected twice.
Thus the point count is
\[
2\cdot192+4\cdot96+8\cdot192=2304.
\]

For strengths $5$ and $7$, the signed index set
\[
L_5\subseteq\{-3279,\ldots,-1,1,\ldots,3279\},\qquad |L_5|=1152,
\]
specified in \cref{app:certificate-4608-9216}, selects $1152$ disjoint $R_4$-invariant eight-point sets.
Each has vanishing moments of degrees $1,2,3,5,6,7$ by $R_4$-invariance; \cref{prop:certificate-five-seven} verifies that their union has vanishing fourth moment.
From each eight-point set, the sign in $L_5$ selects $\{\pm z,\pm\overline z\}$ or its $\mathrm i$-multiple, giving $4608$ points with total second power sum zero.
Each four-point subset has vanishing odd moments and half the fourth moment of the corresponding eight-point set.
The resulting sets are nested rational spherical $5$- and $7$-designs.

\begin{theorem}[Minimum rational circle designs]
\label{thm:minimum-circle-constructions}
There exist finite sets $\mathfrak X_4,\mathfrak X_5,\mathfrak X_7\subseteq\TQ$ satisfying the following:
\begin{enumerate}
\item[(i)]
$\mathfrak X_4$ is a rational spherical $4$-design with $2304$ distinct points.
\item[(ii)]
$\mathfrak X_5$ is an antipodal rational spherical $5$-design with $4608$ distinct points, and $\mathfrak X_7$ is an $R_4$-invariant rational spherical $7$-design with $9216$ distinct points.
\end{enumerate}
\end{theorem}
\begin{remark}\label{rem:exact-seven-coordinate-distinct}
A rational rotation of $\mathfrak X_7$ gives an $R_4$-invariant rational spherical design with $9216$ points, pairwise distinct first coordinates, and strength exactly $7$.
Indeed, $7$ is inert in $\Q(\mathrm i)$, so \cref{lem:inert-prime-moment-divisibility} with $\ell=1$, ${\mathsf{p}}=7$, and ${\mathsf{v}}=1$ would give $49\mid9216$ if $\psi_8(\mathfrak X_7)=0$, whereas $9216\equiv4\pmod{49}$.
Apply \cref{lem:circle-coordinate-separating-rotation}; rotation preserves $R_4$-invariance and the vanishing or nonvanishing of every moment.
\end{remark}

For $z,z'\in\TQ$ and $k\geq1$, the trace satisfies
\[
\operatorname{Tr}(z^k(z')^k)+\operatorname{Tr}(z^k(z')^{-k})=\operatorname{Tr}(z^k)\operatorname{Tr}((z')^k).
\]

\begin{proposition}\label{prop:circle-strength-eleven-quantitative}
There exists an $R_4$-invariant rational spherical $11$-design $\mathfrak X_{11}\subseteq\TQ$ satisfying
\[
|\mathfrak X_{11}|\leq20\,808\,990\,720\,000.
\]
\end{proposition}
\begin{proof}
The positive integers
\[
m_2=1\,835\,937\,500,\qquad m_3=6\,975\,757\,441
\]
satisfy, by direct evaluation with $g_1=(2+\mathrm i)/(2-\mathrm i)$ and $g_3=(4+\mathrm i)/(4-\mathrm i)$,
\begin{equation}\label{eq:circle-eleven-trace-relation}
8\,387\,848\,518+m_2\operatorname{Tr}(g_1^{4k})+m_3\operatorname{Tr}(g_3^{4k})=0
\qquad(k=1,2).
\end{equation}
Set $q:=8\,387\,848\,518+m_2+m_3=17\,199\,543\,459$.

Let $({\mathsf{p}}'_j)_{j\geq1}$ enumerate the primes in $({\mathsf{p}}_j)_{j\geq1}$ other than $5$ and $17$, in increasing order.
Thus ${\mathsf{p}}'_1=13$ and ${\mathsf{p}}'_{100}=1277$.
For $\mathbf k\in J_{100}$, put
\[
{g'_{\mathbf k}}:=\prod_{j=1}^{100}\left(\frac{{p_{{\mathsf{p}}'_j}}}{\overline{{p_{{\mathsf{p}}'_j}}}}\right)^{k_j}\in\TQ.
\]
The corresponding denominator product is
\[
{\mathsf{b}}'_{100}:=\prod_{j=1}^{100}{\mathsf{p}}'_j.
\]
For $\mathfrak F\in\binom{J_{100}}{100}$, consider the pair
\[
\left(\sum_{\mathbf k\in\mathfrak F}\operatorname{Tr}(({g'_{\mathbf k}})^4),\quad
\sum_{\mathbf k\in\mathfrak F}\operatorname{Tr}(({g'_{\mathbf k}})^8)\right).
\]
Since $({\mathsf{b}}'_{100})^4({g'_{\mathbf k}})^4\in\Z[\mathrm i]$ and $|{g'_{\mathbf k}}|=1$, the coordinates lie in $({\mathsf{b}}'_{100})^{-4}\Z\cap[-200,200]$ and $({\mathsf{b}}'_{100})^{-8}\Z\cap[-200,200]$, respectively.
Thus at most
\[
\bigl(400({\mathsf{b}}'_{100})^4+1\bigr)\bigl(400({\mathsf{b}}'_{100})^8+1\bigr)
\]
such pairs occur.
The exact integer inequality
\begin{equation}\label{eq:circle-eleven-sunflower-inequality}
\bigl(|J_{100}|-99\bigr)^{100} >\bigl(400({\mathsf{b}}'_{100})^4+1\bigr)\bigl(400({\mathsf{b}}'_{100})^8+1\bigr) \cdot\bigl(10000(q-1)\bigr)^{100}
\end{equation}
holds with $|J_{100}|=(3^{100}-1)/2$.
It implies
\[
\begin{aligned}
&\frac{\binom{|J_{100}|}{100}}{\bigl(400({\mathsf{b}}'_{100})^4+1\bigr)\bigl(400({\mathsf{b}}'_{100})^8+1\bigr)}\\
&\quad\geq\frac{(|J_{100}|-99)^{100}}{100^{100}\bigl(400({\mathsf{b}}'_{100})^4+1\bigr)\bigl(400({\mathsf{b}}'_{100})^8+1\bigr)}
>100^{100}(q-1)^{100}\geq100!(q-1)^{100}.
\end{aligned}
\]
Hence more than $100!(q-1)^{100}$ candidate members have the same pair of trace sums.

The Erd\H{o}s--Rado sunflower lemma \cite{ErdosRado1960} gives distinct $100$-element sets $\mathfrak F_1,\ldots,\mathfrak F_q$ among these candidates forming a sunflower.
Let $\mathfrak K$ be its core.
Every petal has size $100-|\mathfrak K|$, and distinctness gives $1\leq100-|\mathfrak K|\leq100$.
For $k=1,2$, subtracting the trace sum over $\mathfrak K$ from the common trace sum over each member shows that
\[
\sum_{\mathbf k\in\mathfrak F_i\setminus\mathfrak K}\operatorname{Tr}(({g'_{\mathbf k}})^{4k})
\]
is independent of $i$.

Partition the member positions as
\[
\{1,\ldots,q\}=I_1\sqcup I_2\sqcup I_3,\qquad
(|I_1|,|I_2|,|I_3|)=(8\,387\,848\,518,m_2,m_3),
\]
and define
\[
\begin{aligned}
\mathfrak Y:={}&\bigcup_{i\in I_1}\{({g'_{\mathbf k}})^4:\mathbf k\in\mathfrak F_i\setminus\mathfrak K\}\\
&\cup\bigcup_{i\in I_2}\{g_1^4({g'_{\mathbf k}})^{4{\mathsf{e}}}:\mathbf k\in\mathfrak F_i\setminus\mathfrak K,\ {\mathsf{e}}\in\{-1,1\}\}\\
&\cup\bigcup_{i\in I_3}\{g_3^4({g'_{\mathbf k}})^{4{\mathsf{e}}}:\mathbf k\in\mathfrak F_i\setminus\mathfrak K,\ {\mathsf{e}}\in\{-1,1\}\}.
\end{aligned}
\]
Put $\mathfrak Y^{-1}:=\{z^{-1}:z\in \mathfrak Y\}$.

For the three types of elements $z\in \mathfrak Y$, the Gaussian-prime valuations give
\[
\frac14\left(\operatorname{ord}_{(2+\mathrm i)}(z),\operatorname{ord}_{(4+\mathrm i)}(z), \bigl(\operatorname{ord}_{({p_{{\mathsf{p}}'_j}})}(z)\bigr)_{j=1}^{100}\right) =\begin{cases} (0,0,\mathbf k),&z=({g'_{\mathbf k}})^4,\\ (1,0,{\mathsf{e}}\mathbf k),&z=g_1^4({g'_{\mathbf k}})^{4{\mathsf{e}}},\\ (0,1,{\mathsf{e}}\mathbf k),&z=g_3^4({g'_{\mathbf k}})^{4{\mathsf{e}}}. \end{cases}
\]
The first two coordinates distinguish the multipliers; the remaining coordinates distinguish the indices and signs because the petals are pairwise disjoint and $J_{100}\cap(-J_{100})=\varnothing$.
Thus the points of $\mathfrak Y$ are distinct.
Inversion negates the tuple, replacing $\mathbf k$ by $-\mathbf k\notin J_{100}$ in the first type and producing a $-1$ in one of the first two entries in the other types.
Hence $\mathfrak Y\cap \mathfrak Y^{-1}=\varnothing$, and
\[
|\mathfrak Y|=(100-|\mathfrak K|)(8\,387\,848\,518+2m_2+2m_3)
=26\,011\,238\,400(100-|\mathfrak K|).
\]
For $k=1,2$, the trace product identity and the common petal trace sums give
\[
\sum_{z\in \mathfrak Y}\operatorname{Tr}(z^k) ={}\bigl(8\,387\,848\,518+m_2\operatorname{Tr}(g_1^{4k})+m_3\operatorname{Tr}(g_3^{4k})\bigr) \cdot\sum_{\mathbf k\in\mathfrak F_1\setminus\mathfrak K}\operatorname{Tr}(({g'_{\mathbf k}})^{4k})=0.
\]
Since $\mathfrak Y\subseteq \TQ^{\,4}$, the finite set
\[
\mathfrak X_{11}:=\rho^{-1}(\mathfrak Y\cup \mathfrak Y^{-1})
\]
is $R_4$-invariant, so the moments in \eqref{eq:Pk} satisfy $\psi_k(\mathfrak X_{11})=0$ for $1\leq k\leq11$ with $4\nmid k$.
For $k=1,2$, summation over the fibers gives
\[
\psi_{4k}(\mathfrak X_{11})=4\sum_{z\in \mathfrak Y}\operatorname{Tr}(z^k)=0.
\]
Thus \cref{lem:circle-moment-criterion} shows that $\mathfrak X_{11}$ is a spherical $11$-design.
Its size satisfies
\[
|\mathfrak X_{11}|=8|\mathfrak Y|
=8\cdot26\,011\,238\,400(100-|\mathfrak K|)
\leq20\,808\,990\,720\,000.
\]

\begingroup
\emergencystretch=2em
Enumerate $\binom{J_{100}}{100}$ and group its members $\mathfrak F$ according to the values of
\[
\bigl(\sum_{\mathbf k\in\mathfrak F}\operatorname{Tr}(({g'_{\mathbf k}})^4),\sum_{\mathbf k\in\mathfrak F}\operatorname{Tr}(({g'_{\mathbf k}})^8)\bigr).
\]
For each choice of $q$ distinct members $\mathfrak F_1,\ldots,\mathfrak F_q$ in one group, put $\mathfrak K:=\mathfrak F_1\cap\mathfrak F_2$ and test whether $\mathfrak F_i\cap\mathfrak F_j=\mathfrak K$ for all $1\leq i<j\leq q$.
\par\endgroup

Inequality \eqref{eq:circle-eleven-sunflower-inequality} and the sunflower lemma guarantee that this finite procedure succeeds.
\end{proof}

\subsection{Possible Sizes}\label{subsec:minimum-rational-circle}

The constructions in \cref{thm:minimum-circle-constructions} attain the divisibility bounds of \cref{thm:low-strength-divisibility}.
Together with \cref{prop:conic-cardinality-additivity}, they determine all possible sizes for strengths $4$ through $7$.
For $t\geq1$, write
\begin{equation}\label{def:circle-minimum-cardinality}
\begin{aligned}
\mathcal N_{1,\Q}^{R_2}(t)
&:=\{|\mathfrak X|:\mathfrak X\subseteq\TQ\text{ is a nonempty finite spherical }t\text{-design},\ -\mathfrak X=\mathfrak X\},\\
\mathcal N_{1,\Q}^{R_4}(t)
&:=\{|\mathfrak X|:\mathfrak X\subseteq\TQ\text{ is a nonempty finite spherical }t\text{-design},\ \mathrm i\mathfrak X=\mathfrak X\}.
\end{aligned}
\end{equation}
\Needspace{5\baselineskip}
\begin{theorem}\label{thm:circle-cardinality-spectra}\label{cor:circle-symmetric-spectra}
\[
\mathcal N_{1,\Q}(4)=2304\Z_{>0},
\qquad
\mathcal N_{1,\Q}(5)=4608\Z_{>0},
\qquad
\mathcal N_{1,\Q}(6)=\mathcal N_{1,\Q}(7)=9216\Z_{>0}.
\]
In particular,
\[
\min\mathcal N_{1,\Q}(4)=2304,
\qquad
\min\mathcal N_{1,\Q}(5)=4608,
\qquad
\min\mathcal N_{1,\Q}(6)=\min\mathcal N_{1,\Q}(7)=9216.
\]
Moreover,
\[
\mathcal N_{1,\Q}^{R_2}(4)=\mathcal N_{1,\Q}^{R_2}(5)=4608\Z_{>0},
\qquad
\mathcal N_{1,\Q}^{R_2}(6)=\mathcal N_{1,\Q}^{R_2}(7)=9216\Z_{>0},
\]
and
\[
\mathcal N_{1,\Q}^{R_4}(t)=9216\Z_{>0}
\qquad(4\leq t\leq7).
\]
\end{theorem}

\begin{proof}
\Cref{thm:low-strength-divisibility,cor:fourth-moment-symmetric-divisibility} show that, in each asserted equality between sets, the left-hand side is contained in the right-hand side.
The set $\mathfrak X_4$ in \cref{thm:minimum-circle-constructions} attains the unrestricted minimum at strength $4$.
The antipodal $5$-design $\mathfrak X_5$ attains the unrestricted minimum at strength $5$ and the antipodal minima at strengths $4$ and $5$.
The $9216$-point $R_4$-invariant design $\mathfrak X_7$ attains the unrestricted and antipodal minima at strengths $6,7$, and the $R_4$-invariant minima at strengths $4,5,6,7$, since a $7$-design is a design of every smaller strength and $R_4$-invariance implies antipodality.
For each $\mathfrak X\in\{\mathfrak X_4,\mathfrak X_5,\mathfrak X_7\}$, pairwise disjoint sets of the form ${v}\mathfrak X$, with ${v}\in\TQ$, realize every positive multiple of $|\mathfrak X|$, as in \cref{prop:conic-cardinality-additivity}.
Multiplication in $\TQ$ commutes with $-1$ and $\mathrm i$, so ${v}\mathfrak X$ is antipodal whenever $\mathfrak X$ is antipodal and is $R_4$-invariant whenever $\mathfrak X$ is $R_4$-invariant.
\end{proof}

The unit circle and the ellipses in \cref{subsec:norm-one-local-divisibility} lie in $\R^2$ and are special cases of the rational ellipsoids considered next.
We apply \cref{thm:integer-to-equal-weight} to show that such an ellipsoid admits rational designs of every strength precisely when it has a rational point.

\section{Rational Ellipsoids and Open Problems}\label{sec:problems}

\subsection{Rational Designs on Positive-Definite Ellipsoids}\label{subsec:imaginary-quadratic-ellipsoids}

Following the study of ellipsoidal designs initiated by Pandey \cite{Pandey2022} and developed further by others \cite{MatsumuraSawa2025,InagakiMatsumuraSawaUchida2025}, we investigate the existence of rational designs on ellipsoids.
Our main result characterizes this existence in terms of rational points on the underlying ellipsoid.
The proof combines the density of rational orthogonal transformations with the construction developed in \cref{sec:rational-design-existence}.

Fix an integer \({\mathsf{m}}\geq2\) and a positive-definite symmetric matrix \(M\in \operatorname{Mat}_{\mathsf{m}}(\Q)\).
For \(\mathbb K\in\{\Q,\R\}\), put
\[
\mathrm{O}(M,\mathbb K):=\{g\in \mathrm{GL}_{\mathsf{m}}(\mathbb K):g^{\mathsf T}Mg=M\},
\qquad
\mathrm{SO}(M,\mathbb K):=\{g\in \mathrm{O}(M,\mathbb K):\det g=1\}.
\]
The classical Cayley representation gives the following density statement, which supplies the rational-point density needed in \cref{thm:rational-ellipsoidal-existence}.

\begin{proposition}\label{prop:cayley-density}
The groups \(\mathrm{SO}(M,\Q)\) and \(\mathrm{O}(M,\Q)\) are dense in \(\mathrm{SO}(M,\R)\) and \(\mathrm{O}(M,\R)\), respectively.
\end{proposition}

\begin{proof}
Put
\[
\mathfrak o(M):=\{M_3\in \operatorname{Mat}_{\mathsf{m}}(\R):M_3^{\mathsf T}M+MM_3=0\}.
\]
Since $M$ is symmetric, $M_3\in\mathfrak o(M)$ if and only if $MM_3$ is skew-symmetric.
Hence
\[
\mathfrak o(M)\cap \operatorname{Mat}_{\mathsf{m}}(\Q)
=\{M^{-1}M_4:M_4\in \operatorname{Mat}_{\mathsf{m}}(\Q),\ (M_4)^{\mathsf T}=-M_4\}
\]
is dense in $\mathfrak o(M)$.
For $M_3\in\mathfrak o(M)$ and $\mathbf x\in\R^{\mathsf{m}}$, the equality $(\operatorname{Id}_{\mathsf{m}}+M_3)\mathbf x=\mathbf0$ implies
\[
0=\mathbf x^{\mathsf T}(M_3^{\mathsf T}M+MM_3)\mathbf x=-2\mathbf x^{\mathsf T}M\mathbf x,
\]
so $\mathbf x=\mathbf0$.
Thus $\operatorname{Id}_{\mathsf{m}}+M_3$ is invertible.
Conjugation by the positive-definite square root $M^{1/2}$ takes $\mathfrak o(M)$ to the space of real skew-symmetric matrices.
The Cayley transform
\[
\operatorname{cay}(M_3):=(\operatorname{Id}_{\mathsf{m}}-M_3)(\operatorname{Id}_{\mathsf{m}}+M_3)^{-1}
\]
maps $\mathfrak o(M)$ bijectively onto
\[
\{g\in \mathrm{SO}(M,\R):\det(\operatorname{Id}_{\mathsf{m}}+g)\ne0\},
\]
with inverse
\[
g\mapsto(\operatorname{Id}_{\mathsf{m}}-g)(\operatorname{Id}_{\mathsf{m}}+g)^{-1};
\]
see \cite{Gallier2006}.
Both maps are continuous, and $\operatorname{cay}$ sends rational matrices to rational matrices.
The map $g\mapsto M^{1/2}gM^{-1/2}$ identifies $\mathrm{SO}(M,\R)$ with the connected group $\mathrm{SO}({\mathsf{m}})$ \cite[Chapter~1, pp.~3--30]{Hall2015}.
The function $g\mapsto\det(\operatorname{Id}_{\mathsf{m}}+g)$ is real analytic and nonzero at $\operatorname{Id}_{\mathsf{m}}$, so \cref{lem:topological-inputs}(ii) shows that $\operatorname{cay}(\mathfrak o(M))$ is dense in $\mathrm{SO}(M,\R)$.
Since rational matrices are dense in $\mathfrak o(M)$ and the Cayley map is continuous and preserves rationality, $\mathrm{SO}(M,\Q)$ is dense in this image, and hence in $\mathrm{SO}(M,\R)$.
Finally, for $\mathbf0\ne\mathbf x\in\Q^{\mathsf{m}}$, the matrix
\[
\operatorname{Id}_{\mathsf{m}}-\frac{2\mathbf x\mathbf x^{\mathsf T}M}{\mathbf x^{\mathsf T}M\mathbf x}
\]
belongs to $\mathrm{O}(M,\Q)$ and has determinant $-1$.
Left multiplication by this rational matrix interchanges the two determinant components, so density in $\mathrm{SO}(M,\R)$ gives density in $\mathrm{O}(M,\R)$.
\end{proof}

For \(r\in\Q_{>0}\), put
\[
\mathcal E_M(r):=\{\mathbf x\in\R^{\mathsf{m}}:\mathbf x^{\mathsf T}M\mathbf x=r\}.
\]
As above, $\mathcal E_M(r)(\Q)=\mathcal E_M(r)\cap\Q^{\mathsf{m}}$.
The map $\mathbf x\mapsto r^{-1/2}M^{1/2}\mathbf x$ is a bijection from $\mathcal E_M(r)$ onto $\mathbb S^{{\mathsf{m}}-1}$ that intertwines the $\mathrm{O}(M,\R)$-action with the standard $\mathrm{O}({\mathsf{m}})$-action.
Thus $\mathcal E_M(r)$ is a compact, connected real-analytic manifold, and $\mathrm{O}(M,\R)$ is compact and acts transitively on it.
Let \(\mu_{\mathrm{O}(M,\R)}\) denote normalized Haar measure. For any \(\mathbf x_0\in\mathcal E_M(r)\), define the invariant probability measure \(\mu_{M,r}\) by
\[
 \int_{\mathcal E_M(r)}\varphi\,\d\mu_{M,r}
 =\int_{\mathrm{O}(M,\R)}\varphi(g\mathbf x_0)\,\d\mu_{\mathrm{O}(M,\R)}(g)
 .
\]
Here \(\varphi\) is any continuous real-valued function on \(\mathcal E_M(r)\). The action is transitive, so this measure is independent of \(\mathbf x_0\).
Using a rational change of coordinates that makes the quadratic form diagonal, we apply \cref{thm:integer-to-equal-weight} in the form of \cref{cor:norm-one-equal-weight} to construct rational designs on $\mathcal E_M(r)$ whenever $\mathcal E_M(r)(\Q)\ne\varnothing$.

\begin{theorem}\label{thm:rational-ellipsoidal-existence}
The following are equivalent.
\begin{enumerate}
\item \(\mathcal E_M(r)(\Q)\neq\varnothing\).
\item For every integer \(t\geq1\) and every finite set \(E\subseteq \mathcal E_M(r)(\Q)\), there exists a rational \(\mathcal P_t(\mathcal E_M(r);\R)\)-design on \((\mathcal E_M(r),\mu_{M,r})\) disjoint from \(E\).
\end{enumerate}
\end{theorem}

\begin{proof}
For \textup{(2)} $\Rightarrow$ \textup{(1)}, take $t=1$ and $E=\varnothing$; the resulting design is a nonempty subset of $\mathcal E_M(r)(\Q)$.
Assume \textup{(1)}. Fix an integer $t\geq1$ and a finite $E\subseteq \mathcal E_M(r)(\Q)$, and choose $\mathbf x_0\in \mathcal E_M(r)(\Q)$.
We first make a rational change of coordinates so that a norm-one group associated with an imaginary quadratic field acts on the first two coordinates.
Apply Gram--Schmidt to the standard basis for the inner product $\mathbf x^{\mathsf T}M\mathbf x^{\prime}$, without normalizing the vectors.
At each step, the projection coefficients are ratios of rational inner products to nonzero rational squared lengths.
Inductively, all resulting vectors have rational coordinates.
Let \(M_5\in \mathrm{GL}_{\mathsf{m}}(\Q)\) have these rational orthogonal basis vectors as columns, and put
\[
 M_1:=M_5^{\mathsf T}MM_5.
\]
Then \(M_1\) is diagonal with positive rational diagonal entries.
Choose a positive squarefree integer \(\ell\) such that
\((M_1)_{22}/(\ell(M_1)_{11})\) is a square in \(\Q_{>0}\), and put
\[
 M'_5:=M_5\operatorname{diag}\!\left(1,
 \sqrt{\frac{\ell(M_1)_{11}}{(M_1)_{22}}},1,\ldots,1\right),
 \qquad M':={M'_5}^{\mathsf T}MM'_5.
\]
Thus \(M'_5\in \mathrm{GL}_{\mathsf{m}}(\Q)\), and \(M'\) is diagonal with
\((M')_{22}=\ell(M')_{11}\) and \((M')_{jj}\in\Q_{>0}\) for \(1\leq j\leq {\mathsf{m}}\).
The linear map $\mathbf x^{\prime}\mapsto M'_5\mathbf x^{\prime}$ is a bijection from $\mathcal E_{M'}(r)$ to $\mathcal E_M(r)$ taking $\mathcal E_{M'}(r)(\Q)$ onto $\mathcal E_M(r)(\Q)$.
Conjugation $g\mapsto M'_5g{M'_5}^{-1}$ identifies $\mathrm{O}(M',\R)$ with $\mathrm{O}(M,\R)$ and preserves normalized Haar measure.
Since $M'_5(g{M'_5}^{-1}\mathbf x_0)=(M'_5g{M'_5}^{-1})\mathbf x_0$, it follows that
\[
\int_{\mathcal E_{M'}(r)}\varphi(M'_5\mathbf x^{\prime})\,\d\mu_{M',r}(\mathbf x^{\prime})
=\int_{\mathcal E_M(r)}\varphi(\mathbf x)\,\d\mu_{M,r}(\mathbf x)
\]
for every continuous real-valued function $\varphi$ on $\mathcal E_M(r)$.
Since composition with $M'_5$ preserves polynomial degree, it suffices to construct a rational $\mathcal P_t(\mathcal E_{M'}(r);\R)$-design disjoint from ${M'_5}^{-1}E$.

Put
\[
Z_0:=\{g({M'_5}^{-1}\mathbf x_0):g\in \mathrm{O}(M',\Q)\}.
\]
This set is contained in $\mathcal E_{M'}(r)(\Q)$ and is invariant under $\mathrm{O}(M',\Q)$.
By \cref{prop:cayley-density}, the set $Z_0$ is dense in $\mathcal E_{M'}(r)$.
Apply \cref{prop:rational-invariant-average} with this $Z_0$,
\[
\Gamma=\mathrm{O}(M',\R),\qquad \Gamma_0=\mathrm{O}(M',\Q),
\]
and
\[
\mathcal V_{\Q}:=\{f|_{\mathcal E_{M'}(r)}:f\in\Q[X_1,\ldots,X_{\mathsf{m}}],\ \deg f\leq t\}.
\]
This finite-dimensional $\Q$-vector space contains $\mathbf1$ and is invariant under $\mathrm{O}(M',\Q)$.
Moreover, every function in $\mathcal V_{\Q}$ takes rational values on $Z_0$.
Applying \cref{prop:rational-invariant-average} in every degree gives rational integrals for all polynomials with rational coefficients.
Also, \cref{lem:topological-inputs}(i) gives positive measure to every nonempty relatively open subset of $\mathcal E_{M'}(r)$.

Put \(\mathbb K_\ell:=\Q(\sqrt{-\ell})\). We use multiplication in $\mathbb K_\ell$ to separate the points of an integer-weighted design.
Define $\eta_{M'}:\mathcal E_{M'}(r)\to\C$ by $\eta_{M'}(\mathbf x^{\prime}):=x_1^{\prime}+x_2^{\prime}\sqrt{-\ell}$.
For ${v}=x_1+x_2\sqrt{-\ell}\in T_{\mathbb K_\ell}$, define
\[
{v}\cdot \mathbf x^{\prime}:=(x_1x_1^{\prime}-\ell x_2x_2^{\prime},\ x_2x_1^{\prime}+x_1x_2^{\prime},\ x_3^{\prime},\ldots,x_{\mathsf{m}}^{\prime}).
\]
For every ${v}\in T_{\mathbb K_\ell}$, the matrix of $\mathbf x^{\prime}\mapsto {v}\cdot \mathbf x^{\prime}$ belongs to $\mathrm{O}(M',\Q)$; hence the action preserves $\mathcal E_{M'}(r)$, $\mathcal E_{M'}(r)(\Q)$, and $\mu_{M',r}$.
This action satisfies $\eta_{M'}({v}\cdot \mathbf x^{\prime})={v}\eta_{M'}(\mathbf x^{\prime})$.
The change of variables
\[
(x_1^{\prime},x_2^{\prime})\mapsto(x_1^{\prime}+x_2^{\prime}\sqrt{-\ell},x_1^{\prime}-x_2^{\prime}\sqrt{-\ell})
\]
has determinant \(-2\sqrt{-\ell}\ne0\), and is therefore invertible over \(\C\). Hence, the restrictions to \(\mathcal E_{M'}(r)\) of the real and imaginary parts of
\[
\mathbf x^{\prime}\longmapsto\eta_{M'}(\mathbf x^{\prime})^{k_1}\overline{\eta_{M'}(\mathbf x^{\prime})}^{\,k_2}\prod_{j=3}^{{\mathsf{m}}}(x_j^{\prime})^{k_j},
\qquad
k_1,\ldots,k_{\mathsf{m}}\in\Z_{\geq0},
\qquad
\sum_{j=1}^{\mathsf{m}} k_j\leq t,
\]
span \(\mathcal P_t(\mathcal E_{M'}(r);\R)\).

Since $\overline {v}={v}^{-1}$ for ${v}\in T_{\mathbb K_\ell}$, each displayed monomial transforms by the character \({v}\mapsto {v}^{k_1-k_2}\).
The function $\eta_{M'}$ is real analytic and is nonzero at
\[
(\sqrt{r/(M')_{11}},0,\ldots,0)\in \mathcal E_{M'}(r).
\]
Thus $\eta_{M'}$ is not identically zero on the connected manifold $\mathcal E_{M'}(r)$.
By \cref{lem:topological-inputs}(ii), the set
\[
W:=\{\mathbf x^{\prime}\in \mathcal E_{M'}(r):\eta_{M'}(\mathbf x^{\prime})\ne0\}
\]
is open and dense.
Since $\eta_{M'}({v}\cdot \mathbf x^{\prime})={v}\eta_{M'}(\mathbf x^{\prime})$ and ${v}\ne0$, the set $W$ is $T_{\mathbb K_\ell}$-invariant.
Because $W$ is open and dense and $\mathcal E_{M'}(r)(\Q)$ is dense, $W(\Q)=W\cap \mathcal E_{M'}(r)(\Q)$ is dense in $\mathcal E_{M'}(r)$.
On $W(\Q)$, the map $\eta_{M'}$ takes values in $\mathbb K_\ell^\times$.
Apply \cref{prop:cxx-integer-weighted} using the dense subset $W(\Q)$ to obtain an integer-weighted rational $\mathcal P_t(\mathcal E_{M'}(r);\R)$-design supported on $W(\Q)$.
Apply \cref{cor:norm-one-equal-weight} with $\mathbb K=\mathbb K_\ell$, $\eta=\eta_{M'}|_{W(\Q)}$, and $k_0=1$ to obtain a rational $\mathcal P_t(\mathcal E_{M'}(r);\R)$-design $\mathfrak X$ avoiding ${M'_5}^{-1}E$.
By the integral identity above and the fact that $\varphi\circ M'_5\in\mathcal P_t(\mathcal E_{M'}(r);\R)$ for every $\varphi\in\mathcal P_t(\mathcal E_M(r);\R)$, the set
\[
M'_5\mathfrak X:=\{M'_5\mathbf x^{\prime}:\mathbf x^{\prime}\in\mathfrak X\}
\]
is a rational $\mathcal P_t(\mathcal E_M(r);\R)$-design.
Moreover,
\[
M'_5\mathfrak X\cap E=M'_5(\mathfrak X\cap {M'_5}^{-1}E)=\varnothing.
\]
\end{proof}

\begin{remark}\label{rem:Qn-ellipsoids}
For every integer ${\mathsf{m}}\geq2$, the ellipsoid
\[
\left\{\mathbf x\in\R^{\mathsf{m}}:\sum_{j=1}^{\mathsf{m}}\frac{j(j+1)}2x_j^2=1\right\}
\]
contains the rational point $(1,0,\ldots,0)$, so \cref{thm:rational-ellipsoidal-existence} gives rational designs of every strength on it.
These are the ellipsoids of Inagaki et al.~\cite[Problem~3.2]{InagakiMatsumuraSawaUchida2025}, whose dimension parameter is \({\mathsf{m}}+1\) in their notation. They ask for tight rational designs, with tightness defined in \cite[Definition~3.1]{InagakiMatsumuraSawaUchida2025}.
\end{remark}

\subsection{Open Problems}\label{subsec:rational-design-cardinality-problems}\label{subsec:pte-minimum-existence-problems}

For every integer $t\geq8$, \cref{prop:conic-cardinality-additivity} gives \(\bigl(\min\mathcal N_{1,\Q}(t)\bigr)\Z_{>0}\subseteq\mathcal N_{1,\Q}(t)\).
This inclusion is an equality if and only if
\[
\gcd\bigl(\mathcal N_{1,\Q}(t)\bigr)
=\min\mathcal N_{1,\Q}(t).
\]

\begin{problem}\label{prob:circle-higher-strength-cardinalities}
For every integer \(t\geq8\), \cref{thm:low-strength-divisibility} gives
\[
{\delta_{\TQ}}(t)\mid \gcd\bigl(\mathcal N_{1,\Q}(t)\bigr)\mid \min\mathcal N_{1,\Q}(t).
\]
\begin{enumerate}
\item Determine \(\min\mathcal N_{1,\Q}(8)\).
\item Determine the integers \(t\geq8\) for which \(\gcd\bigl(\mathcal N_{1,\Q}(t)\bigr)=\min\mathcal N_{1,\Q}(t)\).
\item Determine the additive submonoid \(\mathcal N_{1,\Q}(t)\cup\{0\}\) for every \(t\geq8\).
\end{enumerate}
\end{problem}

\begin{center}
* * *
\end{center}

Rational designs for the Chebyshev measure $\pi^{-1}(1-u^2)^{-1/2}du$ on $(-1,1)$ have been extensively studied by Mishima et al~\cite{MishimaLuSawaUchida2026}.
The Chebyshev measure is the $d=2$ case of the Gegenbauer measures in \cref{rem:cxx-alternative-proof}\textup{(ii)}. We study the possible sizes of rational $7$-designs for this measure that are invariant under $u\mapsto-u$.
For a nonnegative integer $t$, a nonempty finite set $\mathfrak{X} \subseteq (-1,1)\cap\mathbb{Q}$ is called a rational $t$-design for the Chebyshev measure
if
\[
\frac{1}{|\mathfrak X|} \sum_{u \in \mathfrak{X}} \varphi(u) = \frac{1}{\pi} \int_{-1}^{1}\frac{\varphi(u)}{\sqrt{1-u^2}}\,\mathrm du \qquad \bigl( \varphi \in \mathbb{R}[X],\ \deg \varphi \le t \bigr).
\]
It is
called
antipodal if $-\mathfrak{X}=\mathfrak{X}$.

Problem~8.5 of Mishima et al.~\cite{MishimaLuSawaUchida2026} asks whether an antipodal rational $7$-design for the Chebyshev measure exists with an even number of points, and which even sizes can occur.
The following result characterizes the possible sizes under the additional condition that $\sqrt{1-u^2}\in\mathbb Q$ for every node. It thus provides explicit examples related to Problem 8.5 of Mishima et al.~\cite{MishimaLuSawaUchida2026}.

\begin{proposition}\label{prop:chebyshev-seven-designs}
For a positive integer $n$, there exists an $n$-point antipodal rational $7$-design $\mathfrak{X}$ for the Chebyshev measure satisfying
\[
\sqrt{1-u^2} \in \mathbb{Q} \qquad (u\in\mathfrak{X})
\]
if and only if $4608\mid n$.
\end{proposition}

\begin{proof}
The unrotated set $\mathfrak{X}_7$ is a $9216$-point rational spherical $7$-design invariant under complex conjugation and the action of $R_4$, and is disjoint from $R_4$.
Each nonempty fiber of real-part projection on $\mathfrak{X}_7$ consists of two points $z,\overline z$.
Hence its real-part image is a $4608$-point antipodal rational $7$-design for the Chebyshev measure, and for $u=\Re z$ we have $\sqrt{1-u^2}=|\Im z|\in\mathbb{Q}$.

Conversely, let $\mathfrak{X}$ be a rational $7$-design for the Chebyshev measure satisfying the square-root condition, and put
\[
\mathfrak{X}' := \{ u \pm \mathrm i \sqrt{1-u^2} : u \in \mathfrak{X} \} \subseteq \mathbb{T}_{\mathbb{Q}}.
\]
This set has $2|\mathfrak{X}|$ points.
The Chebyshev polynomials of the first kind, defined by $T_j(\cos\theta)=\cos(j\theta)$, have zero average with respect to the Chebyshev
measure for $j\ge1$.
The design conditions therefore give
\[
\psi_j(\mathfrak{X}') = 2 \sum_{u \in \mathfrak{X}} T_j(u) = 0 \qquad (1\le j\le7).
\]
By \cref{lem:circle-moment-criterion}, $\mathfrak{X}'$ is a rational spherical $7$-design.
\Cref{thm:circle-cardinality-spectra} gives $9216 \mid 2|\mathfrak{X}|$, and hence $4608 \mid |\mathfrak{X}|$.

To realize every positive multiple, fix a positive integer $k$.
Start with $\mathfrak{X}_7$ if $k$ is odd and with the empty set if $k$ is even, and add pairs
\[
v \mathfrak{X}_7 \sqcup \overline v\,\mathfrak{X}_7, \qquad v \in \mathbb{T}_{\mathbb{Q}},
\]
until $k$ rotated copies have been chosen.
Only finitely many values of $v$ cause an intersection with previously chosen copies or with $R_4$, and an intersection within a pair requires $v^2=z'/z$ for some $z,z'\in\mathfrak{X}_7$, again excluding only finitely many values.
Since $\mathbb{T}_{\mathbb{Q}}$ is infinite, all these values can be avoided.
The resulting union is a $9216k$-point rational spherical $7$-design invariant under complex conjugation and the action of $R_4$, and is disjoint from $R_4$.
Its real-part image is a $4608k$-point antipodal rational $7$-design for the Chebyshev measure satisfying the square-root condition.
\end{proof}

\begin{problem}\label{prob:antipodal-chebyshev-seven-cardinalities}
Determine the spectrum of antipodal rational $7$-designs for the Chebyshev measure, that is, the set
\[
\left\{|\mathfrak{X}|:
\begin{array}{l}
\varnothing \ne \mathfrak{X} \subseteq (-1,1)\cap\mathbb{Q},\quad-\mathfrak{X}=\mathfrak{X},\\
\mathfrak{X} \text{ is a finite rational }7\text{-design for the Chebyshev measure}
\end{array}
\right\}.
\]
\end{problem}

\begin{center}
* * *
\end{center}

For integers $d\geq2$ and $t\geq1$, let $\mathcal N_{\mathbb S^{d-1},\Q}(t)$ be the set of sizes of rational spherical $t$-designs on $\mathbb S^{d-1}$.
This set is nonempty by \cref{thm:rational-spherical-existence}.

\begin{proposition}\label{prop:minimum-rational-sphere-two-five}
The minimum size of a rational spherical $5$-design on $\mathbb S^2$ is $192$.
\end{proposition}

\begin{proof}
Let $G:=\mathrm{SO}(3;\Z)$ be the group of $3\times3$ signed permutation matrices of determinant $1$, and define
\[
V:=\left\{
\begin{array}{l}
(0,63,84),\ (20,40,95),\ (30,45,90),\ (40,55,80),\\
(1,32,100),\ (4,40,97),\ (14,70,77),\ (-14,70,77)
\end{array}\right\}.
\]
Then
\[
\mathfrak X:=\bigcup_{\mathbf v\in V}G\!\left(\frac{\mathbf v}{105}\right)
\]
is a rational spherical $5$-design consisting of $192$ distinct points.
Each vector in $V$ has norm $105$, and the eight orbits $G(\mathbf v/105)$ are pairwise disjoint and have $24$ points each.
Averaging monomials of degree at most $5$ over $G$, and using $x_1^2+x_2^2+x_3^2=1$, reduces the spherical $5$-design conditions to the requirement that the average of $x_1^4+x_2^4+x_3^4$ equals $3/5$.
Indeed,
\[
\mathcal P_5(\mathbb S^2;\R)^G
:=\{f\in\mathcal P_5(\mathbb S^2;\R)\mid f\circ g=f\ (\forall g\in G)\}
=\operatorname{span}_{\R}\{1,x_1^4+x_2^4+x_3^4\}.
\]
Direct calculation gives
\[
\frac{1}{|\mathfrak X|}\sum_{\mathbf x\in\mathfrak X}(x_1^4+x_2^4+x_3^4)
=\frac{1}{8\cdot105^4}\sum_{\mathbf v\in V}(v_1^4+v_2^4+v_3^4)
=\frac35.
\]
We now prove the lower bound.

For $\alpha,\beta\in\Z_{(2)}$ and a positive integer $\mathsf m$, the congruence $\alpha\equiv\beta\pmod{\mathsf m\Z_{(2)}}$ means $\alpha-\beta\in\mathsf m\Z_{(2)}$.
For any $\mathbf x=(x_1,x_2,x_3)\in\mathbb S^2\cap\Q^3$, choose integers $a_1,a_2,a_3$ and a positive integer $D$ such that
\[
\mathbf x=\frac1D(a_1,a_2,a_3),\qquad \gcd(a_1,a_2,a_3,D)=1.
\]
The sphere equation is $a_1^2+a_2^2+a_3^2=D^2$.
Since a square is congruent to $0$ or $1$ modulo $4$, if $D$ were even, all three squares would be congruent to $0$, making $a_1,a_2,a_3,D$ all even.
This contradicts the gcd condition.
Thus $D$ is odd, and since the sum of three squares is congruent to $1$ modulo $4$, exactly one of $a_1,a_2,a_3$ is odd.
Consequently, all coordinates belong to $\Z_{(2)}$, and exactly one lies outside $2\Z_{(2)}$.

Put $f_0(u):=(u^2-1)(u^2-9)$.
Writing $r\in\Z_{(2)}\setminus2\Z_{(2)}$ as $r=2s+1$ with $s\in\Z_{(2)}$, we have
\[
f_0(r)=16(s-1)s(s+1)(s+2).
\]
Two of $s-1,s,s+1,s+2$ belong to $2\Z_{(2)}$, and one of these belongs to $4\Z_{(2)}$, so their product belongs to $8\Z_{(2)}$.
Since $r+1\in2\Z_{(2)}$, we obtain
\[
f_0(r)\in128\Z_{(2)},\qquad (r+1)f_0(r)\in256\Z_{(2)}.
\]
On the other hand, for $r=4s\in4\Z_{(2)}$, we have
\[
(8r+1)f_0(r)\equiv9+32s(s+1)\equiv9\pmod{64\Z_{(2)}}.
\]

Let $\mathfrak Y$ be a rational spherical $5$-design with $N$ points, and put
\[
f_1(\mathbf x):=x_1+x_2+x_3,\qquad f_2(\mathbf x):=x_1+x_2.
\]
Their moments with respect to normalized spherical measure $\mu_{\mathbb S^2}$ are
\[
\int_{\mathbb S^2}f_1^2\,\d\mu_{\mathbb S^2}=1,\quad
\int_{\mathbb S^2}f_1^4\,\d\mu_{\mathbb S^2}=\frac95,\quad
\int_{\mathbb S^2}f_2^2\,\d\mu_{\mathbb S^2}=\frac23,\quad
\int_{\mathbb S^2}f_2^4\,\d\mu_{\mathbb S^2}=\frac45,
\]
and every homogeneous polynomial of odd degree has zero average.
At every point, $f_1(\mathbf x)\in\Z_{(2)}\setminus2\Z_{(2)}$, so $(f_1(\mathbf x)+1)f_0(f_1(\mathbf x))\in256\Z_{(2)}$.
Since this polynomial has degree $5$, the design condition gives
\[
\sum_{\mathbf x\in\mathfrak Y}\bigl(f_1(\mathbf x)+1\bigr)f_0\bigl(f_1(\mathbf x)\bigr)
=N\left(\frac95-10+9\right)=\frac45N\in256\Z_{(2)}.
\]
Hence $64\mid N$.

Put
\[
n_i:=\bigl|\{\mathbf x\in\mathfrak Y:x_i\notin2\Z_{(2)}\}\bigr|\qquad(i=1,2,3).
\]
If $x_3\in2\Z_{(2)}$, then $f_2(\mathbf x)\in\Z_{(2)}\setminus2\Z_{(2)}$.
If $x_3\notin2\Z_{(2)}$, then $x_3^2\equiv1\pmod{8\Z_{(2)}}$; writing $x_1=2r_1$ and $x_2=2r_2$ with $r_1,r_2\in\Z_{(2)}$, the sphere equation gives
\[
r_1+r_2\equiv r_1^2+r_2^2\equiv0\pmod{2\Z_{(2)}}.
\]
Consequently, $f_2(\mathbf x)=2(r_1+r_2)\in4\Z_{(2)}$.
It follows that
\[
\bigl(8f_2(\mathbf x)+1\bigr)f_0\bigl(f_2(\mathbf x)\bigr)
\equiv
\begin{cases}
9,&x_3\notin2\Z_{(2)},\\
0,&x_3\in2\Z_{(2)}
\end{cases}
\pmod{64\Z_{(2)}}.
\]
Summing this congruence and using the degree-$5$ design condition and $64\mid N$, we obtain
\[
\begin{aligned}
9n_3&\equiv\sum_{\mathbf x\in\mathfrak Y}\bigl(8f_2(\mathbf x)+1\bigr)f_0\bigl(f_2(\mathbf x)\bigr)\\
&=N\left(\frac45-\frac{20}3+9\right)=\frac{47}{15}N\equiv0\pmod{64\Z_{(2)}}.
\end{aligned}
\]
Thus $64\mid n_3$, and coordinate permutations show that each $n_i$ is a multiple of $64$.
If all $n_i$ are positive, then $N=n_1+n_2+n_3\geq192$.
If one is zero, we may permute the coordinates so that $n_3=0$.
Then $f_2(\mathbf x)\in\Z_{(2)}\setminus2\Z_{(2)}$ at every point, and similarly the design condition gives
\[
\frac{47}{15}N=\sum_{\mathbf x\in\mathfrak Y}\bigl(f_2(\mathbf x)+1\bigr)f_0\bigl(f_2(\mathbf x)\bigr)\in256\Z_{(2)}.
\]
Hence $256\mid N$, and therefore $N\geq256$.
Thus $N\geq192$ in every case.
\end{proof}

In fixed dimension, the minimum size of real spherical designs grows on the order of $t^{d-1}$, by the lower bound of Delsarte--Goethals--Seidel~\cite{DelsarteGoethalsSeidel1977} and the existence theorem of Bondarenko--Radchenko--Viazovska~\cite[Theorem~1]{BondarenkoRadchenkoViazovska2013}.
For rational designs on the circle, however, \cref{thm:conic-exponential-lower-bound} gives an exponential lower bound.
We ask for a classification of the minimum configurations and for the effect of rationality on the minimum size in higher dimensions.

\begin{problem}\label{prob:higher-dimensional-rational-cardinality}
\begin{enumerate}
\item[\rm(i)] Classify the $192$-point rational spherical $5$-designs on $\mathbb S^2$ up to orthogonal transformations.
For $t=6,7$, determine $\min\mathcal N_{\mathbb S^2,\Q}(t)$ and $\min\mathcal N_{\mathbb S^3,\Q}(t)$, and construct configurations attaining each minimum.
\item[\rm(ii)] For a fixed integer $d\geq3$, determine the growth rate of $\min\mathcal N_{\mathbb S^{d-1},\Q}(t)$ as $t\to\infty$.
In particular, determine the dimensions for which
\[
\min\mathcal N_{\mathbb S^{d-1},\Q}(t)=O_d(t^{d-1}).
\]
\item[\rm(iii)] Give constructions with explicit upper bounds on the number of points and the denominators of their coordinates.
In particular, in the dimensions for which the upper bound in \textup{(ii)} holds, give an explicit constant $C_d>0$ and an algorithm that, on input strength $t$, outputs exact rational coordinates of a configuration with at most $C_dt^{d-1}$ points.
\end{enumerate}
\end{problem}

\begin{center}
* * *
\end{center}

For integers \(t\geq1\) and \(q\geq2\), let \(\mathcal T_t(\TQ)\) be the character family defined in \eqref{eq:circle-character-family}.
These characters take values in $\Q(\mathrm i)^\times$.
Since $\TQ$ is infinite by \cref{lem:norm-one-condition-D}, we can choose $q$ rotations of a minimum rational spherical $t$-design successively, avoiding the finitely many ratios of points that would cause an intersection with a previously chosen block.
By \cref{lem:circle-moment-criterion}, the sum of the $k$-th powers over each rotated block is zero for $1\leq k\leq t$, so these blocks of equal size form a $(\TQ,\mathcal T_t(\TQ))$-PTE set solution.
Consequently,
\[
\kappa^{\mathrm{set}}_{\TQ,\mathcal T_t(\TQ)}(q)\leq \min\mathcal N_{1,\Q}(t).
\]
By \cref{prop:PTEbound-newton}, $\kappa^{\mathrm{set}}_{\TQ,\mathcal T_t(\TQ)}(q)\geq t+1$, and \cref{cor:imaginary-quadratic-minimum} excludes equality unless $t\in\{1,3\}$.
In particular, $\kappa^{\mathrm{set}}_{\TQ,\mathcal T_2(\TQ)}(q)\geq4$.
With $R_k$ as in \eqref{eq:roots-of-unity}, distinct $R_2$-cosets give solutions of size $2$ for $t=1$, and distinct $R_4$-cosets give solutions of size $4$ for $t=2,3$.
Thus, for every $q\geq2$,
\[
\begin{gathered}
\kappa^{\mathrm{set}}_{\TQ,\mathcal T_1(\TQ)}(q)=2,
\qquad
\kappa^{\mathrm{set}}_{\TQ,\mathcal T_2(\TQ)}(q)=\kappa^{\mathrm{set}}_{\TQ,\mathcal T_3(\TQ)}(q)=4,\\
t+2\leq \kappa^{\mathrm{set}}_{\TQ,\mathcal T_t(\TQ)}(q)\leq \min\mathcal N_{1,\Q}(t)\quad(t\geq4).
\end{gathered}
\]

\begin{problem}\label{prob:rational-circle-pte-minimum}
Determine \(\kappa^{\mathrm{set}}_{\TQ,\mathcal T_4(\TQ)}(2)\), and more generally \(\kappa^{\mathrm{set}}_{\TQ,\mathcal T_t(\TQ)}(q)\) for all \(t\geq4\) and \(q\geq2\).
Classify the solutions of size $\kappa^{\mathrm{set}}_{\TQ,\mathcal T_t(\TQ)}(q)$ up to block permutations and the common transformations
\[
z\longmapsto {v} z
\quad\text{or}\quad
z\longmapsto {v}\overline z
\qquad({v}\in\TQ).
\]
\end{problem}

\begin{center}
* * *
\end{center}

The counting argument used to prove Thuswaldner--Ziegler's conjecture (\cref{prop:distinct-unit-sums,cor:thuswaldner-ziegler-conjecture}) also yields the following existence result for free abelian groups of rank at least two.

\begin{proposition}\label{prop:free-abelian-rank-two-existence}
Let $G=(\Z^{\mathsf d},+)$ for an integer $\mathsf d\geq2$, and let $\mathbb K$ be a number field.
For every finite set $\mathcal S\subseteq\Hom(G,\mathbb K^\times)$ satisfying $\mathbf1\in\mathcal S$ and $\bigcap_{\chi\in\mathcal S}\ker\chi=\{0\}$, we have
\[
\kappa^{\mathrm{set}}_{G,\mathcal S}(q)<\infty\qquad(q\geq2).
\]
\end{proposition}

\begin{proof}
Identify the subgroup generated by the first two standard basis vectors $\mathbf e_1,\mathbf e_2$ with $\Z^2$.
For every positive integer $\mathsf M$, apply \cref{lem:number-field-character-sum-tuples} to $(\chi(\mathsf M\mathbf e_1),\chi(\mathsf M\mathbf e_2))_{\chi\in\mathcal S}$.
Its counting bound and the pigeonhole principle, followed by removal of the intersection, give a two-block $(\mathsf M\Z^2,\mathcal S|_{\mathsf M\Z^2})$-PTE set solution $(B_1,B_2)$.

For a finite index set $I$, suppose that blocks $(A_\ell)_{\ell\in I}$ have already been obtained, and put
\[
\Omega:=\bigcup_{\ell\in I}A_\ell.
\]
Choose a positive integer $\mathsf M$ greater than the absolute values of all coordinate differences, so that
\[
\{\mathbf a-\mathbf a':\mathbf a,\mathbf a'\in\Omega\}\cap\mathsf M\Z^2=\{0\}.
\]
Take $B_1,B_2$ as above for this $\mathsf M$, and put
\[
A_\ell+B_j:=\{\mathbf a+\mathbf b:\mathbf a\in A_\ell,\ \mathbf b\in B_j\}.
\]
If $\mathbf a+\mathbf b=\mathbf a'+\mathbf b'$, then $\mathbf a-\mathbf a'=\mathbf b'-\mathbf b\in\mathsf M\Z^2$ gives $\mathbf a=\mathbf a'$ and $\mathbf b=\mathbf b'$.
Thus these sets are pairwise disjoint, and
\[
|A_\ell+B_j|=|A_\ell|\,|B_j|,\qquad
\sum_{\mathbf z\in A_\ell+B_j}\chi(\mathbf z)
=\left(\sum_{\mathbf a\in A_\ell}\chi(\mathbf a)\right)
\left(\sum_{\mathbf b\in B_j}\chi(\mathbf b)\right)
\qquad(\chi\in\mathcal S).
\]
Hence the number of sets can be doubled.
Starting from $A_1=\{0\}$, repeat this operation and choose any $q$ of the resulting sets.
\end{proof}

The remaining problem of existence conditions is the rank-one case.
For $G=(\Z,+)$, we have $\chi(k)=\chi(1)^k$ for $k\in\Z$.
Existence of a two-block set solution is equivalent to existence of a nonzero polynomial $\varphi\in\Z[X]$ with coefficients in $\{-1,0,1\}$ satisfying
\[
\varphi(1)=0,\qquad \varphi\bigl(\chi(1)\bigr)=0\quad(\chi\in\mathcal S).
\]
Indeed, after a common translation of the two blocks to sets of nonnegative integers, we obtain
\[
\varphi(X)=\sum_{a\in A_1}X^a-\sum_{a\in A_2}X^a.
\]
Conversely, take the sets of exponents whose coefficients are $1$ and $-1$.
For related problems on polynomial multiples, see Drungilas--Jankauskas--\v{S}iurys~\cite{DrungilasJankauskasSiurys2018}.

\begin{problem}\label{prob:characterize-free-abelian-pte-existence}
Let $G=(\Z,+)$, let $\mathbb K$ be a number field, and let $\mathcal S\subseteq\Hom(G,\mathbb K^\times)$ be finite with $\mathbf1\in\mathcal S$ and $\bigcap_{\chi\in\mathcal S}\ker\chi=\{0\}$.
Give necessary and sufficient conditions, in terms of the values $\{\chi(1):\chi\in\mathcal S\}$ at the generator, for
\[
\kappa^{\mathrm{set}}_{G,\mathcal S}(q)<\infty\qquad(q\geq2).
\]
Here the common size of the blocks may depend on $q$.
\end{problem}

\appendix

\section{Explicit Circle Constructions and Integer Identities}\label{app:index-data}

This appendix provides the constructions and supporting integer identities used in \cref{thm:minimum-circle-constructions}.
Each signed index list is ordered in increasing absolute value.
The exponent vectors used in the $4608$- and $9216$-point constructions are ordered lexicographically as specified in \eqref{eq:circle-eight-exponent-order}.

\subsection{The $2304$-point construction}\label{app:certificate-2304}
We use $R_4$ and $\psi_k$ as defined in \eqref{eq:roots-of-unity} and \eqref{eq:Pk}, respectively.
For the first seven primes ${\mathsf{p}}_i\equiv1\pmod4$ in \cref{subsec:circle-designs-pte}, choose the Gaussian primes ${p_{{\mathsf{p}}_i}}$ so that
\[
({p_{{\mathsf{p}}_1}},\ldots,{p_{{\mathsf{p}}_7}})
=
(2+\mathrm{i},3+2\mathrm{i},4+\mathrm{i},5+2\mathrm{i},6+\mathrm{i},5+4\mathrm{i},7+2\mathrm{i}).
\]
Put $g_i:={p_{{\mathsf{p}}_i}}/\overline{{p_{{\mathsf{p}}_i}}}$ for $1\leq i\leq7$, as in \eqref{eq:circle-prime-ratios}, and use the common denominator
\begin{equation}\label{eq:circle-seven-prime-denominator}
{\mathsf{b}}_7 = \prod_{i=1}^{7}{\mathsf{p}}_i
=5\cdot13\cdot17\cdot29\cdot37\cdot41\cdot53
=2576450045.
\end{equation}

\begin{lemma}\label{lem:gaussian-parameterization-seven}
\[
\begin{gathered}
\left\{
\frac{{\mathsf{x}_{1}}+{\mathsf{x}_{2}}\mathrm{i}}{{\mathsf{b}}_7}:
({\mathsf{x}_{1}},{\mathsf{x}_{2}})\in\Z^2,
{\mathsf{x}_{1}}^2+{\mathsf{x}_{2}}^2={\mathsf{b}}_7^2
\right\}
\\
=
\left\{
\mathrm i^j\prod_{i=1}^{7}
g_i^{k_i}:
j\in\{0,1,2,3\},
\mathbf k=(k_1,\ldots,k_7)\in\{-1,0,1\}^7
\right\}.
\end{gathered}
\]
The representation on the right is unique, and this set has $4\cdot3^7=8748$ points.
\end{lemma}
\begin{proof}
We use unique factorization in $\Z[\mathrm i]$; see \cite[Chapter~1, Section~4]{IrelandRosen1990}.
Let $({\mathsf{x}_{1}},{\mathsf{x}_{2}})\in\Z^2$ satisfy ${\mathsf{x}_{1}}^2+{\mathsf{x}_{2}}^2={\mathsf{b}}_7^2$. Then the norm $({\mathsf{x}_{1}}+{\mathsf{x}_{2}}\mathrm i)({\mathsf{x}_{1}}-{\mathsf{x}_{2}}\mathrm i)$ is ${\mathsf{b}}_7^2$.
The only Gaussian prime factors that can occur are ${p_{{\mathsf{p}}_i}}$ and $\overline{{p_{{\mathsf{p}}_i}}}$, and their exponents add to $2$ for each $i$.
Thus
\[
{\mathsf{x}_{1}}+{\mathsf{x}_{2}}\mathrm i=\mathrm i^j\prod_{i=1}^{7}
{p_{{\mathsf{p}}_i}}^{1+k_i}\overline{{p_{{\mathsf{p}}_i}}}^{\,1-k_i}
\qquad(j\in\{0,1,2,3\},\ k_i\in\{-1,0,1\}).
\]
Dividing by
\[
{\mathsf{b}}_7=\prod_{i=1}^{7}{p_{{\mathsf{p}}_i}}\overline{{p_{{\mathsf{p}}_i}}}
\]
gives the displayed representation.
Conversely, each product on the right has numerator in $\Z[\mathrm i]$ over denominator ${\mathsf{b}}_7$ and has norm one.
Unique factorization determines all $k_i$ and then $j$ uniquely, so there are $4\cdot3^7$ points.
\end{proof}
After removing the four points on the coordinate axes from the set in \cref{lem:gaussian-parameterization-seven}, the remaining points occur in quadruples under independent sign changes.
Hence
\begin{equation}\label{eq:circle-positive-pair-enumeration}
\{{\boldsymbol{\mathsf{x}}}=({\mathsf{x}_{1}},{\mathsf{x}_{2}}) \in \Z_{>0}^2 : {\mathsf{x}_{1}}^2+{\mathsf{x}_{2}}^2={\mathsf{b}}_7^2\} = \{{\boldsymbol{\mathsf{x}}}^{(j)}=({\mathsf{x}_{1}}^{(j)},{\mathsf{x}_{2}}^{(j)}) : 0 \leq j \leq 2185\},
\end{equation}
where the finite sequence $({\boldsymbol{\mathsf{x}}}^{(j)})_{j=0}^{2185}$ is ordered by ${\mathsf{x}_{1}}^{(0)}<\cdots<{\mathsf{x}_{1}}^{(2185)}$.
The number of pairs is
\[
\frac{8748-4}{4}=2186.
\]
Since ${\mathsf{x}_{2}}$ decreases as ${\mathsf{x}_{1}}$ increases, coordinate exchange $({\mathsf{x}_{1}},{\mathsf{x}_{2}})\mapsto({\mathsf{x}_{2}},{\mathsf{x}_{1}})$ reverses this ordering.
For $0\leq j\leq2185$, coordinate exchange gives
\begin{equation}\label{eq:circle-local-parameters}
({\mathsf{x}_{1}}^{(2185-j)},{\mathsf{x}_{2}}^{(2185-j)})=({\mathsf{x}_{2}}^{(j)},{\mathsf{x}_{1}}^{(j)}).
\end{equation}
For $0\leq j\leq2185$ and ${\mathsf{e}}\in\{-1,1\}$, define the following subsets of $\TQ$:
\begin{equation}\label{eq:circle-local-blocks}
\mathfrak Q_j^{(2,{\mathsf{e}})} := \left\{ \frac{{\mathsf{e}} {\mathsf{x}_{1}}^{(j)}+{\mathsf{x}_{2}}^{(j)}\mathrm{i}}{{\mathsf{b}}_7}, \frac{{\mathsf{e}} {\mathsf{x}_{1}}^{(j)}-{\mathsf{x}_{2}}^{(j)}\mathrm{i}}{{\mathsf{b}}_7} \right\},
\qquad \mathfrak Q_j^{(4)} := \left\{ \frac{{\mathsf{e}} {\mathsf{x}_{1}}^{(j)}+{\mathsf{e}}' {\mathsf{x}_{2}}^{(j)}\mathrm{i}}{{\mathsf{b}}_7}: {\mathsf{e}},{\mathsf{e}}'\in\{-1,1\} \right\}.
\end{equation}

\begin{lemma}\label{lem:local-moment-formulas-four}
For $0\leq j\leq2185$ and ${\mathsf{e}}\in\{-1,1\}$, the sets in \eqref{eq:circle-local-blocks} have the following moments:
\begingroup
\setlength{\arraycolsep}{2pt}
\renewcommand{\arraystretch}{2.2}
\[
\begin{array}{c|cccc}
 & \psi_1 & \psi_2 & \psi_3 & \psi_4 \\
\hline
\mathfrak Q_j^{(2,{\mathsf{e}})}
& \frac{2{\mathsf{e}} {\mathsf{x}_{1}}^{(j)}}{{\mathsf{b}}_7}
& \frac{2(({\mathsf{x}_{1}}^{(j)})^2-({\mathsf{x}_{2}}^{(j)})^2)}{{\mathsf{b}}_7^2}
& \frac{2{\mathsf{e}}(4({\mathsf{x}_{1}}^{(j)})^3-3{\mathsf{x}_{1}}^{(j)}{\mathsf{b}}_7^2)}{{\mathsf{b}}_7^3}
& \frac{4(({\mathsf{x}_{1}}^{(j)})^2-({\mathsf{x}_{2}}^{(j)})^2)^2}{{\mathsf{b}}_7^4}-2 \\
\mathfrak Q_j^{(4)}
& 0
& \frac{4(({\mathsf{x}_{1}}^{(j)})^2-({\mathsf{x}_{2}}^{(j)})^2)}{{\mathsf{b}}_7^2}
& 0
& \frac{8(({\mathsf{x}_{1}}^{(j)})^2-({\mathsf{x}_{2}}^{(j)})^2)^2}{{\mathsf{b}}_7^4}-4 \\
\mathfrak Q_j^{(4)}\sqcup \mathfrak Q_{2185-j}^{(4)}
& 0 & 0 & 0
& \frac{16(({\mathsf{x}_{1}}^{(j)})^2-({\mathsf{x}_{2}}^{(j)})^2)^2}{{\mathsf{b}}_7^4}-8
\end{array}
\]
\endgroup
\end{lemma}
\begin{proof}
For $\mathfrak Q_j^{(2,{\mathsf{e}})}$, expand the conjugate pair in \eqref{eq:circle-local-blocks} and use $({\mathsf{x}_{1}}^{(j)})^2+({\mathsf{x}_{2}}^{(j)})^2={\mathsf{b}}_7^2$.
Since $\mathfrak Q_j^{(4)}=\mathfrak Q_j^{(2,1)}\sqcup \mathfrak Q_j^{(2,-1)}$, we have $\psi_k(\mathfrak Q_j^{(4)})=\psi_k(\mathfrak Q_j^{(2,1)})+\psi_k(\mathfrak Q_j^{(2,-1)})$ for $1\leq k\leq4$.
Since $j\ne 2185-j$, the sets $\mathfrak Q_j^{(4)}$ and $\mathfrak Q_{2185-j}^{(4)}$ in \eqref{eq:circle-local-blocks} are disjoint.
By \eqref{eq:circle-local-parameters}, coordinate exchange negates $({\mathsf{x}_{1}}^{(j)})^2-({\mathsf{x}_{2}}^{(j)})^2$. Adding the moments of the two disjoint four-point sets therefore gives the stated moments of their union.
\end{proof}
The following lists specify the index sets
\[
L_4^{(2)}\subseteq\{-2185,\ldots,-1,1,\ldots,2185\},
\qquad L_4^{(4)},L_4^{(8)}\subseteq\{0,\ldots,2185\}.
\]
An index $i\in L_4^{(2)}$ selects $\mathfrak Q_{|i|}^{(2,\operatorname{sgn}(i))}$, an index $j\in L_4^{(4)}$ selects $\mathfrak Q_j^{(4)}$, and an index $j\in L_4^{(8)}$ selects $\mathfrak Q_j^{(4)}\sqcup \mathfrak Q_{2185-j}^{(4)}$.
\compactindexlist{L_4^{(2)}}{%
9,10,-24,-52,-55,-57,-75,79,-86,-92,94,116,-152,174,183,-193,199,-206,-225,-228,231,-245,250,-282,
-295,296,-301,-306,-328,337,356,-365,403,405,424,-432,440,450,451,455,467,-471,-482,-500,-519,-522,-534,-538,
-569,577,-587,600,602,609,-623,624,-634,-639,-642,662,678,689,-691,-694,698,-701,-706,734,-749,751,757,-760,
-765,778,813,-821,-829,833,-839,847,871,-878,-880,-882,895,-902,939,970,-1017,1031,1070,1078,1079,-1085,-1087,1106,
-1107,1123,1145,-1150,-1155,1172,1182,-1186,-1189,-1206,1207,1219,1236,1240,-1275,1278,-1299,1317,-1320,1350,1352,1372,-1378,-1413,
-1428,-1454,1459,-1479,1498,-1515,1524,-1551,-1555,1579,1585,1592,-1598,1621,-1626,-1628,1631,-1633,1647,-1649,-1655,-1664,1669,1679,
-1703,-1716,1735,1749,1753,-1756,-1761,1786,1794,1797,1803,1827,1840,1854,-1884,-1889,1893,1903,-1909,-1945,1957,-1960,-1966,-1969,
1972,-1975,1986,-1992,-2008,-2011,2018,-2021,-2031,2037,-2039,-2044,2059,-2065,-2086,2089,2091,-2109,2110,2126,-2128,-2137,2166,2179%
}
\compactindexlist{L_4^{(4)}}{%
1,4,47,102,103,134,175,187,240,302,317,323,332,333,389,421,429,438,464,502,503,529,531,552,
557,564,565,571,590,592,599,651,652,685,693,759,761,776,789,799,802,817,912,921,934,967,997,1151,
1188,1196,1234,1242,1247,1254,1267,1292,1301,1302,1304,1336,1338,1343,1349,1375,1386,1395,1508,1541,1630,1651,1673,1682,
1689,1696,1702,1709,1750,1754,1771,1788,1789,1832,1835,1878,1886,1914,1917,1923,1948,1965,1997,2007,2035,2067,2111,2177%
}
\compactindexlist{L_4^{(8)}}{%
3,13,21,25,34,35,65,81,83,85,87,95,100,104,105,115,117,121,122,128,132,133,144,157,
159,162,169,179,181,182,184,186,195,200,222,224,226,227,229,242,255,258,263,264,270,272,274,278,
283,286,298,308,310,313,319,324,329,340,342,359,360,361,367,371,377,380,381,387,392,393,394,402,
410,411,416,423,426,428,437,445,448,452,453,454,457,458,468,498,504,505,513,523,526,535,541,573,
574,578,582,588,594,596,604,610,611,622,629,640,648,650,658,666,668,673,674,676,682,683,688,699,
707,709,710,714,718,720,722,728,730,737,739,763,764,768,781,784,785,803,804,806,814,816,820,823,
826,828,832,845,848,850,856,863,866,870,872,874,876,887,906,908,909,913,917,922,924,927,940,944,
955,960,963,964,971,973,983,986,998,1004,1014,1015,1020,1023,1027,1037,1041,1046,1063,1077,1083,1084,1088,1089%
}
Define $\mathfrak X_4$ by
\[
\mathfrak X_4
:=
\bigsqcup_{i\in L_4^{(2)}}\mathfrak Q_{|i|}^{(2,\operatorname{sgn}(i))}
\sqcup
\bigsqcup_{j\in L_4^{(4)}}\mathfrak Q_j^{(4)}
\sqcup
\bigsqcup_{j\in L_4^{(8)}}(\mathfrak Q_j^{(4)}\sqcup \mathfrak Q_{2185-j}^{(4)}).
\]

\begin{proposition}\label{prop:certificate-four}
The sets $L_4^{(2)},L_4^{(4)},L_4^{(8)}$ satisfy
\[
\lvert L_4^{(2)}\rvert=192,\qquad
\lvert L_4^{(4)}\rvert=96,\qquad
\lvert L_4^{(8)}\rvert=192.
\]
Also, all indices selected by $L_4^{(2)}$, $L_4^{(4)}$, and $L_4^{(8)}$ are pairwise distinct, where $i\in L_4^{(2)}$ selects $|i|$, $j\in L_4^{(4)}$ selects $j$, and $j\in L_4^{(8)}$ selects both $j$ and $2185-j$.
Moreover,
\begin{gather*}
\sum_{i\in L_4^{(2)}}\operatorname{sgn}(i){\mathsf{x}_{1}}^{(|i|)}
=
\sum_{i\in L_4^{(2)}}\operatorname{sgn}(i)({\mathsf{x}_{1}}^{(|i|)})^3
=
0,\\
2\sum_{i\in L_4^{(2)}}(({\mathsf{x}_{1}}^{(|i|)})^2-({\mathsf{x}_{2}}^{(|i|)})^2)
+
4\sum_{j\in L_4^{(4)}}(({\mathsf{x}_{1}}^{(j)})^2-({\mathsf{x}_{2}}^{(j)})^2)
=
0,\\
2\sum_{i\in L_4^{(2)}}(({\mathsf{x}_{1}}^{(|i|)})^2-({\mathsf{x}_{2}}^{(|i|)})^2)^2 +4\sum_{j\in L_4^{(4)}}(({\mathsf{x}_{1}}^{(j)})^2-({\mathsf{x}_{2}}^{(j)})^2)^2 +8\sum_{j\in L_4^{(8)}}(({\mathsf{x}_{1}}^{(j)})^2-({\mathsf{x}_{2}}^{(j)})^2)^2 =1152{\mathsf{b}}_7^4.
\end{gather*}
\end{proposition}
\begin{proof}
The sizes and distinctness of the indices follow by counting and comparing the listed integers.
To evaluate the four sums, use \cref{lem:gaussian-parameterization-seven} to enumerate the positive integer pairs in \eqref{eq:circle-positive-pair-enumeration}, with coordinate exchange as in \eqref{eq:circle-local-parameters}.
Integer summation with the listed indices and signs gives $0$, $0$, $0$, and $1152{\mathsf{b}}_7^4$, respectively.
\end{proof}

\begin{proof}[Proof of \cref{thm:minimum-circle-constructions}\textup{(i)}]
By \cref{prop:certificate-four}, the sets defining $\mathfrak X_4$ are pairwise disjoint and
\[
|\mathfrak X_4|=2|L_4^{(2)}|+4|L_4^{(4)}|+8|L_4^{(8)}|=2304.
\]
In \cref{lem:local-moment-formulas-four}, only the sets $\mathfrak Q_{|i|}^{(2,\operatorname{sgn}(i))}$ contribute to the first and third moments.
The equalities involving ${\mathsf{x}_{1}}^{(|i|)}$ and $({\mathsf{x}_{1}}^{(|i|)})^3$ in \cref{prop:certificate-four}, together with \cref{lem:local-moment-formulas-four}, give $\psi_1(\mathfrak X_4)=\psi_3(\mathfrak X_4)=0$.
For the second moment, the third identity in \cref{prop:certificate-four} and \cref{lem:local-moment-formulas-four} give $\psi_2(\mathfrak X_4)=0$.
For the fourth moment, the nonconstant terms are twice the left-hand side of the fourth identity in \cref{prop:certificate-four}, divided by ${\mathsf{b}}_7^4$; the constant terms sum to $-|\mathfrak X_4|$. Thus
\[
\psi_4(\mathfrak X_4)=\frac{2\cdot1152{\mathsf{b}}_7^4}{{\mathsf{b}}_7^4}-2304=0.
\]
Therefore $\mathfrak X_4$ is a spherical $4$-design by \cref{lem:circle-moment-criterion}.
\end{proof}
\subsection{The $4608$- and $9216$-point constructions}\label{app:certificate-4608-9216}
Retain the Gaussian primes ${p_{{\mathsf{p}}_1}},\ldots,{p_{{\mathsf{p}}_7}}$ from \cref{app:certificate-2304}, and choose ${p_{{\mathsf{p}}_8}}={p_{61}}=5+6\mathrm i$.
Then
\begin{equation}\label{eq:circle-eight-prime-denominator}
{\mathsf{b}}_8=61{\mathsf{b}}_7=157163452745.
\end{equation}
Recall from \eqref{eq:circle-prime-ratios} and \eqref{eq:circle-products} that $g_i={p_{{\mathsf{p}}_i}}/\overline{{p_{{\mathsf{p}}_i}}}$ and ${g_{\mathbf k}}=\prod_{i=1}^8g_i^{k_i}$.
Specializing \eqref{eq:circle-exponent-set} to eight coordinates gives
\[
J_8=\{\mathbf k\in\{-1,0,1\}^8\setminus\{\mathbf0\}:
\exists j\in\{1,\ldots,8\}\ (k_j=1,\ k_i=0\text{ for }1\leq i<j)\}.
\]
By \eqref{eq:circle-exponent-set}, exactly one of $\mathbf k$ and $-\mathbf k$ belongs to $J_8$ for each nonzero $\mathbf k\in\{-1,0,1\}^8$, so $|J_8|=(3^8-1)/2=3280$.
Enumerate the elements of $J_8$ as $\mathbf k^{(0)},\ldots,\mathbf k^{(3279)}$ in lexicographic order, using the coordinate order $1,\ldots,8$ and the entry order $-1<0<1$:
\begin{equation}\label{eq:circle-eight-exponent-order}
\mathbf k^{(0)}<\cdots<\mathbf k^{(3279)}.
\end{equation}
For $\mathbf k\in J_8$, define
\[
\mathfrak Q_{\mathbf k}^{(8)}:=\{\mathrm i^j {g_{\mathbf k}},\mathrm i^j\overline{{g_{\mathbf k}}}:j\in\{0,1,2,3\}\}.
\]
Suppose that $\mathbf k,\mathbf k'\in J_8$, $j,j'\in\{0,1,2,3\}$, and ${\mathsf{e}}',{\mathsf{e}}\in\{-1,1\}$ satisfy $\mathrm i^j {g_{\mathbf k}}^{{\mathsf{e}}'}=\mathrm i^{j'} {g_{\mathbf k'}}^{\mathsf{e}}$.
Taking valuations at ${p_{{\mathsf{p}}_1}},\ldots,{p_{{\mathsf{p}}_8}}$ gives ${\mathsf{e}}'\mathbf k={\mathsf{e}}\mathbf k'$.
Since the first nonzero coordinate of each vector in $J_8$ is $1$, it follows that $\mathbf k=\mathbf k'$, ${\mathsf{e}}'={\mathsf{e}}$, and then $j=j'$.
Thus every $\mathfrak Q_{\mathbf k}^{(8)}$ has eight points and the sets $\mathfrak Q_{\mathbf k}^{(8)}$ are pairwise disjoint.
For a nonzero integer $i$ with $|i|\leq3279$, its sign selects the point
\begin{equation}\label{eq:circle-signed-point-choice}
\mathrm i^{(1-\operatorname{sgn}(i))/2}{g_{\mathbf k^{(|i|)}}}\in \mathfrak Q_{\mathbf k^{(|i|)}}^{(8)}.
\end{equation}
The following list specifies the signed index set
\[
L_5\subseteq\{-3279,\ldots,-1,1,\ldots,3279\}.
\]
For $i\in L_5$, the absolute value $|i|$ selects the exponent vector $\mathbf k^{(|i|)}$ in \eqref{eq:circle-eight-exponent-order}, and its sign selects the antipodal four-point subset containing the point in \eqref{eq:circle-signed-point-choice}.
\compactindexlist{L_5}{%
2,5,-8,-11,-14,16,20,22,-23,-24,26,-29,30,-35,36,-38,39,45,46,-47,-50,53,59,-66,
70,71,-74,-77,83,-86,89,91,92,93,-95,-100,-101,102,-107,110,-111,-113,-116,-125,128,-131,-132,-137,
141,-142,143,-148,-149,152,-153,-154,-155,-158,-161,164,-167,171,173,-176,179,-182,-185,-191,-197,-200,-203,206,
-209,215,218,221,-222,224,227,-230,-233,-236,239,242,-243,245,250,254,-257,260,262,266,269,272,-274,-275,
-278,281,-284,288,290,-291,-293,299,-302,305,-309,311,-313,-314,315,316,-317,319,320,-323,-325,-332,335,338,
341,-344,-347,348,350,353,-356,358,359,361,362,-365,366,-368,-376,-377,-379,-380,-383,-385,386,-389,392,395,
400,401,-404,407,-410,-413,-416,-417,-420,-422,-424,428,-429,431,434,-439,440,-443,-445,446,-447,452,-454,-455,
458,-464,467,-470,473,-476,477,479,-482,485,-490,-491,-494,-506,-509,-512,-515,-516,-518,-519,521,-524,527,529,
530,533,-535,-536,-538,-542,-545,-550,-551,-552,554,557,560,-563,-569,-575,578,581,-584,-588,593,596,-598,599,
602,-605,609,-611,-614,623,629,-632,-635,-641,-642,644,646,-647,650,-653,-655,-656,-659,662,665,668,671,672,
673,-674,-677,-678,-680,-686,-687,689,-691,-692,-695,-701,704,-705,-708,710,713,-714,716,-718,-719,-723,728,731,
-734,-740,-752,753,-755,-757,-758,-759,-761,764,-768,770,772,-773,776,-782,785,-788,792,794,-800,803,806,809,
810,-811,-812,815,818,821,-824,827,828,-830,836,839,-842,-845,848,851,-854,857,-858,-866,-869,872,-875,878,
881,-884,889,-898,899,-902,905,908,-911,914,-916,917,920,-923,-924,-929,932,933,935,938,950,-953,956,959,
960,-965,968,-971,-973,-974,977,-980,983,989,992,995,-1001,1004,1010,-1012,-1016,-1019,1025,1034,1037,1043,1049,1060,
-1064,1070,1073,-1076,-1079,-1080,1082,-1084,-1085,1088,1089,1091,-1093,-1096,-1100,-1101,1102,1103,1106,1109,-1112,-1114,1118,-1125,
-1127,1130,-1133,-1136,1137,1138,-1141,-1145,-1146,1148,-1154,1157,-1159,-1160,-1163,-1165,-1166,-1169,1172,1175,-1177,-1184,-1187,1190,
1196,-1199,-1202,-1205,1211,1214,-1216,1217,-1219,-1220,-1223,-1224,-1225,1226,1229,1232,1238,1241,1244,1246,1247,1250,-1252,-1256,
1259,1262,1264,1265,-1271,1272,-1274,1277,1281,1283,1286,1292,1295,1298,1301,-1304,1307,-1316,1319,-1322,1330,-1331,1333,1337,
-1340,-1341,-1343,1345,-1346,-1347,1348,1349,1351,-1355,-1359,-1361,-1364,1368,1370,-1373,1376,1379,-1386,1388,1389,-1390,-1391,1392,
1396,-1397,-1398,-1400,-1403,1409,1414,1416,-1418,-1420,-1423,-1424,-1427,-1430,1432,-1433,1436,1439,-1442,1446,-1448,1449,1451,1454,
1457,-1460,-1466,-1469,-1472,-1473,1475,-1478,1481,-1482,1484,-1487,-1490,-1491,1499,-1502,1503,1504,-1505,1508,1509,-1511,1514,-1516,
1517,-1520,1523,1526,-1529,-1535,-1538,-1539,-1540,1541,1543,-1546,-1547,-1549,1550,-1556,1557,1559,-1560,1565,1567,1574,1577,-1580,
1589,1592,1594,1595,-1597,1598,1602,-1604,-1607,1608,-1610,1612,-1613,1619,1620,-1622,1625,-1628,-1631,1634,1637,1640,1646,1649,
-1652,1655,-1658,-1664,-1667,-1670,1676,-1677,1679,-1685,-1688,-1691,-1694,-1696,-1697,-1700,1702,1703,-1705,1706,-1709,-1712,1721,-1724,
-1727,-1730,1731,-1733,-1736,1739,-1742,-1743,-1744,-1745,-1748,1749,1750,-1754,1755,1756,-1763,-1764,1766,-1769,1772,1775,-1780,-1781,
-1784,-1787,-1792,1793,1796,1802,-1804,1808,1811,-1814,1817,-1818,-1820,-1823,1825,1826,1827,-1829,1832,-1835,-1838,1841,-1844,-1847,
-1850,1853,-1856,1859,-1860,-1861,-1865,-1867,1869,-1871,-1874,1877,1881,-1887,1889,1895,1898,1899,-1901,-1907,-1910,1911,-1916,-1917,
-1919,1920,1925,-1928,1931,-1934,-1936,-1937,1939,-1940,-1942,-1943,-1952,-1953,-1955,-1956,-1957,1958,1962,1964,-1967,1970,-1976,1979,
-1982,1983,1985,-1987,-1988,1991,-1994,1998,-1999,2000,2003,2006,-2009,2012,2014,2015,2017,2018,-2021,2024,2025,2027,2030,2033,
2036,2042,-2043,2050,-2051,2054,2059,2066,-2072,-2075,2078,-2084,-2087,2089,2090,-2092,-2093,-2098,2099,2105,-2111,2113,2114,2117,
-2120,-2123,2127,2129,2131,-2135,-2141,-2144,2147,-2150,-2153,2154,-2156,2159,-2161,-2162,-2165,2171,-2180,-2181,2183,-2186,-2189,-2192,
-2195,2197,2198,-2200,2201,-2202,-2204,-2207,-2208,-2211,-2213,2214,-2219,-2222,-2225,-2228,2230,-2231,-2234,-2237,2243,-2245,-2249,-2255,
-2258,2259,-2261,2262,-2264,-2267,-2269,-2270,2271,2273,2274,2275,2276,2279,-2282,2283,-2287,2291,2294,-2297,2298,-2306,-2307,2309,
-2315,-2318,-2320,2324,2333,-2339,2342,2343,-2345,-2346,2348,2350,-2351,-2354,2355,-2357,-2360,-2363,-2366,2372,2375,2379,2381,-2384,
2387,-2390,2396,2399,-2402,-2405,2408,-2411,-2412,-2414,2417,-2420,2423,-2424,2429,-2430,2432,-2435,-2438,2441,-2444,2453,-2455,-2456,
2459,2462,2465,2468,2474,-2475,-2477,-2480,-2483,2486,2489,-2491,2492,2493,2494,-2495,2498,-2501,2504,2510,2513,-2515,2516,2519,
2521,-2522,2528,2529,2531,2534,2537,-2540,-2542,2543,-2546,2548,-2558,-2559,-2561,2563,-2564,2567,-2573,2576,-2579,-2580,2585,-2588,
2589,-2591,-2600,-2603,-2605,-2609,2615,2618,2621,2624,2625,2627,2630,2632,-2636,-2642,2645,2651,-2654,-2657,2659,-2660,-2661,2668,
-2669,2675,-2678,-2681,2686,-2687,-2690,2699,-2702,2704,2705,-2708,2711,-2712,-2714,-2717,2720,2723,2726,2729,-2730,-2735,2737,-2738,
2742,2747,-2750,2753,2754,-2756,-2759,2762,-2765,2768,-2771,2777,2780,-2782,-2783,-2789,-2792,-2795,-2798,2801,2807,2808,-2810,-2813,
2816,-2818,-2825,2828,-2829,-2830,2831,-2834,2837,2842,-2846,-2852,-2855,2858,-2861,2864,2865,-2867,-2868,2870,-2873,-2876,-2879,-2888,
-2893,2896,-2897,2900,2902,-2903,-2906,2908,2909,2915,2927,2930,-2933,2936,2937,-2938,2939,-2942,2945,2948,-2953,-2954,-2960,2968,
-2969,-2972,2978,2981,-2984,2987,-2989,2993,-2996,-2999,3002,-3005,3011,3012,-3017,3018,3022,3023,3026,3029,-3032,-3035,3038,-3039,
3041,3044,3047,-3048,3049,3050,3053,-3054,-3056,-3059,-3060,-3062,3063,-3065,3068,-3071,-3077,3080,-3082,-3083,-3086,3087,3089,3095,
-3096,3097,-3098,3099,3100,3103,-3104,3107,3110,3112,-3113,3116,-3119,3122,3125,-3128,3131,3132,-3133,-3134,3136,-3137,3138,3140,
-3149,3152,3155,3157,3158,3161,-3164,3170,3172,3175,3176,-3177,-3179,3182,3188,3191,3194,3197,-3198,3200,3203,3206,3212,-3215,
-3218,3219,-3221,3224,3233,-3234,3237,-3238,-3239,-3242,-3245,3248,-3253,3254,-3257,-3260,3263,-3265,-3266,3269,-3270,-3271,3272,-3275%
}
The selected eight-point sets and their antipodal four-point subsets define
\[
\mathfrak X_7
:=
\bigsqcup_{i\in L_5}
\mathfrak Q_{\mathbf k^{(|i|)}}^{(8)},
\qquad
\mathfrak X_5
:=
\bigsqcup_{i\in L_5}
\left\{{\mathsf{e}}\mathrm i^{(1-\operatorname{sgn}(i))/2}{g_{\mathbf k^{(|i|)}}}^{{\mathsf{e}}'}:{\mathsf{e}},{\mathsf{e}}'\in\{-1,1\}\right\}.
\]
\begin{proposition}\label{prop:certificate-five-seven}
The set $L_5$ has $1152$ elements with pairwise distinct absolute values.
For $i\in L_5$, write
\[
{g_{\mathbf k^{(|i|)}}}=({\mathsf{x}_{1}}^{\prime\,(i)}+{\mathsf{x}_{2}}^{\prime\,(i)}\mathrm{i})/{\mathsf{b}}_8
\]
with ${\mathsf{x}_{1}}^{\prime\,(i)},{\mathsf{x}_{2}}^{\prime\,(i)}\in\Z$. Then
\[
\sum_{i\in L_5} 2\left(({\mathsf{x}_{1}}^{\prime\,(i)})^4-6({\mathsf{x}_{1}}^{\prime\,(i)})^2({\mathsf{x}_{2}}^{\prime\,(i)})^2+({\mathsf{x}_{2}}^{\prime\,(i)})^4\right) = 0, \qquad \sum_{i\in L_5} \operatorname{sgn}(i) \left(({\mathsf{x}_{1}}^{\prime\,(i)})^2-({\mathsf{x}_{2}}^{\prime\,(i)})^2\right) = 0.
\]
Consequently,
\[
\sum_{i\in L_5}
\left({g_{\mathbf k^{(|i|)}}}^4+\overline{{g_{\mathbf k^{(|i|)}}}}^{\,4}\right)
=
0,
\qquad
\sum_{i\in L_5}
\operatorname{sgn}(i)\Re({g_{\mathbf k^{(|i|)}}}^2)
=
0.
\]
\end{proposition}
\begin{proof}
Counting and comparing the listed integers gives $|L_5|=1152$ and pairwise distinct absolute values.
For each $i\in L_5$, expand
\[
{\mathsf{x}_{1}}^{\prime\,(i)}+{\mathsf{x}_{2}}^{\prime\,(i)}\mathrm{i}
=
\prod_{j=1}^{8}
{p_{{\mathsf{p}}_j}}^{1+k_j^{(|i|)}}
\overline{{p_{{\mathsf{p}}_j}}}^{\,1-k_j^{(|i|)}}
\]
in $\Z[\mathrm{i}]$.
Summing
\[
2(({\mathsf{x}_{1}}^{\prime\,(i)})^4-6({\mathsf{x}_{1}}^{\prime\,(i)})^2({\mathsf{x}_{2}}^{\prime\,(i)})^2+({\mathsf{x}_{2}}^{\prime\,(i)})^4)
\qquad\text{and}\qquad \operatorname{sgn}(i)(({\mathsf{x}_{1}}^{\prime\,(i)})^2-({\mathsf{x}_{2}}^{\prime\,(i)})^2)
\]
over $i\in L_5$ gives the two integer identities.
Dividing the first identity by ${\mathsf{b}}_8^4$ and the second by ${\mathsf{b}}_8^2$ gives the last two identities.
\end{proof}
\begin{proof}[Proof of \cref{thm:minimum-circle-constructions}\textup{(ii)}]
The disjointness proved for the sets $\mathfrak Q_{\mathbf k}^{(8)}$, together with \cref{prop:certificate-five-seven}, shows that $\mathfrak X_7$ consists of $8\cdot1152=9216$ distinct points.
Each four-point set defining $\mathfrak X_5$ is a subset of the corresponding $\mathfrak Q_{\mathbf k^{(|i|)}}^{(8)}$, so these sets are pairwise disjoint and $|\mathfrak X_5|=4\cdot1152=4608$.

Each set $\mathfrak Q_{\mathbf k^{(|i|)}}^{(8)}$ is invariant under multiplication by $R_4$, so its moments of degrees $1,2,3,5,6,7$ vanish.
Its fourth moment is
\[
\sum_{j=0}^3\left((\mathrm i^j {g_{\mathbf k^{(|i|)}}})^4+(\mathrm i^j\overline{{g_{\mathbf k^{(|i|)}}}})^4\right)
=
4\left({g_{\mathbf k^{(|i|)}}}^4+\overline{{g_{\mathbf k^{(|i|)}}}}^{\,4}\right).
\]
Summing over $i\in L_5$ and applying \cref{prop:certificate-five-seven} shows that the fourth moment of $\mathfrak X_7$ is zero.
Hence $\mathfrak X_7$ is a spherical $7$-design by \cref{lem:circle-moment-criterion}.

Each four-point set defining $\mathfrak X_5$ is invariant under multiplication by $-1$, so its odd moments vanish.
For the second and fourth moments, we have
\[
\begin{aligned}
\sum_{{\mathsf{e}},{\mathsf{e}}'\in\{-1,1\}}
\left({\mathsf{e}}\mathrm i^{(1-\operatorname{sgn}(i))/2}{g_{\mathbf k^{(|i|)}}}^{{\mathsf{e}}'}\right)^2
&=4\operatorname{sgn}(i)\Re({g_{\mathbf k^{(|i|)}}}^2),\\
\sum_{{\mathsf{e}},{\mathsf{e}}'\in\{-1,1\}}
\left({\mathsf{e}}\mathrm i^{(1-\operatorname{sgn}(i))/2}{g_{\mathbf k^{(|i|)}}}^{{\mathsf{e}}'}\right)^4
&=4\Re({g_{\mathbf k^{(|i|)}}}^4).
\end{aligned}
\]
Summing over $i\in L_5$ and applying \cref{prop:certificate-five-seven} gives $\psi_2(\mathfrak X_5)=\psi_4(\mathfrak X_5)=0$.
Hence $\mathfrak X_5$ is a spherical $5$-design by \cref{lem:circle-moment-criterion}.
\end{proof}

\section{Proofs for the Degree-Five Family}

\subsection{Proof of Theorem~\ref{thm:new-degree-five-family}}\label{app:degree-five-family}

The following lemma gives a sufficient condition for two six-element multisets to have equal power sums through degree five.
We use it in the proof of \cref{thm:new-degree-five-family}.

\begin{lemma}\label{lem:degree-five-block-criterion}
Consider
\begin{equation}
 \begin{gathered}
\mathbb{K} = \mathbb{Q}, \qquad G = \{\pm1\}^2, \qquad \mathcal{S} = \{\mathbf{1},\chi_1,\chi_2\}, \\
 A'_1=\{(1,1),(-1,-1)\}, \qquad A'_2=\{(1,-1),(-1,1)\},
\end{gathered}
 \label{eq:c2-source-concise}
\end{equation}
where $G$ has componentwise multiplication, ${\mathsf{e}},{\mathsf{e}}'\in\{-1,1\}$, and \(\chi_1({\mathsf{e}},{\mathsf{e}}')={\mathsf{e}}\), \(\chi_2({\mathsf{e}},{\mathsf{e}}')={\mathsf{e}}'\).
The sets $A'_1,A'_2$ are disjoint, have size $2$, and satisfy
\[
\sum_{g\in A'_1}\chi_j(g)=\sum_{g\in A'_2}\chi_j(g)=0\qquad(j=1,2).
\]
Thus $(A'_1,A'_2)$ is a $(G,\mathcal S)$-PTE solution of size $2$.
Let $h_i,a_{i,1},a_{i,2} \in \mathbb{Q}$ $(1\leq i\leq3)$.
Define \(\upsilon:G\to\operatorname{MSet}(\mathbb Q)\) by letting $\upsilon({\mathsf{e}},{\mathsf{e}}')$ be the three-element multiset of the values $h_i+{\mathsf{e}}a_{i,1}+{\mathsf{e}}'a_{i,2}$ for $1\leq i\leq3$, counted with multiplicity.
The subscripts $i$ on $h_i,a_{i,1},a_{i,2}$ are read cyclically modulo $3$.
Suppose that
\begin{equation}
\begin{gathered}
 \sum_{i=1}^3 h_i=0,\qquad
 (a_{i,1}a_{i,2})_{i=1}^3\in\operatorname{span}_{\mathbb Q}\{(h_2-h_3,h_3-h_1,h_1-h_2)\},\\
 a_{i,1}^2+a_{i,2}^2+3h_i^2=a_{1,1}^2+a_{1,2}^2+3h_1^2\qquad(1\leq i\leq3).
\end{gathered}
 \label{eq:c2-criterion-concise}
\end{equation}
Then $\upsilon_*A'_1$ and $\upsilon_*A'_2$ each have size six and have equal power sums through degree five.
\end{lemma}

\begin{proof}
The four characters $\mathbf1,\chi_1,\chi_2,\chi_1\chi_2$ form a $\mathbb Q$-basis of the functions $G\to\mathbb Q$.
For $0\leq k\leq5$, let $d_k$ be the coefficient of $\chi_1\chi_2$ in
\[
(\psi_k\circ\upsilon)({\mathsf{e}},{\mathsf{e}}')=\sum_{i=1}^3(h_i+{\mathsf{e}}a_{i,1}+{\mathsf{e}}'a_{i,2})^k.
\]
Then $d_0=d_1=0$, and direct expansion gives
\[
\begin{gathered}
d_2=2\sum_{i=1}^3a_{i,1}a_{i,2},\qquad
d_3=6\sum_{i=1}^3h_ia_{i,1}a_{i,2},\\
d_4=4\sum_{i=1}^3a_{i,1}a_{i,2}(3h_i^2+a_{i,1}^2+a_{i,2}^2),\qquad
d_5=20\sum_{i=1}^3h_ia_{i,1}a_{i,2}(h_i^2+a_{i,1}^2+a_{i,2}^2).
\end{gathered}
\]
Under \eqref{eq:c2-criterion-concise}, $d_2=d_3=d_4=0$ because
\[
\sum_{i=1}^3(h_{i+1}-h_{i+2})
=\sum_{i=1}^3h_i(h_{i+1}-h_{i+2})=0,
\]
and $3h_i^2+a_{i,1}^2+a_{i,2}^2$ is independent of $i$.
In the expression for $d_5$, the common value of $3h_i^2+a_{i,1}^2+a_{i,2}^2$ contributes zero because $\sum_{i=1}^3h_ia_{i,1}a_{i,2}=0$.
The remaining sum $-40\sum_{i=1}^3h_i^3a_{i,1}a_{i,2}$ vanishes by the proportionality condition in \eqref{eq:c2-criterion-concise} and the identity
\[
\sum_{i=1}^3h_i^3(h_{i+1}-h_{i+2})
=-(h_1+h_2+h_3)(h_1-h_2)(h_2-h_3)(h_3-h_1)=0.
\]
Consequently, $d_k=0$ for $0\leq k\leq5$, so $\psi_k\circ\upsilon\in\operatorname{span}_{\mathbb Q}\{\mathbf1,\chi_1,\chi_2\}$.
The conclusion follows from \cref{thm:ideal-block-transfer}(a).
\end{proof}

For a six-element multiset $A$ on $\mathbb Q$, set $\psi_0^\circ(A):=6$ and, for each integer $k\geq1$, put
\[
\psi_k^\circ(A):=\sum_{z\in A}\left(z-\frac{\psi_1(A)}6\right)^k,
\]
where the sum counts multiplicities.

\begin{proof}
\emph{\textup{(a)}.} For $1\leq i\leq3$, let
\[
  h'_i=h_i-\frac{h_1+h_2+h_3}{3}, \qquad a_{i,1}=\frac{{a'_{i,1}}+{a'_{i,2}}}{2}, \qquad a_{i,2}=\frac{{a'_{i,1}}-{a'_{i,2}}}{2}.
\]
Since $u\in\mathbb Q^\times$, the identity $y_2=(u^2+1)/u$ gives $y_2\ne0$.
With the subscripts of $h'_i$ read cyclically modulo $3$, substituting \eqref{eq:degree-five-data-concise} and using $y_2^2=y_1^2+4$ gives
\[
\begin{gathered}
 \sum_{i=1}^3 h'_i=0,\qquad
 a_{i,1}a_{i,2}=\frac{(y_1^2-y_1+4)(y_1^2+2y_1+4)^2}{y_2}(h'_{i+1}-h'_{i+2}),\\
 a_{i,1}^2+a_{i,2}^2+3(h'_i)^2
 =y_2^4(y_1-2)^2(y_1+1)^2(y_1+4)^2+2\sum_{j=1}^3(h'_j)^2
 \qquad(1\leq i\leq3).
\end{gathered}
\]
Apply \cref{lem:degree-five-block-criterion} with $h'_i$ in place of $h_i$, and let $\upsilon$ be the resulting replacement map.
Then
\[
\upsilon_*A'_1=A_1(u)-\frac{h_1+h_2+h_3}{3},\qquad
\upsilon_*A'_2=A_2(u)-\frac{h_1+h_2+h_3}{3}.
\]
Since a common translation preserves equal power sums, $\psi_k(A_1(u))=\psi_k(A_2(u))$ for $1\leq k\leq5$.
Newton's identities show that the polynomials $\prod_{z\in A_1(u)}(X-z)$ and $\prod_{z\in A_2(u)}(X-z)$ differ by a constant, so $\prod_{z\in A_1(u)}z-\prod_{z\in A_2(u)}z$ is invariant under a common translation.
Direct factorization gives
\begin{align}
  \prod_{z\in A_1(u)}z-\prod_{z\in A_2(u)}z
 ={}&128y_1^2(y_1-2)^3(y_1+1)^2(y_1+2)(y_1+4)^2
  (y_1^2+2)(y_1^2+8)\notag\\
  &\times(y_1^2-y_1+4)(y_1^2+2y_1+4)^2(2y_1^2+y_1+8)
  (y_1^3+2y_1^2+4y_1+12)\notag\\
  &\times(3y_1^3+4y_1^2+8y_1+16).
\label{eq:degree-five-product-concise}
\end{align}
The factors $y_1^2+2$, $y_1^2+8$, $y_1^2-y_1+4$, $y_1^2+2y_1+4$, and $2y_1^2+y_1+8$ are positive.
The polynomials $X^3+2X^2+4X+12$ and $3X^3+4X^2+8X+16$ have no zeros modulo $5$.
Their leading coefficients are prime to $5$, so any rational root written in lowest terms has a denominator prime to $5$ by the rational-root theorem.
Reducing such a root modulo $5$ would give a contradiction.
A zero of one of $y_1$, $y_1-2$, $y_1+1$, $y_1+2$, and $y_1+4$ requires $y_1\in\{0,2,-1,-2,-4\}$.
For $y_1\in\{2,-1,-2,-4\}$, we have $y_1^2+4\in\{5,8,20\}$, none of which is a square in $\mathbb Q$.
Since $y_2\in\mathbb Q$ and $y_2^2=y_1^2+4$, these four values of $y_1$ are impossible. Thus the right-hand side of \eqref{eq:degree-five-product-concise} can vanish only when $y_1=0$, hence $u=\pm1$.

For $u\ne\pm1$,
\[
\prod_{z\in A_1(u)}z\ne\prod_{z\in A_2(u)}z,
\]
so \cref{thm:ideal-block-transfer}(b) shows that $A_1(u)\ne A_2(u)$, and \cref{thm:ideal-block-transfer}(c) shows that their supports are disjoint.
Hence \([A_1(u)]={}_5^6[A_2(u)]\).

Direct expansion gives
\[
  \psi_3^\circ(A_1(u)) = -96y_1y_2(y_1^2+4)(y_1^2+2y_1+4)^2(y_1^3+6y_1+16)(2y_1^3+3y_1^2+8).
\]
The polynomials $X^3+6X+16$ and $2X^3+3X^2+8$ have no zeros modulo $5$.
Their leading coefficients are prime to $5$, so any rational root written in lowest terms has a denominator prime to $5$ by the rational-root theorem.
Reducing such a root modulo $5$ would give a contradiction.
Since $y_1^2+4>0$, $(y_1^2+2y_1+4)=(y_1+1)^2+3>0$, and $y_2\neq0$, $\psi_3^\circ(A_1(u))=0$ is possible only for $u=\pm1$.
For $u\ne\pm1$ and any invertible rational affine map ${\alpha}$, we have
\[
\psi_3^\circ({\alpha_{*}}A_1(u))={\bigl(\alpha(1)-\alpha(0)\bigr)}^3\psi_3^\circ(A_1(u))\neq0.
\]
A six-element multiset $A$ satisfying $-A=A$ has $\psi_3^\circ(A)=0$, so no equivalent solution can have both blocks invariant under negation.
Thus the solution is nonsymmetric.

\emph{\textup{(b)}.}
For $u\neq\pm1$, the sum of real squares $\psi_2^\circ(A_1(u))$ is positive: if it were zero, all centered entries would vanish, contrary to $\psi_3^\circ(A_1(u))\neq0$.
For a six-element multiset \(A\) satisfying \(\psi_2^\circ(A)\neq0\), put
\[
 \iota(A):=\frac{\psi_3^\circ(A)^2}{\psi_2^\circ(A)^3}.
\]
For every invertible rational affine map ${\alpha}$, the identities
\[
 \psi_k^\circ({\alpha_{*}}A)={\bigl(\alpha(1)-\alpha(0)\bigr)}^k\psi_k^\circ(A)\quad(k=2,3)
\]
give $\iota({\alpha_{*}}A)=\iota(A)$.
Since $A_1(u)$ and $A_2(u)$ have the same power sums through degree $5$, they have the same value of $\iota$; thus $\iota$ is also unchanged by interchanging the blocks.
Exact substitution gives $\iota(A_1(2))\equiv1$ and $\iota(A_1(-3))\equiv4\pmod{11}$, with denominators not divisible by $11$.
Hence $\iota(A_1(u))$ is a nonconstant rational function of $u$.
After substituting $y_1=u-u^{-1}$ and $y_2=u+u^{-1}$, each entry of $A_1(u)$ and $A_2(u)$ is an integer Laurent polynomial in $u$ with exponents between $-5$ and $5$.
For the integer blocks $A_1,A_2$ in \cref{ex:degree-five-ideal-solution}, direct evaluation at $u=-3$ gives
\begin{equation}\label{eq:degree-five-test-point-affine-map}
\left( -\frac{243}{32}A_1(-3) + 845, -\frac{243}{32}A_2(-3) + 845 \right) = (A_1,A_2).
\end{equation}
The twelve integers in \eqref{eq:degree-five-test-point-concise} are distinct, so the twelve entries at $u=-3$ are also distinct.
Thus every pairwise difference is a nonzero rational function of $u$.
Each such difference has only finitely many zeros, so only finitely many parameters make two entries equal.
The preimage of each value under the nonconstant rational function $u\mapsto \iota(A_1(u))$ is finite, since clearing denominators gives a nonzero polynomial equation.
After excluding the finitely many parameters at which two entries coincide, choose recursively a sequence $(u_i)_{i\geq1}$ of distinct integers in $\mathbb Z\setminus\{0,\pm1\}$ for which the values $\iota(A_1(u_i))$ are pairwise distinct.
For every $i\geq1$, $u_i^5A_1(u_i)$ and $u_i^5A_2(u_i)$ are sets of integers.
Since $\iota$ is affine invariant, these solutions are pairwise affinely inequivalent over $\mathbb{Q}$.
\end{proof}

\subsection{Proof of Proposition~\ref{prop:degree-five-separation}}\label{app:degree-five-separation}
\begin{proof}
We use the six-element multisets $A_1(u),A_2(u)$ in \eqref{eq:degree-five-family-concise} and the integer blocks $A_1,A_2$ in \eqref{eq:degree-five-test-point-concise}.
We use the centered power sums $\psi_k^\circ$ defined above.
Substitution into the first block $A$ of \cite[Eq.~(3.8)]{Choudhry2003} gives \(4\psi_4^\circ(A)=(\psi_2^\circ(A))^2\).
The second block has the same mean and power sums through degree four, so it satisfies the same identity.
For every invertible rational affine map ${\alpha}$, the replacement $A\mapsto{\alpha_{*}}A$ multiplies $4\psi_4^\circ(A)-(\psi_2^\circ(A))^2$ by ${\bigl(\alpha(1)-\alpha(0)\bigr)}^4$.
Thus the identity holds for both blocks of every solution affinely equivalent to this family, including after interchange of the blocks.
For \(A_1\) in \eqref{eq:degree-five-test-point-concise}, we have \(4\psi_4^\circ(A_1)-(\psi_2^\circ(A_1))^2\equiv1\pmod5\).
The affine map in \eqref{eq:degree-five-test-point-affine-map} therefore shows that $4\psi_4^\circ(A_1(u))-(\psi_2^\circ(A_1(u)))^2$ is nonzero at $u=-3$.
It is therefore a nonzero rational function of $u$ and has only finitely many zeros in $\mathbb Q^\times$.
Equivalence to the specified family can occur only at these zeros.

Given a two-block solution $(A'_1,A'_2)$ with blocks of size six, partition the six occurrences in each block into three unordered pairs, match the pairs between the blocks, and order each pair.
Index the pairs of $A'_1$ by $i=1,2,3$, and write the $i$-th ordered pair and the corresponding ordered pair of $A'_2$ as \(({s_{1,i}},{s'_{1,i}})\) and \(({s_{2,i}},{s'_{2,i}})\), respectively.
Put
\[
 o:=\sum_{i=1}^3
 \bigl(({s_{1,i}}-{s'_{1,i}})+({s_{2,i}}-{s'_{2,i}})\bigr)
 \bigl(({s_{1,i}}-{s'_{1,i}})-({s_{2,i}}-{s'_{2,i}})\bigr)^3.
\]
Take the twelve rational entries in \cite[Eq.~(3.30)]{Choudhry2017}, pairing positions $i$ and $i+3$ within each block for $1\leq i\leq3$.
Denote the two progression centers in the cited construction by $h_{i,1},h_{i,2}$ and its half-length parameters by $w_i$.
These parameters satisfy the relations of \cite[Eqs.~(3.27)--(3.29)]{Choudhry2017}.
The paired entries satisfy
\[
{s_{1,i}}+{s'_{1,i}}={s_{2,i}}+{s'_{2,i}}, \qquad ({s_{1,i}}-{s'_{1,i}})+({s_{2,i}}-{s'_{2,i}})=2(h_{i,1}-h_{i,2}) \qquad(1\leq i\leq3).
\]
The differences $({s_{1,i}}-{s'_{1,i}})-({s_{2,i}}-{s'_{2,i}})$ are the same multiple of $w_i$ for $1\leq i\leq3$.
Dividing the identities in \cite[Eqs.~(3.14) and (3.18)]{Choudhry2017} by $2$ gives
\[
 \sum_{i=1}^3w_i(h_{i,1}-h_{i,2})
 =\sum_{i=1}^3w_i(4w_i^2-1)(h_{i,1}-h_{i,2})=0.
\]
Expanding the second sum and using the first gives $\sum_{i=1}^3w_i^3(h_{i,1}-h_{i,2})=0$; the paired-entry identities then give $o=0$.
The equalities \({s_{1,i}}+{s'_{1,i}}={s_{2,i}}+{s'_{2,i}}\) for \(1\leq i\leq3\) are preserved by any common invertible rational affine map \(\alpha\). Such a map multiplies \(o\) by \((\alpha(1)-\alpha(0))^4\); interchanging the blocks changes the sign of \(o\).
Consequently, every solution equivalent over $\mathbb Q$ to one obtained from \cite[Eq.~(3.30)]{Choudhry2017} with twelve rational entries has partitions of its two blocks, a bijection between the two sets of pairs, and orderings satisfying
\[
 {s_{1,i}}+{s'_{1,i}}={s_{2,i}}+{s'_{2,i}}\quad(1\leq i\leq3),\qquad o=0.
\]

For the integer blocks $A_1,A_2$ in \eqref{eq:degree-five-test-point-concise}, the common sums of unordered pairs are $-2560$, $870$, and $1690$, each occurring once in each block.
The three corresponding pairs partition each block, so the matched partitions satisfying ${s_{1,i}}+{s'_{1,i}}={s_{2,i}}+{s'_{2,i}}$ for $1\leq i\leq3$ are unique up to reordering the pairs.
For their $2^6$ possible pair orderings, direct calculation gives $o\bmod29\in\{3,14,18,21,25,28\}$, so $o\ne0$ for every ordering.
By the affine map in \eqref{eq:degree-five-test-point-affine-map}, no choice of partitions, matching, and orderings at $u=-3$ satisfies ${s_{1,i}}+{s'_{1,i}}={s_{2,i}}+{s'_{2,i}}$ for $1\leq i\leq3$ and $o=0$.
For the blocks $A_1(u),A_2(u)$ in \eqref{eq:degree-five-family-concise}, fix partitions, a bijection between the sets of pairs, and orderings, and denote the resulting paired occurrences by $({s_{1,i}},{s'_{1,i}})$ and $({s_{2,i}},{s'_{2,i}})$.
By the formulas of Section~\ref{subsec:pte-block-transfer}, the three differences ${s_{1,i}}+{s'_{1,i}}-{s_{2,i}}-{s'_{2,i}}$ and $o$ are evaluations at $X=u$ of four elements of $\mathbb Q[X,X^{-1}]$.
For each fixed choice, at least one of these Laurent polynomials is nonzero at $X=-3$.
That polynomial has finitely many roots in $\mathbb Q^\times$, so only finitely many parameters make all four evaluations vanish.
Since there are finitely many pairs of partitions, bijections, and orderings, only finitely many parameters yield equivalence over $\mathbb Q$ to a solution from \cite[Eq.~(3.30)]{Choudhry2017} with twelve rational entries.
Combining the finite parameter sets obtained for \cite[Eq.~(3.8)]{Choudhry2003} and \cite[Eq.~(3.30)]{Choudhry2017} proves the proposition.
\end{proof}

\subsection*{Acknowledgements}

F.~\v{S}piranec thanks Ziqing Xiang for reading an earlier draft and for noting that the use of the sunflower lemma is a conceptual innovation in addressing the problem of rational spherical $4$-designs on the unit circle posed by Cui--Xia--Xiang~\cite[Section~1.3]{CuiXiaXiang2019}.

This work was partially supported by JSPS KAKENHI Grant Numbers 26K00617 and 26KJ0570.

R.~Misawa and F.~\v{S}piranec used ChatGPT (OpenAI, GPT-5.6 Pro) to check computations and to assist with the organization and presentation of the manuscript.
In particular, the AI produced a proof that the size of every rational spherical $4$-design on the unit circle is divisible by $2304$ (\cref{prop:fourth-moment-divisibility}), and R.~Misawa verified the proof.
The authors take full responsibility for the content of this paper.


\end{document}